\documentclass[pdflatex,sn-mathphys-num]{sn-jnl}

\usepackage{graphicx}%
\usepackage{multirow}%
\usepackage{amsmath,amssymb,amsfonts}%
\usepackage{amsthm}%
\usepackage{mathrsfs}%
\usepackage[title]{appendix}%
\usepackage{xcolor}%
\usepackage{textcomp}%
\usepackage{manyfoot}%
\usepackage{booktabs}%
\usepackage{algorithm}%
\usepackage{algorithmicx}%
\usepackage{algpseudocode}%
\usepackage{listings}%
\usepackage{tikz}
\usepackage{tikz-cd}
\usepackage{amscd}
\usepackage[all,cmtip]{xy}

\theoremstyle{thmstyleone}%
\newtheorem{thm}{Theorem}
\newtheorem{prop}[thm]{Proposition}%
\newtheorem{cor}[thm]{Corollary}
\theoremstyle{thmstyletwo}%

\theoremstyle{thmstylethree}%
\newtheorem{defn}{Definition}%
\newtheorem{lem}{Lemma}
\newcommand\A{\mathbb{A}}

\newcommand\Hom{Hom}
\newcommand\Ind{Ind}

\newcommand\supp{supp}
\newcommand\ind{ind}
\newcommand\Tr{Tr}

\newcommand \ONB{ONB}

\newcommand\GL{GL}

\newcommand\OO{O}
\newcommand\PGL{PGL}

\newcommand\SO{SO}
\newcommand\GSO{GSO}

\newcommand\Sp{Sp}

\newcommand\GO{GO}
\newcommand\GSp{GSp}
\newcommand\PGSp{PGSp}
\newcommand\PGSO{PGSO}
\newcommand\PGO{PGO}
\newcommand\U{U}
\newcommand\SU{SU}

\newcommand \dd{d}
\newcommand \rest{rest}
\newcommand \temp{temp}

\newcommand{\normmm}[1]{{\left\vert\kern-0.25ex\left\vert\kern-0.25ex\left\vert #1 
		\right\vert\kern-0.25ex\right\vert\kern-0.25ex\right\vert}}

\begin{document}

\title[Article Title]{The Sakellaridis-Venkatesh Conjecture for $  U(2)\backslash SO_{2,3}$}


\author{\fnm{} \sur{Wan Xiaolei}}

\affil[]{\orgdiv{Beijing International Center for Mathematical Research}, \orgname{Peking University}, \orgaddress{\street{Yiheyuan Road}, \city{Beijing}, \postcode{100871}, \country{China}}}


\abstract{This article gives the Plancherel decomposition of $L^2\left(SO_{2,3}(F)\slash U(2)(F)\right)$, where $F$ is a local field with characteristic $0$. Finally, we obtain a factorization of the global period $P_{U(2)}$.}




\maketitle

\section{Introduction}\label{sec1}

\subsection{\bf Strongly tempered varieties}
Let $B=MN$ be a Borel subgroup of $G$ and $K$ a maximal compact subgroup of $G$ which is special when $F$ is $p$-adic. Let $ dk$ be the probability measure on $K$. Choose maps 
$$
m_{B}:G\longrightarrow M;\qquad n_{B}:G\longrightarrow N;\qquad k_{B}:G\longrightarrow K
$$
such that $g=m_{B}(g)\cdot n_{B}(g)\cdot k_{B}(g)$ for all $g \in G$. Pick the element $e$ in the normalized induction representation $ ind_{B}^{G}1$ such that $e(k)=1$ for all $k \in K$. 


\begin{defn}
	The Harish-Chandra Xi function is defined to be
	$$
	\Xi^{G}(g):=\int_{K} R(g) e(k)\cdot\overline{e(k)} dk.
	$$
\end{defn}
Then we have the equality:
$$
\Xi^{G}(g)=\int_{K}\delta_{B}(m_{B}(kg))^{1/2}  dk,
$$
where $\delta_{B}$ is the modular character of $B$.

One can verify $\Xi^{G}(g^{-1})=\Xi^{G}(g)$, $\Xi^{G}(k_1gk_2)=\Xi^{G}(g)$ for $k_1, k_2\in K$, and $\Xi^{G}(g)>0$. 
For $f_1,f_2 \in C(G)$, we say that $f_1$ and $f_2$ are equivalent if there exist $c>0$ and $C>0$ such that 
$$
c\cdot f_2 \leqslant f_1 \leqslant C\cdot f_2.
$$ 
The function $\Xi^{G}$ depends on the choices of $B$, $K$ and the measure $dk$ we made, but it does not matter because different choices would yield equivalent functions. Using the Harish-Chandra Xi function, one can define the Harish-Chandra Schwartz space $\mathcal{C}(G)$ of $G$, and the weak Harish-Chandra Schwartz space $\mathcal{C}^{\omega}(G)$ of $G$. For example, one can find the definitions in \cite[\S2.4]{BP1}. It is the space of smooth vectors in the space of tempered functions. One has 
$$
C_c^{\infty}(G)\subset \mathcal{C}(G)\subset \mathcal{C}^{\omega}(G)
$$
and $C_c^{\infty}(G)$ is dense in $ \mathcal{C}^{\omega}(G)$. All tempered matrix coefficients of $G$ lie in the weak Harish-Chandra Schwartz space $\mathcal{C}^{\omega}(G)$.

Recall the definition of strongly tempered varieties. 
\begin{defn}
	We say that a spherical variety $X$ is strongly tempered if, for any $x \in X(F)$ with stabilizer $H$, the restriction of any irreducible tempered matrix coefficient to $H$ is in $L^{1}(H)$, i.e.\   for every irreducible tempered unitary representation $\pi$ of $G$ and $v,w \in \pi^{\infty}$, one has
	$$
	\int_{H}\left|\langle h\cdot v,w \rangle_{\pi}\right|  dh < \infty.
	$$

\end{defn}

For example, $SO(2)\backslash SO(3)$, and $GL_n\backslash GL_n\times GL_{n+1}$  are strongly tempered.

\subsection{\bf Whittaker varieties}
Let $B=MN$ be a Borel subgroup of $G$, $\psi$ a non-degenerate unitary character of $N$. By \cite[\S 7.1]{BP1}, 
the functional 
$$
C_c^{\infty}( G) \longrightarrow \mathbb{C}_{\psi}
$$
given by
$$
f\longmapsto \int_{N}f(n)\cdot \overline{\psi(n)} d n
$$
extends continuously to a (necessarily unique) functional 
$$
\mathcal{C}^{\omega}(G)\longrightarrow \mathbb{C}_{\psi}
$$
to be denoted by
$$
f\longmapsto \int_{N}^{\ast}f(n)\cdot \overline{\psi(n)} d n.
$$

\subsection{\bf Bernstein's method}

Let $G$ be a reductive group over a local field $F$ acting transitively on a homogeneous space $X$ of polynomial growth (for definition of polynomial growth, see \cite[\S 0 Pg.665]{B}). We fix a base point $x_0 \in X$ with stabilizer $H \subset G$, so that $g\mapsto g^{-1}\cdot x_0$ gives an identification $H\backslash G\cong X$. Note that the the space of $F$-points of a spherical variety is of polynomial growth (see \cite[\S11.2 Pg.155]{SV}).

Fix a $G$-invariant measure $d x$ on $X$ and a $G$-invariant inner product in $L^2(X)$
\[\langle \phi_1 , \phi_2 \rangle_X := \int_X \phi_1(x)\cdot \overline{\phi_2(x)}d x,\]
where $\phi_1,\phi_2 \in L^2(X)$. Then $L^2(X)$ is a unitary representation of $G$ by right translation. Such a unitary representation admits a direct integral decomposition
\begin{equation} 
	\iota:L^2(X) \cong \int_{\Omega}\sigma(\omega)d \mu_{\Omega}(\omega)  \label{bern:eq7},
\end{equation}
where $\Omega$ is a Borel measurable space with a measure $d \mu_{\Omega}$ and $\sigma: \Omega \to \widehat{G}$ is a measurable map from $\Omega$ to the unitary dual $\widehat{G}$ of $G$ which is equipped with a standard Borel structure on $\widehat{G}$ \cite[Proposition 4.6.1]{Dix}. 
In \cite[Theorem 1.5]{B}, Bernstein showed that the natural embedding $C_c^{\infty}(X)\hookrightarrow L^2(X)\cong \int_{\Omega}\sigma(\omega)d\mu_{\Omega}(\omega)$ is pointwise-defined on $C_c^{\infty}(X)$, i.e.\ there is a family of $G$-equivariant continuous maps
\[\alpha_{\omega}^{X}: C_c^{\infty}(X)  \longrightarrow \sigma(\omega)\]
defined for $\mu_{\Omega}$-almost all $\omega$ such that for each $\phi \in C_c^{\infty}(X)$, the function
\[\omega \mapsto \alpha^{X}_{\omega}(\phi)\]
represents the vector $\iota(\phi) \in \int_{\Omega}\sigma(\omega)d \mu_{\Omega}(\omega)$. Moreover, in \cite[\S3.5 Pg.689]{B}, Bernstein proved that the natural embedding $\mathcal{C}(X)\hookrightarrow L^2(X)\cong \int_{\Omega}\sigma(\omega)d\mu_{\Omega}(\omega)$ is also pointwise-defined on $\mathcal{C}(X)$, where the $LF$ space  $\mathcal{C}(X) \supset C_c^{\infty}(X)$ is called the Harish-Chandra Schwartz space of $X$ (for definition, see \cite[\S3.5 Pg.689]{B}). In other words, the family of $G$-equivariant maps $\alpha_{\omega}^{X}$ can extend continuously to $\mathcal{C}(X)$
\[\alpha_{\omega}^{X}: \mathcal{C}(X)  \longrightarrow \sigma(\omega)\]
such that for each $\phi \in \mathcal{C}(X)$, the function
\[\omega \mapsto \alpha^{X}_{\omega}(\phi)\]
represents the vector $\iota(\phi) \in \int_{\Omega}\sigma(\omega)d \mu_{\Omega}(\omega)$.

When $\iota$ in equation (\ref{bern:eq7}) is fixed, such a family $\{\alpha^{X}_{\omega}\mid \omega \in supp(\mu_{\Omega})\}$ is essentially unique, i.e.\ any two such families differ only on a subset of $\Omega$ with measure $0$ with respect to $\mu_{\Omega}$, and the maps $\alpha^{X}_{\omega}$ are nonzero for $\mu_{\Omega}$-almost all $\omega\in supp(\mu_{\Omega})$. If $\alpha_{\omega}^{X}$ is nonzero, then by duality, one has a $G$-equivariant embedding
$$
\overline{\beta}^{X}_{\omega}:\sigma(\omega)^{\vee}\cong \overline{\sigma(\omega)}\longrightarrow C^{\infty}(X),
$$
where the isomorphism $\sigma(\omega)^{\vee}\cong \overline{\sigma(\omega)}$ is induced by $\langle \cdot,\cdot \rangle_{\sigma}$ and the duality between $C_c^{\infty}(X)$ and $C^{\infty}(X)$ is given by the integration with respect to the measure $d x$. Take complex conjugate on both sides, one obtains a family of $G$-equivariant maps 
$$
\beta^{X}_{\omega}:\sigma(\omega)^{\infty}\longrightarrow C^{\infty}(X),
$$
such that
\begin{equation}
	\langle\alpha_{\omega}^{X}(f), v \rangle_{\sigma(\omega)}=\langle f,\beta_{\omega}^{X}(v)  \rangle_{X} \label{bern:eq4},
\end{equation}
for $f \in C_c^{\infty}(X)$ and $v \in \sigma(\omega)^{\infty}$.

When $F$ is $p$-adic, we endow $\sigma(\omega)^{\infty}$ with the finest locally convex topology. When $F=\mathbb{R}$, we endow $\sigma(\omega)^{\infty}$ with the topology defined by the semi-norms
$$
\|v\|_{u}:=\|\pi^{\infty}(u)v\|,
$$
where $u \in \mathcal{U}(\mathfrak{g})$ and $\|\cdot\|$ is the norm induced by the inner product on $\sigma(\omega)$. If we compose $\beta_{\omega}^{X}$ with the evaluation-at-$x_0$ map $ev_{x_0}$, one obtains a family of continuous $H$-invariant maps
\[\ell_{\omega}^{X}:=ev_{x_0}\circ \beta_{\omega}^{X} \in Hom_{H}\left(\sigma(\omega)^{\infty},\mathbb{C}\right). \]

These families depend on the isomorphism $\iota$ in the equation \eqref{bern:eq7}, and are defined for $\mu_{\Omega}$-almost all $\omega\in \Omega$. Changing $\iota$ will result in another families which differ from the original one by a measurable function $f: supp(\mu_{\Omega})\longrightarrow \mathbb{S}^1$. Thus the family
\[\{\alpha_{\omega}^{X}\otimes\overline{\alpha}_{\omega}^{X}:\omega \in supp(\mu_{\Omega})\} \]
is independent of the choice of the isomorphism $\iota$ in \eqref{bern:eq7}. Likewise, the family
\[\{\beta_{\omega}^{X}\otimes\overline{\beta}_{\omega}^{X}:\omega \in supp(\mu_{\Omega})\} \] 
and
\[\{\ell_{\omega}^{X}\otimes\overline{\ell}_{\omega}^{X}:\omega \in supp(\mu_{\Omega})\} \] 
are independent of $\iota$.

Define the positive semi-definite Hermitian form $J^{X}_{w}(\cdot,\cdot)$ of $C_c^{\infty}(X)$ by
$$
J^{X}_{w}(f_1,f_2):=\langle \alpha_{\omega}^{X}(f_1),\alpha_{\omega}^{X}(f_2) \rangle_{\sigma(\omega)},
$$
where $f_1,f_2 \in C_c^{\infty}(X)$. Then one knows $J_{\omega}^{X}$ factors through $\sigma(\omega)^{\infty}\otimes \overline{\sigma(\omega)^{\infty}}$ and
\begin{equation} \label{pointdecomp: eq}
	\langle f_1,f_2 \rangle_{X}=\int_{\Omega}J^{X}_{w}(f_1,f_2) d\mu_{\Omega}(\omega).
\end{equation}

By a theorem of Dixmier-Malliavin \cite{Di}, one obtains
\begin{lem} \label{zf12: pointwiselemma}
	\begin{enumerate}
		\item For any $G$-spherical variety $X$, one has
		$$
		L(C_c^{\infty}(G))\mathcal{C}(X)=\mathcal{C}(X),
		$$
		where $L$ is the left translation representation.
		\item For any irreducible tempered representation $\sigma$, one has
		$$
		\sigma^{\infty}\left(C_c^{\infty}(G)\right) \sigma^{\infty}=\sigma^{\infty}.
		$$
	\end{enumerate}
\end{lem}

We have a more precise decomposition than Equation (\ref{pointdecomp: eq}).
\begin{lem} \label{zf12: pointwiselemma2}
	For any $f \in \mathcal{C}(X)$ and any $x \in X$, one has
	$$
	f(x)=\int_{\Omega}\beta^{X}_{\omega}\circ \alpha_{\omega}^{X}(f)(x) d\mu_{\Omega}(\omega).
	$$
\end{lem}
\begin{proof}
	If the equation
	$$
	f(x_0)=\int_{\Omega}\beta^{X}_{\omega}\circ \alpha_{\omega}^{X}(f)(x_0) d\mu_{\Omega}(\omega)
	$$
	holds for any $f \in \mathcal{C}(X)$, then we can replace $f$ by $L(g)f$ to get
	$$
	f(g^{-1}\cdot x_0)=\int_{\Omega}\beta^{X}_{\omega}\circ \alpha_{\omega}^{X}(f)(g^{-1}\cdot x_0) d\mu_{\Omega}(\omega).
	$$
	
	Thus we need only to prove the equation
	$$
	f(x_0)=\int_{\Omega}\beta^{X}_{\omega}\circ \alpha_{\omega}^{X}(f)(x_0) d\mu_{\Omega}(\omega)
	$$
	holds for any $f \in \mathcal{C}(X)$. By Lemma \ref{zf12: pointwiselemma}, we need only to check that for any $\varphi \in C_c^{\infty}(G)$ and any $f \in \mathcal{C}(X)$, we have
	\begin{align*}
		L(\varphi)f(x_0)=&\int_{\Omega}\beta^{X}_{\omega}\circ \alpha_{\omega}^{X}(L(\varphi)f)(x_0) d\mu_{\Omega}(\omega)  \\
		=&\int_{\Omega}L(\varphi)\beta^{X}_{\omega}\circ \alpha_{\omega}^{X}(f)(x_0) d\mu_{\Omega}(\omega).
	\end{align*}
	
	Let $\varphi^{\prime}(x):=\int_{H}\varphi(h\cdot x)d h \in C_c^{\infty}(X)$. One has
	\begin{align*}
		L(\varphi)f(x_0)=&\int_{X}\varphi^{\prime}(x)\cdot f(x)d x \\
		=& \int_{\Omega}\langle \alpha_{\omega}^{X}(\varphi^{\prime}), \alpha_{\omega}^{X}(f) \rangle_{\sigma(\omega)}d \mu_{\Omega}(\omega) \\
		=&\int_{\Omega}\langle \varphi^{\prime}, \beta_{\omega}^{X}\circ\alpha_{\omega}^{X}(f) \rangle_{X}d \mu_{\Omega}(\omega) \\
		=&\int_{\Omega}L(\varphi)\beta^{X}_{\omega}\circ \alpha_{\omega}^{X}(f)(x_0) d\mu_{\Omega}(\omega).
	\end{align*}
\end{proof}

\subsection{\bf Relative character} \label{subsection: relativecharacter}
In \cite[\S7]{GW}, one sees various incarnations of relative characters. We choose one of them: define the relative character $\vartheta^{X}_{\omega}$ of $C_c^{\infty}(X)$ associated to $\omega \in \Omega$ by
$$
\vartheta^{X}_{\omega}(f):=\ell_{\omega}^{X}\circ \alpha_{\omega}^{X}(f)
$$  
for $f \in C_c^{\infty}(X)$. Note that although $\alpha_{\omega}^{X}$ and $\ell^{X}_{\omega}$ is dependent of the isomorphism in equation (\ref{bern:eq7}), the relative character $\vartheta^{X}_{\omega}$ is not. 

\begin{lem}
	For almost all $\omega$, the relative character $\vartheta^{X}_{\omega}$ can extend to $\mathcal{C}(X)$ continuously.
\end{lem}
This relative character depends on $\omega \in \Omega$ and $f \in \mathcal{C}(X)$.

\subsection{\bf Plancherel decompositions}
Let $X=H\backslash G$ be a spherical variety. Motivated by and refining the work of Gaitsgory-Nadler in the geometric Langlands program, Sakellaridis and Venkatesh \cite[\S3]{SV} associated two data to $X$: a dual group $X^{\vee}$ and a distinguished morphism
\begin{equation}  
	\kappa: X^{\vee} \times \mathrm{SL}_{2}(\mathbb{C}) \longrightarrow    G^{\vee}. 
\end{equation}
In some favorable cases, $\kappa$ gives rise to a map
$$
\kappa_{\ast}: \widehat{G}_{X}^{temp}\longrightarrow \widehat{G}
$$
from the tempered unitary dual of $G_{X}$ to the unitary dual of $G$,where $G_{X}$ is a reductive group defined over $F$ such that $G_{X}^{\vee}=X^{\vee}$. Under these assumptions, one hopes to have a spectral decomposition 
\begin{equation} \label{zf12: planchereleq1}
	L^2\left( X \right)\cong \int_{\widehat{G}_{X}^{temp}}m(\sigma)\cdot \kappa_{\ast}(\sigma) d \mu_{G_{X}}(\sigma),
\end{equation}
where $\mu_{G_X}$ is the Plancherel measure of $G_{X}$, and $m(\sigma)$ is some multiplicity space.

If $\dim m(\sigma)\leqslant 1$ for all $\sigma\in \widehat{G}_{X}^{temp}$, then in the isomorphism (\ref{bern:eq7}), we have 
$$
\Omega=\{\sigma \in \widehat{G}_{X}^{temp}\mid \dim m(\sigma)=1 \}\subset \widehat{G}_{X}^{temp}
$$
and $\mu_{\Omega}$ is the restriction of $\mu_{G_X}$ on $\Omega$.
Thus the decomposition (\ref{zf12: planchereleq1}) can also be written as
$$
L^2(X)\cong \int_{\Omega}\kappa_{\ast}(\sigma(\omega)) d \mu_{\Omega}(\omega).
$$
If so, one would like to have explicit formulae of the positive semi-definite Hermitian forms  $J^{X}_{\omega}(\cdot,\cdot)$.

	\subsection{\bf Whittaker-torus periods}
	Consider the homogeneous $PGSO_{2,2}\cong PGL_2\times PGL_2$-space 
	$$
	X:=\left( N\times \overline{T}_{E}\backslash PGSO_{2,2}, \psi\boxtimes \mathbb{C} \right)\cong \left(N\backslash PGL_2,\psi\right)\times \left(\overline{T}_{E}\backslash PGL_2,\mathbb{C}\right).
	$$
	
	Then one has
	\begin{align*}
		&L^{2}\left( N\times \overline{T}_{E}\backslash PGSO_{2,2}, \psi\boxtimes \mathbb{C}\right) \\
		\cong &L^2\left(N\backslash PGL_2,\psi\right)\widehat{\boxtimes} L^2\left(\overline{T}_{E}\backslash PGL_2\right) \\
		\cong & \int_{\widehat{PGL}_2^{temp}}\int_{\widehat{PGL}_{2}^{temp}}  Hom_{N\times\overline{T}_{E}}\left(\sigma_1\widehat{\boxtimes}\sigma_2,\psi\boxtimes\mathbb{C}   \right)\cdot \sigma_1\widehat{\boxtimes}\sigma_2  \ d\mu_{PGL_2}(\sigma_1)d\mu_{PGL_{2}}(\sigma_2).
	\end{align*}
	Comparing 
	\begin{align*}
		L^2\left(PGSO_{2,2} \right)\cong &\int_{\widehat{PGSO}_{2,2}^{temp}}(\sigma_1\widehat{\boxtimes} \sigma_2)\widehat{\boxtimes} (\sigma_1^{\vee}\widehat{\boxtimes} \sigma_2^{\vee})d \mu_{PGSO_{2,2}}(\sigma_1\widehat{\boxtimes} \sigma_2) 
	\end{align*}
	with
	\begin{align*}
		L^2\left(PGSO_{2,2} \right) \cong & L^2\left(PGL_2\right)\widehat{\boxtimes}L^2\left(PGL_2\right) \\
		\cong & \int_{\widehat{PGL}_2^{temp}}\int_{\widehat{PGL}_{2}^{temp}} 
		(\sigma_1\widehat{\boxtimes}\sigma_1^{\vee})\widehat{\boxtimes}(\sigma_2\widehat{\boxtimes}\sigma_2^{\vee})
		d\mu_{PGL_2}(\sigma_1)d\mu_{PGL_{2}}(\sigma_2) \\
		\cong &\int_{\widehat{PGL}_2^{temp}}\int_{\widehat{PGL}_{2}^{temp}}
		(\sigma_1\widehat{\boxtimes} \sigma_2)\widehat{\boxtimes} (\sigma_1^{\vee}\widehat{\boxtimes} \sigma_2^{\vee})
		d\mu_{PGL_2}(\sigma_1)d\mu_{PGL_{2}}(\sigma_2),
	\end{align*}
	we have
	\begin{align*}
		&L^{2}\left( N\times \overline{T}_{E}\backslash PGSO_{2,2}, \psi\boxtimes \mathbb{C}\right) \\
		\cong & \int_{\widehat{PGL}_2^{temp}}\int_{\widehat{PGL}_{2}^{temp}}  Hom_{N\times\overline{T}_{E}}\left(\sigma_1\widehat{\boxtimes}\sigma_2,\psi\boxtimes\mathbb{C}   \right)\cdot \sigma_1\widehat{\boxtimes}\sigma_2  d\mu_{PGL_2}(\sigma_1)d\mu_{PGL_{2}}(\sigma_2) \\
		\cong &
		\int_{\widehat{PGSO}_{2,2}^{temp}}
		Hom_{N\times\overline{T}_{E}}\left(\sigma,\psi\boxtimes\mathbb{C}   \right)\cdot \sigma
		d \mu_{PGSO_{2,2}}(\sigma),
	\end{align*}
	where $\mu_{PGSO_{2,2}}$ is the Plancherel measure of $PGSO_{2,2}$ and $\widehat{PGSO}_{2,2}^{temp}$ is the tempered unitary dual of $PGSO_{2,2}$.

	In this case, for any $\sigma\cong \sigma_1\widehat{\boxtimes}\sigma_2\in \widehat{PGSO}_{2,2}^{temp}\cong \widehat{PGL}_2^{temp}\boxtimes\widehat{PGL}_2^{temp}$, we have
	$$
	\alpha_{\sigma}^{X}=\alpha_{\sigma_1}^{N,\psi\backslash PGL_2}\boxtimes \alpha_{\sigma_2}^{\overline{T}_{E}\backslash PGL_2};\qquad \ell_{\sigma}^{X}=\ell_{\sigma_1}^{N,\psi\backslash PGL_2}\boxtimes\ell_{\sigma_2}^{\overline{T}_{E}\backslash PGL_2}.
	$$

	\subsection{\bf Weil representation} \label{zf13:21a}
	Let $V$ be a split quadratic space of dimension four over a local field $F$. We will fix the Witt decomposition $V=V_1\oplus V_2$ and a basis
	$$
	\boldsymbol{e}:=\{e_1,e_2,e_1^{\prime},e_2^{\prime}  \}
	$$
	for $V$ with $V_1$ the span of $e_1,e_2$, $V_2$ the span of $e_1^{\prime},e_2^{\prime}$, and $\langle e_i,e_j^{\prime} \rangle_{V}=\delta_{ij}$. With this basis $\boldsymbol{e}$, we have an identification
	$$
	GO(V)=\left\{ h=\begin{pmatrix}
		A&B\\
		C&D
	\end{pmatrix}\in GL_4(F)\middle\vert\, {^{t}h}\begin{pmatrix}
		&I_2 \\
		I_2& 
	\end{pmatrix}h=\lambda\cdot \begin{pmatrix}
		&I_2 \\
		I_2& 
	\end{pmatrix} \text{ for some $\lambda \in GL_1(F)$}  \right\}.
	$$ 
	We will write $\lambda_{V}(\cdot)$ for the similitude character of $GO(V)$. Then
	$$
	O(V)=\left\{ h\in GO(V)\mid \lambda_{V}(h)=1  \right\}.
	$$
	We will denote the stabilizer of $V_1$ in $GO(V)$ by $P$, with Levi decomposition $P=MN$. Then
	\begin{align}
		M=&\left\{ m=\begin{pmatrix}
			A&\\
			&D
		\end{pmatrix}\in GL_4(F)\middle\vert\, {^{t}m}\begin{pmatrix}
			&I_2 \\
			I_2& 
		\end{pmatrix}m=\lambda\cdot \begin{pmatrix}
			&I_2 \\
			I_2& 
		\end{pmatrix} \text{ for some $\lambda \in GL_1(F)$}  \right\} \notag \\
		=&\left\{ m(A,\lambda)=\begin{pmatrix}
			A&\\
			&\lambda\cdot {^{t}A^{-1}}
		\end{pmatrix}\middle\vert\, A \in GL_2(F)  \right\} \notag
	\end{align}
	and
	\begin{align}
		N=&\left\{ n=\begin{pmatrix}
			I_2&X\\
			&I_2
		\end{pmatrix}\in GL_4(F)\middle\vert\, {^{t}n}\begin{pmatrix}
			&I_2 \\
			I_2& 
		\end{pmatrix}n= \begin{pmatrix}
			&I_2 \\
			I_2& 
		\end{pmatrix}  \right\} \notag \\
		=&\left\{ n(x)=\begin{pmatrix}
			I_2&x\cdot J\\
			&I_2
		\end{pmatrix}\middle\vert\, x \in F \right\}, \notag
	\end{align}
	where $J:=\begin{pmatrix}
		&-1 \\
		1&
	\end{pmatrix}$.
	We denote the stabilizer of $V_1$ in $O(V)$ by $P_1$, with Levi decomposition $P_1=M_1N$. Then
	$$
	M_1=\left\{ m(A,1)=\begin{pmatrix}
		A&\\
		& {^{t}A^{-1}}
	\end{pmatrix}\middle\vert\, A \in GL_2(F)  \right\}.
	$$

	Let $W$ be a symplectic space of dimension four over $F$. We fix a basis
	$$
	\boldsymbol{f}:=\{f_1,f_2,f_3,f_4  \}
	$$
	for $W$ such that 
	$$
	(\langle f_i,f_j \rangle_{W})_{4\times 4}=\begin{pmatrix}
		&&&1 \\
		&&1& \\
		&-1&& \\
		-1&&&
	\end{pmatrix}.
	$$
	Let $W_1$ be the span of $f_1,f_2$; $W_2$ be the span of $f_3,f_4$.
	With this basis $\boldsymbol{f}$, one obtains an identification
	$$
	GSp(W)=\left\{ g=\begin{pmatrix}
		A&B\\
		C&D
	\end{pmatrix}\in GL_4(F)\middle\vert\, {^{t}g}(\langle f_i,f_j \rangle_{W})g=\lambda\cdot (\langle f_i,f_j \rangle_{W}) \text{ for some $\lambda \in GL_1(F)$}  \right\}.
	$$ 
	We will write $\lambda_{W}(\cdot)$ for the similitude character of $GSp(W)$. Then
	$$
	Sp(W)=\left\{ g\in GSp(W)\mid \lambda_{W}(g)=1  \right\}.
	$$

	Fix a non-trivial additive unitary character $\psi$ of $F$. Then $\psi$ can be seen as a non-trivial multiplicative unitary character of $N$ by
	$$
	\psi\left(n(x) \right):=\psi(x).
	$$   
	The adjoint action of $M$ on $N$ 
	$$
	m(A,\lambda)^{-1}\cdot n\left(x\right)\cdot m(A,\lambda)=n\left(\lambda\det(A)^{-1}\cdot x\right)
	$$
	induces an action of $M$ on the set of multiplicative unitary characters of $N$, and the stabilizer of $\psi$ in $M$ is
	\begin{align}
		M_{\psi}:=&\left\{ m(A,\lambda)=\begin{pmatrix}
			A&\\
			&\lambda\cdot {^{t}A^{-1}}
		\end{pmatrix}\middle\vert\, m(A,\lambda)^{-1}\cdot n\left(x\right)\cdot m(A,\lambda)=n(x)  \right\} \notag \\
		=&\left\{ m\left(A,\det(A)\right)=\begin{pmatrix}
			A&\\
			&\det(A)\cdot {^{t}A^{-1}}
		\end{pmatrix}\middle\vert\, A \in GL_2(F)  \right\}. \notag
	\end{align}

	Let 
	$$
	GL_2^{(1)}:=\left\{ h(A)=\begin{pmatrix}
		a\cdot I_2& b\cdot J\\
		c\cdot (-J)&d\cdot I_2
	\end{pmatrix}\in GL_4(F)\middle\vert\, A=\begin{pmatrix}
		a&b\\
		c&d
	\end{pmatrix}\in GL_2(F) \right\}\subset GSO(V).
	$$
	Let 
	$$
	B:=\left\{h(A)\in GL_2^{(1)}\middle\vert\, A= \begin{pmatrix}
		a&b\\
		0&d
	\end{pmatrix}\in GL_2(F) \right\}.
	$$
	Denote the center of $GSO(V)$ by $Z_V$. Then one has
	$$
	GL_2^{(1)}\bigcap M_{\psi}=Z_V
	$$
	and
	$$
	Z_V^{\bigtriangledown}\backslash GL_2^{(1)}\times M_{\psi}\cong GSO(V),
	$$
	where $Z_V^{\bigtriangledown}:=\{(\lambda\cdot I_4,\lambda^{-1}\cdot I_4)\in GL_2^{(1)}\times M_{\psi}\mid \lambda \in GL_1(F)\}$. Thus
	$$
	PGSO(V)\cong \overline{GL_2^{(1)}}\times \overline{M}_{\psi}\cong PGL_2\times PGL_2,
	$$
	where $\overline{ GL}_2^{(1)}:=Z_V\backslash GL_2^{(1)}$ and $\overline{M}_{\psi}:=Z_V\backslash M_{\psi}$. Let
	$$
	\overline{B}:=Z_{V}\backslash B \subset \overline{GL_2^{(1)}}.
	$$

	Let 
	$$
	\boldsymbol{t}:=\begin{pmatrix}
		1&&& \\
		&&&-1 \\
		&&1& \\
		&-1&&
	\end{pmatrix}.
	$$
	Then $\boldsymbol{t} \in GO(V)$ but $\notin GSO(V)$, so that $GO(V)=GSO(V)\rtimes \langle\boldsymbol{t}\rangle$. Moreover, one has
	$$
	\boldsymbol{t}^{-1}\cdot h(A)\cdot \boldsymbol{t}=m(A,\det(A))
	$$
	and
	$$
	\boldsymbol{t}^{-1}\cdot m(A,\det(A))\cdot \boldsymbol{t}=h(A).
	$$
	Thus an irreducible (unitary) representation of $GSO(V)$ is of the form $\sigma_1\boxtimes\sigma_2$ with an irreducible (unitary) representation $\sigma_i$ of $GL_2(F)$ such that the central characters of $\sigma_1$ and $\sigma_2$ are same. Moreover, the action of $\boldsymbol{t}$ sends $\sigma_1\boxtimes \sigma_2$ to $\sigma_2\boxtimes\sigma_1$.

	\subsection{\bf Weil representation}
	The Weil representation $\omega_{\psi}$ of $O(V)\times Sp(W)$ can be realized on the space $S(V_2\otimes W)$ of Bruhat-Schwartz functions and the action of $(M_1\ltimes N)\times Sp(W)$ on $S(V_2\otimes W)$ is given by
	
	\begin{equation} \notag
		\left\{
		\begin{aligned}
			&\omega_{\psi}(h)\Phi(T)=\Phi(h^{-1}T),     &\text{for $h \in Sp(W)$;}                  \\
			&\omega_{\psi}\left(m(A,1)\right)\Phi(T)= \left| \det(A)  \right|^{2}\cdot \Phi\left(TA\right),  &\text{for $A \in GL_2$;}   \\
			&\omega_{\psi}\left(n(x)\right)\Phi(T)=\psi \left(xQ(T)\right)\Phi(T),     &\text{for $x \in F$,}                  
		\end{aligned}
		\right.
	\end{equation}
	where $Q(T):=a_1b_4+a_2b_3-a_3b_2-a_4b_1$ for
	$$
	T=\begin{pmatrix}
		a_1& b_1 \\
		a_2&b_2 \\
		a_3&b_3 \\
		a_4&b_4
	\end{pmatrix} \in V_2\otimes W\cong M_{4\times 2}.
	$$
	Let 
	$$
	R:=GO(V)\times GSp(W)
	$$
	and
	$$
	R_0:=(GO(V)\times GSp(W))^{\circ}:=\{(h,g)\in GO(V)\times GSp(W)\mid \lambda_{V}(h)\cdot \lambda_{W}(g)=1 \}. 
	$$
	The Weil representation $\omega_{\psi}$ extends to the group $R_0$ via
	\begin{equation} \label{equation321}
		\omega_{\psi}(h,g)\Phi(T)=\left|\lambda_{W}(g) \right|^{-2}\omega_{\psi}(h_1,1)(\Phi\circ g^{-1})
	\end{equation}
	where
	$$
	h_1:=h \begin{pmatrix}
		\lambda_{V}(h)^{-1} & 0 \\
		0& 1
	\end{pmatrix}\in O(V).
	$$
	Then the action of $\left((M \ltimes N)\times GSp(W)\right)^{\circ}:=\{\left(m(A,\lambda_{W}(g)^{-1})n(x),g\right)\in (M \ltimes N)\times GSp(W) \}$ on $\omega_{\psi}$ is given by
	\begin{align}
		&\omega_{\psi}\left(m(A,\lambda_{W}(g)^{-1})n(x), g\right) \Phi(T)  \label{zf2: weilformula2} \\
		=&\left|\lambda_W(g)\det A\right|^2\psi\left(\lambda_{W}(g)\det(A) x\cdot Q(T)\right)\Phi\left(\lambda_{W}(g)g^{-1}TA\right). \notag 
	\end{align}
	Note in particular that the central elements $(t,t^{-1})\in R_0$ act trivially on $\omega_{\psi}$. The Weil representation $\Omega$ of $GO(V)\times GSp(W)$ is defined by the following compact induction
	$$
	\Omega:=ind_{R_0}^{R}\omega_{\psi}.
	$$
	Then one can define the similitude theta correspondence for $GO(V)\times GSp(W)$. Define the smooth big theta lift of $\sigma$ to $GSp(W)$ with respect to $\Omega$ by
	$$
	\Theta(\sigma):=\left(\Omega\otimes \sigma^{\vee} \right)_{GO(V)}.
	$$  
	Define the smooth small theta lift of $\sigma$ to $GSp(W)$ with respect to $\Omega$ by the maximal semisimple quotient of $\Theta(\sigma)$. By the works of \cite{H2,W3,GT2}, one knows $\theta(\sigma)$ is irreducible or 0. Likewise, define the smooth big (small) theta lift of $\pi$ to $GO(V)$ with respect to $\Omega$.

	By \cite[Lemma 2.2]{GT1}, one obtains
	\begin{lem}
		If $\sigma$ is an irreducible representation of $GO(V)$ ($resp.$ $GSp(W)$), and the restriction of $\sigma$ to the relevant isometry group is $\oplus_i\tau_i$, then as representation of $Sp(W)$ ($resp.$ $O(V)$),
		$$
		\Theta(\sigma)\cong \bigoplus_i\Theta_{\psi}(\tau_i).
		$$ 
		In particular, $\Theta(\sigma)$ is admissible of finite length. Moreover, if $\Theta_{\psi}(\tau_i)=\theta_{\psi}(\tau_i)$ for each $i$, then 
		$$
		\Theta(\sigma)=\theta(\sigma).
		$$
	\end{lem}
	By \cite[Corollary C.4]{GI2}, one obtains
	\begin{lem} \label{lemma2.2}
		For any irreducible tempered representation $\sigma$ of $GSO(V)$, one has
		$$
		\Theta(\sigma)=\theta(\sigma) \neq 0.
		$$ 
	\end{lem}
	
	If $\sigma_1\cong\sigma_2\cong \sigma$ as $GL_2$-representations, there are two extensions of the irreducible  representation $\sigma\boxtimes\sigma$ of $GSO(V)$ to $GO(V)$, which we denote by $(\sigma\boxtimes\sigma)^{\pm}$. By \cite[Theorem 7.8, Corollary 7.9]{R2}, one can deduce that exactly one of them participates in the theta correspondence with $GSp_2$, and we denote this distinguished extension by $(\sigma\boxtimes\sigma)^{+}$. If $\sigma_1 \neq \sigma_2$, then
	$$
	ind_{GSO(V)}^{GO(V)}(\sigma_1\boxtimes\sigma_2)\cong ind_{GSO(V)}^{GO(V)}(\sigma_2\boxtimes\sigma_1)
	$$ 
	is irreducible, and we denote $ind_{GSO(V)}^{GO(V)}(\sigma_1\boxtimes\sigma_2)$ by $(\sigma_1\boxtimes\sigma_2)^{+}$.
	We define the theta correspondence from $GSO(V)$ to $GSp(W)$ by
	$$
	\theta(\sigma_1\boxtimes\sigma_2):=\theta((\sigma_1\boxtimes\sigma_2)^{+}),
	$$
	where the theta correspondence on the right hand side is the theta correspondence from $GO(V)$ to $GSp(W)$. Hence this theta correspondence from $GSO(V)$ to $GSp(W)$ is generically $2$ to $1$.  
	
	\subsection{\bf Spherical function } \label{subsection: measure: weil}
	Assume $F$ is non-Archimedean and the conductor of $\psi$ is $\mathcal{O}_{F}$ in this subsection. Consider a characteristic function in $S(M_{4\times 2}(F))\cong S(V_2\otimes W)\cong \omega_{\psi}$ given by
	$$
	\Phi_{0,1}(T)=\left\{ 
	\begin{aligned}
		& 1,   \qquad \text{ if $T\in M_{4\times 2}(\mathcal{O}_{F})$;} \\
		& 0,  \qquad \text{ otherwise.}
	\end{aligned}
	\right.
	$$
	Let $L=M_{4\times 2}(\mathcal{O}_{F})\subset M_{4\times 2}(F)$. Then $\Phi_{0,1}$ is the characteristic function of $L$.
	
	One may use another model of the Weil representation 
	$$
	\omega_{\psi}\cong S(V\otimes W_2) \cong S(M_{2\times 4}).
	$$
	Then the isomorphism
	$$
	S(M_{4\times 2})\cong S(M_{2\times 4})
	$$
	is given by the partial Fourier transform defined by
	$$
	\widehat{\Phi}(u,v):=\int_{M_{2\times 2}(F)}\Phi\begin{pmatrix}
		x \\
		u
	\end{pmatrix}\psi(Tr(v\cdot (^{t}x)))d x,
	$$
	where the measure $d x$ is self-dual with respect to the pairing $(x,y)=\psi(Tr(\cdot,\cdot) )$.
	Then $\widehat{\Phi}_{0,1}$ is the characteristic function of $M_{2\times 4}(\mathcal{O}_{F})$. 
	We fix a measure $d_{\psi} T^{\prime}$ on $M_{2\times 4}(F)$ which is self-dual with respect to the pairing
	$$
	(x,y)=\psi\left(Tr(x\begin{pmatrix}
		& I_2 \\
		I_2 &  
	\end{pmatrix})\cdot (^{t}y) \right).
	$$ 
	Then one can fix a measure $d_{\psi} T$ on $M_{4\times 2}(F)$ such that
	$$
	\langle \Phi_{0,1},\Phi_{0,1} \rangle_{M_{4\times 2}(F)}=\langle \widehat{\Phi}_{0,1},\widehat{\Phi}_{0,1} \rangle_{M_{2\times 4}(F)}.
	$$
	By the formula of Weil representation action (for example, see \cite[Proposition 4.3]{K1}), one has
	$$
	O(V)(\mathcal{O}_{F})\times Sp(W)(\mathcal{O}_{F})\cdot \widehat{\Phi}_{0,1}=\widehat{\Phi}_{0,1}.
	$$
	Hence
	$$
	O(V)(\mathcal{O}_{F})\times Sp(W)(\mathcal{O}_{F})\cdot \Phi_{0,1}=\Phi_{0,1}.
	$$
	
	We denote $\langle \cdot ,\cdot \rangle_{M_{4\times 2}}$ by $\langle \cdot,\cdot \rangle_{\omega_{\psi}}$.

	\subsection{\bf Weil representation for $PGO(V)\times PGSp(W)$} \label{zf13: sectionweil}

	Let $Z_{V}$ and $Z_{W}$ be the center of $GO(V)$ and $GSp(W)$ respectively.  Fix the isomorphism $Z_{V}\cong Z_{W} \cong F^{\times}$. Let $Z^{\bigtriangledown}:=\{(z,z^{-1})\mid z \in F^{\times}\}\subset Z_{V}\times Z_{W}$. Observe that 
	$$
	R_0\cdot Z_{V}=R_0 \cdot Z_{W}=R_0\cdot(Z_{V}\times Z_{W}) \subset GO(V) \times GSp(W)
	$$
	with 
	$$
	R_0 \cap  (Z_{V}\times Z_{W})=\{(z,z^{-1})\mid z \in F^{\times}\}\cdot(Z_{O(V)}\times Z_{Sp(W)})=  Z^{\bigtriangledown}\cdot(Z_{O(V)}\times Z_{Sp(W)})
	$$
	and
	$$
	R_{0}\backslash R_0\cdot(Z_{V}\times Z_{W})=Z^{\bigtriangledown}\cdot(Z_{O(V)}\times Z_{Sp(W)}) \backslash Z_{V}\times Z_{W}.
	$$

	Note that the composite map
	$$
	\begin{CD}
		GO(V) @>\lambda_{V}>> F^{\times} @>>> F^{\times 2}\backslash F^{\times} 
	\end{CD}	
	$$
	induces an isomorphism
	$$
	O(V)\cdot Z_V \backslash GO(V) \cong F^{\times 2}\backslash F^{\times}.
	$$
	For each $a \in F^{\times 2}\backslash F^{\times} $, fix a representative $\tilde{a}\in F^{\times}$ such that
	\begin{enumerate}
		\item when $a=1 \pmod {F^{\times 2}}$, choose $\tilde{a}=1$;
		\item when $a \in \mathcal{O}_{F}^{\times} \pmod {F^{\times 2}}$, choose $\tilde{a} \in \mathcal{O}_{F}^{\times}$.
	\end{enumerate}

	Fix 
	$$
	h_a=\begin{pmatrix}
		\tilde{a}\cdot I_{2}& \\
		&I_{2}
	\end{pmatrix} \in GO(V).
	$$
	Then one has
	$$
	\lambda_{V}(h_a)  \equiv a \pmod {F^{\times 2}},
	$$
	and
	$$
	GO(V)=\bigcup_{a \in F^{\times 2}\backslash F^{\times }} O(V)\cdot Z_{V}\cdot h_a.
	$$
	
	For any $a \in \mathcal{O}_{F}^{\times 2}\backslash \mathcal{O}_{F}^{\times}$, we can define $\Phi_{0,a}$ in $\omega_{\psi_a}:=\omega_{\psi_{\widetilde{a}}}$ which is similar to the definition of $\Phi_{0,1}$ in Subsection \ref{subsection: measure: weil}. For any $a \in F^{\times 2}\backslash F^{\times}$, we can also define $\langle \cdot,\cdot \rangle_{\omega_{\psi_a}}$. 
	\subsection{\bf Smooth Weil representation for $PGO(V)\times PGSp(W)$}

	The smooth extended Weil representation for the similitude dual pair $GO(V)\times GSp(W)$ is defined by
	\begin{align}
		\Omega=&ind_{R_0}^{R}\omega_{\psi} \notag \\
		\cong& ind_{R_0\cdot(Z_{V}\times Z_{W})}^{R} ind_{R_0}^{R_0\cdot(Z_{V}\times Z_{W})}\omega_{\psi}. \notag
	\end{align}

	Let $\omega^{+}_{\psi}$ be the $1$-eigenspace of $\epsilon:=-1 \in Z_{O(V)}\cong \mu_2$ and $\omega^{-}_{\psi}$ be the $-1$-eigenspace of $\epsilon \in Z_{O(V)}\cong \mu_2$. Then $\omega_{\psi}=\omega^{+}_{\psi}\oplus \omega^{-}_{\psi}$ is a decomposition as $R_0$-module. Denote the non-trivial character of $\mu_2\cong Z_{O(V)}\cong Z_{Sp(W)}$ by $sign$. 
	
	\begin{lem}
		The $Z_{V}$-coinvariant of $\Omega$ is 
		$$
		\Omega_{Z_V}\cong ind_{R_0\cdot Z_{V}}^{R}\left(\omega_{\psi}^{+}\boxtimes  1 \right).
		$$
	\end{lem}
	\begin{proof}
		We denote $Z^{\bigtriangledown}\cdot(Z_{O(V)}\times Z_{Sp(W)})=Z^{\bigtriangledown}\times Z_{O(V)}$ by $\mathcal{Z}$. Let $\mathcal{Z}^{\Delta}=\{(z,z)\in \mathcal{Z}\times \mathcal{Z}\}$. Then one has 
		\begin{align}
			&ind_{R_0}^{R_0\cdot(Z_{V}\times Z_{W})}\omega_{\psi} \notag \\
			\cong &ind_{R_0}^{R_0\cdot(Z_{V}\times Z_{W})}\omega_{\psi}^{+}\bigoplus ind_{R_0}^{R_0\cdot(Z_{V}\times Z_{W})}\omega_{\psi}^{-} \notag \\
			\cong & ind_{\mathcal{Z}^{\Delta}\backslash R_0\times \mathcal{Z}}^{\mathcal{Z}^{\Delta}\backslash R_0\times (Z_V\times Z_W)}(\omega_{\psi}^{+}\boxtimes (1\boxtimes 1))\bigoplus ind_{\mathcal{Z}^{\Delta}\backslash R_0\times \mathcal{Z}}^{\mathcal{Z}^{\Delta}\backslash R_0\times (Z_V\times Z_W)}(\omega_{\psi}^{-}\boxtimes (1\boxtimes sign)) \notag \\
			\cong & ind_{ R_0\times \mathcal{Z}}^{ R_0\times (Z_V\times Z_W)}(\omega_{\psi}^{+}\boxtimes (1\boxtimes 1))\bigoplus ind_{ R_0\times \mathcal{Z}}^{ R_0\times (Z_V\times Z_W)}(\omega_{\psi}^{-}\boxtimes (1\boxtimes sign)) \notag \\
			\cong &\omega_{\psi}^{+}\boxtimes \left( ind_{  \mathcal{Z}}^{  Z_V\times Z_W} 1\boxtimes 1 \right)\bigoplus \omega_{\psi}^{-}\boxtimes \left( ind_{  \mathcal{Z}}^{  Z_V\times Z_W} 1\boxtimes sign \right) \notag \\
			\cong &\omega_{\psi}^{+}\boxtimes \left(  ind_{  Z^{\bigtriangledown}\times Z_{O(V)}}^{  Z^{\bigtriangledown}\times Z_{V}} 1\boxtimes 1 \right)\bigoplus \omega_{\psi}^{-}\boxtimes \left(  ind_{  Z^{\bigtriangledown}\times Z_{O(V)}}^{  Z^{\bigtriangledown}\times Z_{V}} 1\boxtimes sign \right) \notag \\
			\cong &\omega_{\psi}^{+}\boxtimes \left(  ind_{   Z_{O(V)}}^{   Z_{V}}  1 \right)\bigoplus \omega_{\psi}^{-}\boxtimes \left( ind_{   Z_{O(V)}}^{   Z_{V}}  sign \right) \notag 
		\end{align}

		Then we have the decomposition of $\Omega$:
		\begin{align}  
			\Omega\cong& ind_{R_0\cdot(Z_{V}\times Z_{W})}^{R} ind_{R_0}^{R_0\cdot(Z_{V}\times Z_{W})}\omega_{\psi} \notag \\
			\cong& ind_{R_0\cdot Z_{V}}^{R}\left(\omega_{\psi}^{+}\boxtimes ind_{Z_{O(V)}}^{Z_{V}} 1 \right)\bigoplus  ind_{R_0\cdot Z_{V}}^{R}\left(\omega_{\psi}^{-}\boxtimes ind_{Z_{O(V)}}^{Z_{V}} sign \right)\notag 
		\end{align}

		For any $(f_1,f_2) \in ind_{R_0\cdot Z_{V}}^{R}\left(\omega_{\psi}^{+}\boxtimes ind_{Z_{ O(V)}}^{Z_{V}} 1 \right)\bigoplus  ind_{R_0\cdot Z_{V}}^{R}\left(\omega_{\psi}^{-}\boxtimes ind_{Z_{O(V)}}^{Z_{V}} sign \right)$, we have
		$$
		\epsilon\cdot (f_1,f_2)-(f_1,f_2)=(f_1,-f_2)-(f_1,f_2)=(0,-2f_2).
		$$ 
		Hence
		\begin{align}
			\Omega_{Z_V}\cong &\left( ind_{R_0\cdot Z_{V}}^{R}\left(\omega_{\psi}^{+}\boxtimes  ind_{Z_{O(V)}}^{Z_V}1 \right)\right)_{Z_V} \notag \\
			\cong&ind_{R_0\cdot Z_{V}}^{R}\left(\omega_{\psi}^{+}\boxtimes  ind_{Z_{O(V)}}^{Z_V}1 \right)_{Z_V} \notag \\
			\cong& ind_{R_0\cdot Z_{V}}^{R}\left(\omega_{\psi}^{+}\boxtimes  \left(ind_{Z_{ O(V)}}^{Z_V}1 \right)_{Z_V}\right) \notag 
			\end{align}
		
		Since 
		\begin{align}
			&Hom\left(\left(ind_{Z_{O(V)}}^{Z_V}1 \right)_{Z_V}, \mathbb{C} \right) \notag \\
			=&Hom_{Z_V}\left(ind_{Z_{O(V)}}^{Z_V}1 , 1 \right) \notag \\
			=&Hom_{Z_V}\left(1,Ind_{Z_{O(V)}}^{Z_V}1\right)\notag \\
			=&Hom_{Z_{O(V)}}\left(1,1\right)\notag \\
			=&\mathbb{C},\notag 
		\end{align}
		one has
		$$
		\left(ind_{Z_{O(V)}}^{Z_V}1 \right)_{Z_V}=1.
		$$
		
		Therefore we have
		$$
		\Omega_{Z_V}\cong ind_{R_0\cdot Z_{V}}^{R}\left(\omega_{\psi}^{+}\boxtimes  1\right).
		$$
		
	\end{proof}

		We call the smooth $PGO(V)\times PGSp(W) $-representation
		$$
		\Omega(1):=ind_{R_0\cdot Z_V}^{R}\left(\omega_{\psi}^{+}\boxtimes 1\right)
		$$
		the Weil representation of $PGO(V)\times PGSp(W)$.

		For $\sigma \in Irr\left(PGO(V)\right)$, the smooth big theta lift of $\sigma$ to $PGSp(W)$ with respect to $\Omega(1)$ is:
		$$
		\Theta^{1 }(\sigma):=\left(\Omega(1)\otimes \sigma^{\vee} \right)_{PGO(V)}.
		$$
		Let us recall the smooth big theta lift of $\sigma$ to $GSp(W)$ with respect to $\Omega$:
		$$
		\Theta(\sigma):=\left(\Omega\otimes \sigma^{\vee} \right)_{GO(V)}.
		$$  
		Likewise, for $\pi \in Irr\left(PGSp(W)\right)$, the smooth big theta lift of $\pi$ to $PGO(V)$ with respect to $\Omega(1)$ is:
		$$
		\Theta^{1}(\pi):=\left(\Omega(1)\otimes \pi^{\vee} \right)_{PGSp(W)}.
		$$
		The smooth big theta lift of $\pi$ to $GO(V)$ with respect to $\Omega$ is:
		$$
		\Theta(\pi):=\left(\Omega\otimes \pi^{\vee} \right)_{GSp(W)}.
		$$  
		The theta correspondence defined by $\Omega(1)$ is the same as the theta correspondence defined by $\Omega$. 
		\begin{lem} \label{zf13: lemma300}
			\begin{enumerate}
				\item If $\sigma \in Irr\left(PGO(V) \right)$, then
				$$
				\Theta^{1}(\sigma)\cong\Theta(\sigma).
				$$
				Hence they have the isomorphic maximal semisimple quotients, i.e.\ $\theta^{1}(\sigma)\cong \theta(\sigma)$.
				\item If $\pi \in Irr\left(PGSp(W) \right)$, then
				$$
				\Theta^{1}(\pi)\cong\Theta(\pi).
				$$
				Hence they have the isomorphic maximal semisimple quotients, i.e.\ $\theta^{1}(\pi)\cong \theta(\pi)$.
			\end{enumerate}
		\end{lem}
		\begin{proof}
			Assume $\sigma \in Irr^{1}\left(PGO(V) \right)$. Then one has
			\begin{align}
				\Theta(\sigma)&=\left(\Omega\otimes \sigma^{\vee} \right)_{GO(V)}\notag \\
				=& \left(\left(\Omega\otimes \sigma^{\vee} \right)_{Z_V}\right)_{PGO(V)} \notag \\
				\cong& \left(\Omega_{Z_V}\otimes \sigma^{\vee}\right)_{PGO(V)} \notag \\
				\cong&\left(\Omega(1)\otimes\sigma^{\vee}\right)_{PGO(V)} \notag \\
				=&\Theta^{1}(\sigma). \notag
			\end{align}
		\end{proof}

		For $a \in F^{\times 2}\backslash F^{\times}$, let $\omega_{\psi_a}=\omega_{\psi}=S(V_2\otimes W)$ as vector space, and be equipped with the action of $O(V)\times Sp(W)$ on $\omega_{\psi_a}$ given by
		$$
		\omega_{\psi_a}(h,g) \Phi=\omega_{\psi}\left(h_a(h,g)h_a^{-1}\right) \Phi
		$$
		where the left action $\cdot$ is in $\omega_{\psi_a}$ while the right action $\cdot$ is in $\omega_{\psi}$.
		
		\begin{lem} \label{zf13: lemma1}
			The $O(V)\times Sp(W)$ map 
			$$
			\Omega(1)|_{O(V)\times Sp(W)}\cong\bigoplus_{a \in F^{\times 2}\backslash F^{\times }}\omega_{\psi_a}^{+}
			$$
			given by 
			$$
			\Phi \longmapsto \oplus_{a \in F^{\times 2}\backslash F^{\times}}\Phi(h_a)=:\oplus_{a \in F^{\times 2}\backslash F^{\times}}\Phi_a
			$$
			is a $O(V)\times Sp(W)$-equivariant isomorphism. Moreover, this isomorphism is $R_0$-equivariant.
		\end{lem}
		\begin{proof}
			For $(h,g)\in R_0$, one has
			\begin{align}
				(h,g)\cdot \Phi\mapsto &\oplus_{a}\Phi\left(h_a(h,g)\right) \notag \\
				=&\oplus_{a}\Phi\left(h_a(h,g)h_a^{-1}h_a\right) \notag \\
				=&\oplus_{a}h_a(h,g)h_a^{-1}\cdot \left(\Phi\left(h_a\right)\right) \notag \\
				=&(h,g)\cdot \oplus_a \Phi(h_a). \notag 
			\end{align}
		\end{proof}

		\subsection{\bf The actions} \label{subsection2.2b}
		Let $\Phi \in \Omega(1)\cong \oplus_{a\in F^{\times 2}\backslash F^{\times}}\omega_{\psi_a}^{+}\cong \oplus_{a\in F^{\times 2}\backslash F^{\times}}S\left(V_2^{(a)}\otimes W^{(a)}\right)^{+ }$, where $V_2^{(a)}\otimes W^{(a)}=V_2\otimes W$ and $S\left(V_2^{(a)}\otimes W^{(a)}\right)$ is the realization of the Weil representation with respect to $\psi_{a}:=\psi_{\widetilde{a}}$. Then $S\left(V_2^{(a)}\otimes W^{(a)}\right)= \omega_{\psi_a}$ as $M_1N\times Sp(W)$-modules.  For any $a \in F^{\times 2}\backslash F^{\times}$, fix 
		$$
		g_{a}:=\begin{pmatrix}
			\tilde{a}\cdot I_2 & \\
			& I_2
		\end{pmatrix} \in GSp(W).
		$$
		
		For $a,c \in F^{\times 2}\backslash F^{\times}$, there exists $d \in F^{\times}$ such that $\tilde{a}\cdot \tilde{ac}^{-1}=d^2\cdot \tilde{c}^{-1}$. Then one has
		\begin{equation} \label{zf2: action}
			\left\{\begin{aligned}
				& z\cdot \oplus_{a}\Phi_{a}(T)=\oplus_{a}\Phi_a(T),  \quad \text{ for $z \in Z_{V}$;}      \\
				& m(A,1)\cdot \oplus_{a}\Phi_{a}(T)=\left|\det(A)\right|^2\cdot\oplus_{a}\Phi_{a}( TA), \quad\text{for any $m(A,1) \in SO(V)$; } \\
				& n(x)\cdot \oplus_{a}\Phi_{a}(T)=\oplus_{a}\psi_{\tilde{a}}\left(xQ(T)\right)\Phi_{a}(T);\\
				& h_c\cdot \oplus_{a}\Phi_a(T)=\oplus_{a}\left|d \right|^4\cdot\Phi_{ac}\left(d\cdot T\right);      \\ 
				& z\cdot \oplus_{a}\Phi_{a}(T)=\oplus_{a}\Phi_a(T),  \quad \text{ for $z \in Z_{W}$;}      \\       
				& g\cdot \oplus_{a}\Phi_{a}(T)=\oplus_{a}\Phi_{a}(g^{-1}\cdot T), \quad\text{for any $g \in Sp(W)$; } \\
				& g_c\cdot \oplus_{a}\Phi_a(T)=\oplus_{a} \left|\tilde{a}\cdot\tilde{ac}^{-1} \right|^2\cdot\Phi_{ac}\left(d\cdot g_{c}^{-1}\cdot T\right).
			\end{aligned}\right.
		\end{equation}


		\subsection{\bf $L^2$-Weil representation for $PGO(V)\times PGSp(W)$} \label{zf13: subsectionL2weil}

		The $L^2$-Weil representation of the similitude dual pair $GO(V)\times GSp(W)$ is defined by
		$$
		\Omega:=Ind_{R_0}^{R}\omega_{\psi}.
		$$
		Then
		\begin{align}
			\Omega=&Ind_{R_0}^{R}\omega_{\psi} \notag \\
			\cong & Ind_{R_0\cdot(Z_{V}\times Z_{W})}^{R}Ind_{R_0}^{R_0\cdot(Z_{V}\times Z_{W})}\omega_{\psi}. \notag
		\end{align}

		Let $\omega^{+}_{\psi}$ be the $O(V) \times Sp(W)$-submodule of $\omega_{\psi}$ with $Z_{O(V)}\cong \mu_2$ acting trivially, and let  $\omega^{-}_{\psi}$ be the $O(V) \times Sp(W)$-submodule of $\omega_{\psi}$ with $Z_{O(V)}\cong \mu_2$ acting non-trivially. Then $\omega_{\psi}=\omega^{+}_{\psi}\oplus \omega^{-}_{\psi}$ is also a decomposition as $R_0$-module. Let $\widehat{Z}_{V}^{temp,+}$ and $\widehat{Z}_{V}^{temp,-}$ be the subset of the tempered unitary dual $\widehat{Z}_{V}^{temp}$ of $Z_V$ such that $Z_{O(V)}\cong \mu_2$ acts trivially and non-trivially, and $\mu_{Z_{V}}^{+}$ and $\mu_{Z_{V}}^{-}$ be the Plancherel measure on $\widehat{Z}_{V}^{temp,+}$ and $\widehat{Z}_{V}^{temp,-}$. Let $sign$ be the non-trivial character of $\mu_2$. Then we have the $Z_{V}$-spectral decomposition:
		\begin{align}  
			&Ind_{R_0}^{R_0\cdot(Z_{V}\times Z_{W})}\omega_{\psi} \notag \\
			\cong & \omega_{\psi}^{+}\widehat{\boxtimes}\left(Ind_{Z^{\bigtriangledown}\cdot(Z_{O(V)}\times Z_{Sp(W)})}^{Z_{V}\times Z_{W}} 1 \right)\bigoplus  \omega_{\psi}^{-}\widehat{\boxtimes} \left(Ind_{Z^{\bigtriangledown}\cdot(Z_{O(V)}\times Z_{Sp(W)})}^{Z_{V}\times Z_{W}} sign \right)  \notag \\
			\cong &\omega_{\psi}^{+}\widehat{\boxtimes} \left(Ind_{Z_{O(V)}}^{Z_{V}} 1 \right)\bigoplus  \omega_{\psi}^{-}\widehat{\boxtimes} \left(Ind_{Z_{O(V)}}^{Z_{V}} sign \right)\notag \\
			\cong&\omega_{\psi}^{+}\widehat{\boxtimes} \left(\int_{\widehat{Z}_{V}^{temp,+}} \chi d \mu_{Z_{V}}^{+}(\chi) \right)\bigoplus  \omega_{\psi}^{-}\widehat{\boxtimes} \left(\int_{\widehat{Z}_{V}^{temp,-}} \chi d \mu_{Z_{V}}^{-}(\chi) \right)\notag \\
			\cong&\int_{\widehat{Z}_{V}^{temp}} \omega_{\psi}^{\chi(-1)}\widehat{\boxtimes}\chi d \mu_{Z_{V}}(\chi).   \notag  
		\end{align}
		Then one has
		\begin{align*}
			\Omega\cong &Ind_{R_0\cdot(Z_{V}\times Z_{W})}^{R}Ind_{R_0}^{R_0\cdot(Z_{V}\times Z_{W})}\omega_{\psi} \notag \\
			\cong & Ind_{R_0\cdot Z_{V}}^{R}\int_{\widehat{Z}_{V}^{temp}} \omega_{\psi}^{\chi(-1)}\widehat{\boxtimes}\chi d \mu_{Z_{V}}(\chi) \notag \\
			\cong&\int_{\widehat{Z}_{V}^{temp}}Ind_{R_0\cdot Z_{V}}^{R}\left( \omega_{\psi}^{\chi(-1)}\widehat{\boxtimes}\chi  \right) d \mu_{Z_{V}}(\chi). \notag
		\end{align*}

		We call the unitary $PGO(V)\times PGSp(W)$-representation
		$$
		\Omega(1):=Ind_{R_0\cdot Z_{V}}^{R}\left( \omega_{\psi}^{+}\boxtimes 1 \right)\cong \oplus_{a \in F^{\times 2}\backslash F^{\times}}\omega_{\psi_a}^{+}
		$$
		the $L^2$-Weil representation of $PGO(V)\times PGSp(W)$. The inner product of $\Omega(1)$ is defined as 
		\begin{equation} \label{innerproductomega1}
			\langle \cdot,\cdot \rangle_{\Omega(1)}=\frac{1}{2}\cdot \left(\oplus_{a \in F^{\times 2}\backslash F^{\times}}\langle \cdot,\cdot  \rangle_{\omega_{\psi_a}}\right)
		\end{equation}
		
		\subsection{\bf The map $A_{\sigma}$ } \label{subsection22a}
		For every irreducible unitary representation $ \sigma$ of $PGSO(V)$, we fix a nonzero element $\theta_{\sigma} \in Hom_{R}\left(\Omega(1),\sigma\boxtimes\theta(\sigma)\right)$. Then define a $PGSO(V)$-invariant and $PGSp(W)$-equivariant map
		$$
		A_{\sigma}: \Omega(1)\otimes \sigma^{\vee}\twoheadrightarrow \theta(\sigma)
		$$
		by
		$$
		A_{\sigma}(\Phi,v)=\langle \theta_{\sigma}(\Phi),v \rangle_{\sigma}.
		$$

		\subsection{\bf Local doubling zeta integral} \label{subsection2.3b}

		For any tempered representation $\sigma \in \widehat{PGSO}(V)^{temp}$, $\Phi_1,\Phi_2\in \Omega(1)$, and $v_1,v_2\in \sigma $, define the local doubling zeta integral
		$$
		Z_{\sigma}^{PGSO(V)}(\Phi_1,\Phi_2,v_1,v_2):=\int_{PGSO(V)}\langle \Omega(1)(h)\Phi_1,\Phi_2  \rangle_{\Omega(1)}\cdot \overline{\langle \sigma(h)v_1,v_2 \rangle}_{\sigma}\  d^{PGSO(V)} h.
		$$
		By \cite[Lemma 9.5]{GI} and \cite[Lemma 3.1.4]{S3}, one knows $Z_{\sigma}^{PGSO(V)}(\Phi_1,\Phi_2,v_1,v_2)$ is absolutely convergent.

		\begin{lem}
			For every $\sigma \in \widehat{PGSO}(V)^{temp}$, one has
			\begin{align*}
				\theta(\sigma)\neq 0 \ \text{if and only if}\   Z_{\sigma}\neq 0.
			\end{align*}
		\end{lem}
		\begin{proof}
			One knows that $Z_{\sigma} \neq 0 $ if and only if there are $\Phi_1,\Phi_2 \in \omega_{\psi_a}$ for some $a \in F^{\times 2}\backslash F^{\times}$, and there are $v_1,v_2\in \tau$ for some constituent $\tau$ of $\sigma|_{SO(V)}$ such that 
			$$
			\int_{SO(V)}\langle \omega_{\psi_a}(h) \Phi_1, \Phi_2 \rangle_{\omega_{\psi_a}}\cdot \overline{\langle \tau(h)v_1,v_2 \rangle_{\tau}} d h \neq 0.
			$$
			By \cite[Proposition 11.5]{GQT}, this is equivalent to say that for some $a \in F^{\times 2}\backslash F^{\times}$ and some constituent $\tau$ of $\sigma|_{SO(V)}$, one has
			$$
			Hom_{SO(V)\times SO(V)}\left(\omega_{\psi_a}\otimes\overline{\omega_{\psi_a}}\otimes\overline{\tau}\otimes\tau,\mathbb{C} \right) \neq 0. 
			$$
			One knows the above holds if and only if for some $a \in F^{\times 2}\backslash F^{\times}$ and some constituent $\tau$ of $\sigma|_{SO(V)}$, one has
			$$
			\theta_{\psi_a}(\tau)\neq 0,
			$$
			which is equivalent to
			$$
			\theta(\sigma)\neq 0.
			$$

		\end{proof}
		
		\subsection{\bf Inner product on $\theta(\sigma)$} \label{subsection2.3c}
		For every $\Phi_1,\Phi_2 \in \Omega(1)$ and $v_1,v_2 \in\sigma$, the local doubling zeta integral
		\begin{equation} \label{local: eq2}
			Z_{\sigma}^{PGSO(V)}(\Phi_1,\Phi_2,v_1,v_2)=\langle A_{\sigma}(\Phi_1,v_1),A_{\sigma}(\Phi_2,v_2) \rangle_{\theta(\sigma)}
		\end{equation}
		for some Hermitian form $\langle \cdot,\cdot \rangle_{\theta(\sigma)}$ on $\Theta(\sigma)=\theta(\sigma)$ (by Lemma \ref{lemma2.2}). By the proof of \cite[Proposition 3.3.1]{S3} or \cite[Appendix A]{HLS}, $\langle \cdot,\cdot \rangle_{\theta(\sigma)}$ descends to a nonzero inner product on $\theta(\sigma)$. We fix this inner product on $\theta(\sigma)$. This equation is the local analogue of the Rallis inner product formula.

		\subsection{\bf The map $B_{\sigma}$} \label{zf13d: 1}
		Similarly, we can define a $PGSp(W)$-invariant and $PGSO(V)$-equivariant map
		$$
		B_{\sigma}: \Omega(1) \otimes \theta(\sigma)^{\vee} \longrightarrow \sigma
		$$
		by
		$$
		B_{\sigma}(\Phi,w)=\langle \theta_{\sigma}(\Phi), w \rangle_{\theta(\sigma)}.
		$$
		Then for any $\Phi\in \Omega(1)$, $v \in \sigma$ and $w \in \theta(\sigma)$, one has
		\begin{equation}  \label{local: eq300}
			\langle A_{\sigma}(\Phi,v),w  \rangle_{\theta(\sigma)}=\langle \theta_{\sigma}(\Phi), v\otimes w\rangle_{\sigma\boxtimes \theta(\sigma)}=\langle B_{\sigma}(\Phi,w),v  \rangle_{\sigma}.
		\end{equation}

		\subsection{\bf Plancherel formula for $\Omega(1)$}\label{zf13: 24}
		
		For any $\Phi_1,\Phi_2 \in \Omega(1)$ and $\sigma \in \widehat{PGSO}(V)^{temp}$, let 
		$$
		J^{\theta}_{\sigma}\left(\Phi_1,\Phi_2 \right)=\sum_{v\in ONB(\sigma)}\int_{PGSO(V)}\langle \Omega(1)(h)\Phi_1,\Phi_2 \rangle_{\Omega(1)}\overline{\langle \sigma(h) v,v \rangle}_{\sigma}d h.
		$$
		By the same method of \cite[Proposition 3.3.1]{S3}, one obtains
		\begin{prop} \label{zf13: lemma100}
			For every tempered representation $\sigma$ of $PGSO(V)$, the Hermitian form $J^{\theta}_{\sigma}$ is positive semi-definite, and
			$$
			\langle \Phi_1,\Phi_2  \rangle_{\Omega(1)}= \int_{\widehat{PGSO}(V)^{temp}} J^{\theta}_{\sigma}\left(\Phi_1,\Phi_2 \right)  d \mu_{PGSO(V)}(\sigma).
			$$
		\end{prop}
		Now we want to obtain the Plancherel decomposition of $\Omega(1)$.
		\begin{lem}
			For every irreducible representation $\sigma$ of $PGSO(V)$, one has
			$$
			J_{\sigma}^{\theta}\neq 0 \ \text{if and only if}\  Z_{\sigma}^{PGSO(V)}\neq 0.
			$$
		\end{lem}
		\begin{proof}
			It is easy to see that $J_{\sigma}^{\theta}\neq 0$ implies $Z_{\sigma}^{PGSO(V)}\neq 0$.
			
			If $Z_{\sigma}^{PGSO(V)}\left(\Phi_1,\Phi_2,v_1,v_2\right)=\langle A_{\sigma}(\Phi_1,v_1),A_{\sigma}(\Phi_2,v_2) \rangle_{\theta(\sigma)}\neq 0$, then for some $\Phi\in \Omega(1)$ and $v\in \sigma$, we have $A_{\sigma}(\Phi,v) \neq 0$, which implies 
			$$
			Z_{\sigma}^{PGSO(V)}\left(\Phi,\Phi,v,v\right)=\langle A_{\sigma}(\Phi,v),A_{\sigma}(\Phi,v) \rangle_{\theta(\sigma)} >0.
			$$ 
			Thus $J_{\sigma}^{\theta}(\Phi,\Phi)> 0$.
		\end{proof}
		For $\Phi_1,\Phi_2\in \omega_{\psi}$ and $v_1,v_2\in \sigma$, where $\sigma$ is an irreducible tempered representation of $PGSO(V)$, define the local doubling zeta integral
		\begin{equation} \label{zeta2}
			Z_{\sigma}^{SO(V)}(\Phi_1,\Phi_2,v_1,v_2):=\int_{SO(V)}\langle \omega_{\psi}(h)\Phi_1,\Phi_2  \rangle_{\omega_{\psi}}\cdot \overline{\langle \sigma(h)v_1,v_2 \rangle}_{\sigma}\  d^{SO(V)} h.
		\end{equation}
		For $\Phi_1,\Phi_2\in \omega_{\psi}$ and $v_1,v_2\in \tau$, where $\tau$ is an irreducible representation of $SO(V)$, define the local doubling zeta integral
		\begin{equation}\label{zeta3}
			Z_{\tau}^{SO(V)}(\Phi_1,\Phi_2,v_1,v_2):=\int_{SO(V)}\langle \omega_{\psi}(h)\Phi_1,\Phi_2  \rangle_{\omega_{\psi}}\cdot \overline{\langle \tau(h)v_1,v_2 \rangle}_{\tau}\  d^{SO(V)} h.
		\end{equation}
		\begin{lem}
			For any irreducible representation $\sigma$ of $PGSO(V)$, one has
			$$
			Z_{\sigma}^{PGSO(V)}\neq 0 \iff Z_{\sigma}^{SO(V)}\neq 0
			$$
		\end{lem}
		\begin{proof}
			If $Z_{\sigma}^{PGSO(V)}\neq 0$, then there are $\Phi_1, \Phi_2 \in \omega_{\psi}^{+}$  and $v_1,v_2 \in \sigma$ such that 
			\begin{align*}
				&\int_{PGSO(V)}\langle \Omega(1)(h)\Phi_1,\Phi_2 \rangle_{\Omega(1)}\overline{\langle \sigma(h) v_1,v_2 \rangle}_{\sigma}d^{PGSO(V)} h \\
				=&\int_{SO(V)}\langle \omega_{\psi}(h)\Phi_1,\Phi_2 \rangle_{\omega_{\psi}}\overline{\langle \sigma(h) v_1,v_2 \rangle}_{\sigma}d^{PGSO(V)} h \neq 0.
			\end{align*}
			Thus one has $Z_{\sigma}^{SO(V)}\neq 0$.

			If for some $\Phi_1, \Phi_2 \in \omega_{\psi}$  and $v_1,v_2 \in \sigma$, $Z_{\sigma}^{SO(V)}(\Phi_1,\Phi_2,v_1,v_2)\neq 0$, then let $\Phi_i^{\prime}:=\Phi_i+\omega_{\psi}(\epsilon)\Phi_i \in \omega_{\psi}^{+}$, where $\mu_2=\{1,\epsilon\}$ is the center of $SO(V)$. One has
			\begin{align*}
				&Z_{\sigma}^{SO(V)}\left(\Phi_1^{\prime},\Phi_2^{\prime},v_1,v_2 \right) \\
				=&2\cdot Z_{\sigma}^{SO(V)}\left(\Phi_1,\Phi_2,v_1,v_2\right)+2\cdot Z_{\sigma}^{SO(V)}\left(\omega_{\psi}(\epsilon)\Phi_1,\Phi_2,v_1,v_2\right) \\
				=&4\cdot Z_{\sigma}^{SO(V)}\left(\Phi_1,\Phi_2,v_1,v_2\right)\neq 0.
			\end{align*} 
			Thus $Z_{\sigma}^{PGSO(V)}(\Phi_1^{\prime},\Phi_2^{\prime},v_1,v_2)\neq 0$.
		\end{proof}
		
		\begin{lem} \label{zf13: lemma200}
			We have
			$$
			\Omega(1)\cong \int_{\widehat{PGSO}(V)^{temp}}\sigma\widehat{\boxtimes} \theta(\sigma)d \mu_{PGSO(V)}(\sigma)
			$$
		\end{lem}
		\begin{proof}
			From Lemma \ref{zf13: lemma100}, one has
			$$
			\Omega(1)\cong \int_{\widehat{PGSO}(V)^{temp}}\sigma\boxtimes M(\sigma)d \mu_{PGSO(V)}(\sigma)
			$$
			where $M(\sigma)$ is some unitary representation of $PGSp(W)$. Then 
			$M(\sigma)\neq 0$ holds if and only if $J_{\sigma}^{\theta}\neq 0$ holds. This is equivalent to $Z_{\sigma}^{PGSO(V)}\neq 0$, namely $\theta(\sigma)\neq 0$.
		
			Besides, the map
			$$
			\Omega(1) \twoheadrightarrow \sigma^{\infty}\boxtimes M(\sigma)^{\infty}
			$$
			factors through
			$$
			\sigma^{\infty}\boxtimes\Theta(\sigma^{\infty})\twoheadrightarrow \sigma^{\infty}\boxtimes M(\sigma)^{\infty}.
			$$
			Since $\Theta(\sigma^{\infty})$ has finite length, one knows $\Theta(\sigma^{\infty})$ is admissible. Because $M(\sigma)^{\infty}$ is a quotient of $\Theta(\sigma^{\infty})$, we know $M(\sigma)^{\infty}$ is also admissible. Moreover, $M(\sigma)$ is unitary, so we know $M(\sigma)^{\infty}$ is semisimple. Thus
			the above map factors through
			$$
			\sigma^{\infty}\boxtimes\theta(\sigma)^{\infty}\twoheadrightarrow \sigma^{\infty}\boxtimes M(\sigma)^{\infty}.
			$$
			Hence we have $M(\sigma)=\theta(\sigma)$.
		\end{proof}

		Similarly, one has the spectral decompositions 
		$$
		\Omega(1)\cong \int_{\widehat{PGO}(V)^{temp}}\sigma\widehat{\boxtimes} \theta(\sigma)d \mu_{PGO(V)}(\sigma),
		$$
		and
		$$
		\Omega\cong \int_{\widehat{GO}(V)^{temp}}\sigma\widehat{\boxtimes} \theta(\sigma)d \mu_{GO(V)}(\sigma).
		$$
		Then we have the following compatibility:
		\begin{lem}
			Let $\sigma \in \widehat{PGO}(V)^{temp}$. Then the theta correspondence of $\sigma$ with respect to $\Omega(1)$ is equal to the theta correspondence of $\sigma$ with respect to $\Omega$.
			\end{lem}


\section{Plancherel Formula}
\subsection{\bf Plancherel formula of $L^2\left(\SU(2)\backslash \SO_{2,3} \right)$}
\label{zf2: section2}

For $b \in F$ and $a \in F^{\times 2}\backslash F^{\times}$, consider the $a$-th copy $Y_{b}^{(a)} \subset V_{2}^{(a)}\otimes W^{(a)}$ of $Y_{b}$ ( for definition of $Y_b$, see Subsection \ref{subsection: 3.1a} ). Define
$$
^{b}Y:=\sqcup_{a} Y_{\tilde{a}\cdot b}^{(a)}. 
$$
We endow $\sqcup_{a} V_2^{(a)}\otimes W^{(a)}$ with the measure $\dd T:=\sqcup_{a} \dd_{\psi_a}T$ which is $\PGSp(W)$-invariant.
We endow $^{b}Y=\sqcup_{a} Y_{\tilde{a}\cdot b}^{(a)}$ with a measure $\left|\omega_{b} \right|=\sqcup_{a}\left|\omega_{b}^{(a)} \right| $ such that
for any $f \in C_c^{\infty}\left(\sqcup_{a}V_2^{(a)}\otimes W^{(a)}\setminus \sqcup_{a}Y^{(a)}_{0}\right)$, one has
$$
\int_{\sqcup_{a}V_2^{(a)}\otimes W^{(a)}}f(T)\dd T=\int_{F^{\times}}\left(\int_{^{b}Y}f(T)\cdot \left|\omega_{b}\right|\right)\dd_{\psi}b
$$
where the function of $b \in F^{\times}$ defined by the inner integral on the right-hand-side is smooth and compactly supported.

	\subsection{\bf Spectral decomposition of $L^2\left(\SU(2)\backslash \PGSp(W)\right)$} \label{subsection4.1b}
	
	Let 
	$$
	L^2\left({^{1}Y} \right)^{+}:=\{ f\in L^2\left({^{1}Y} \right)\mid f(\epsilon\cdot T)=f(T) \}
	$$
	where $\epsilon$ is the nontrivial element in $Z_{\Sp(W)}$.
	Recall $B\subset \overline{\GL}_2^{(1)}$ is the subgroup of upper triangular matrices, and $\overline{\GL}_2^{(1)}\times\overline{M}_{\psi}=\PGSO(V)$. We endow the space $L^2\left({^{1}Y}, \left|\omega_{1}\right| \right)^{+}$ with the natural $\overline{M}_{\psi}\times \PGSp(W)$-action inherited by $\Omega(1)$. 
	
	We will fix the Haar measure on $F$ characterized by the requirement that $\dd x$ is self-dual with respect to the Fourier transform relative to the pairing $(x,y)\mapsto \psi(xy)$ on $F$. We will fix the Haar measure on $F^{\times}$ given by 
	$$
	\dd^{\times} x = \frac{\dd x}{|x|}.
	$$
	Recall 
	$$
	\overline{B}=S\cdot N,
	$$
	where $\overline{B}:=F^{\times}\backslash B$, and
	$$
	S=\left\{s(y)=\begin{pmatrix}
		y& \\
		&1 
	\end{pmatrix} \middle\vert\, y \in F^{\times}\right\}.
	$$ 
	We have Haar measures on $N(F)=F$ and $S(F)=F^{\times}$ and hence a right-invariant measure on $\overline{B}$. The Haar measures we have fixed endow the latter space with a unitary structure with inner product 
	$$
	\langle f_1,f_2 \rangle_{N\backslash \overline{B}}=\int_{F^{\times}}f_1(s(y))\cdot \overline{f_2(s(y))}\dd^{\times} y.
	$$

	We define the action of $F^{\times}$ on $\sqcup_{a}V_{2}^{(a)}\otimes W^{(a)}\slash \{\pm 1\}$ by
	$$
	y\ast \left(v_{a} \right)_{a}=\left( x\cdot v_{a y} \right)_{a}
	$$ 
	where $x^2=\widetilde{a}\cdot y\cdot \widetilde{ay}^{-1}$. For $b \in F^{\times}$, the multiplication $b \ast$ gives an isomorphism of varieties 
	$$
	\lambda_{b}: {^{c}Y}\slash \{\pm 1 \} \longrightarrow {^{cb}Y}\slash \{\pm 1 \}.
	$$ 
	We equip the quotient measure of $|\omega_{c}|$ on $^{c}Y\slash \{\pm 1 \}$ under the natural map
	$$
	^{c}Y \longrightarrow {^{c}Y}\slash \{\pm 1 \}.
	$$
	We still denote by $|\omega_{c}|$ this quotient measure. Then one may consider the pushforward measure $(\lambda_{b})_{\ast}(|\omega_{c}|)$ on $^{cb}Y\slash \{\pm 1 \}$.
	\begin{lem}
		One has
		$$
		(\lambda_{b})_{\ast}\left(\left|\omega^{(a)}_{c\cdot b^{-1}}\right|\right)=\left|b \right|\cdot \left|\widetilde{a}\cdot b\cdot \widetilde{ab}^{-1}\right|^{-4}\cdot\left|\omega_{c}^{(a)}\right|.
		$$
		
	\end{lem}
	\begin{proof}
		For $ f \in C_c^{\infty}\left(\sqcup_{a}V_2^{(a)}\otimes W^{(a)}\setminus \sqcup_{a}Y_0^{(a)}\right)^{+}$ and  $b \in F^{\times}$, one has
		\begin{align*}
			&\int_{\sqcup_{a}V_2^{(a)}\otimes W^{(a)}}f( b\ast T)\dd T  \\
			=&\sum_{a}\left|\widetilde{a}\cdot b\cdot \widetilde{ab}^{-1}\right|^{-4} \int_{V_2^{(ab)}\otimes W^{(ab)}}f(T)\dd T  \\
			=&\int_{\sqcup_{a}V_2^{(a)}\otimes W^{(a)}}\widetilde{f}( T)\dd T\\
			=&\int_{F^{\times}}\left(\int_{^{c}Y}\widetilde{f}\left(T\right)\cdot \left|\omega_{c}\right|\right)\dd_{\psi}c 
		\end{align*}
		where $\widetilde{f}(T)=\left|\widetilde{a}\cdot b\cdot \widetilde{ab}^{-1}\right|^{-4}\cdot f(T)$ when $T \in V_2^{(ab)}\otimes W^{(ab)}$.

		On the other hand, one has
		\begin{align*}
			&\int_{\sqcup_{a}V_2^{(a)}\otimes W^{(a)}}f( b\ast T)\dd T \\
			=&\int_{F^{\times}}\left(\int_{^{c}Y}f\left( b\ast T\right)\cdot \left|\omega_{c}\right|\right)\dd_{\psi}c \\
			=&\int_{F^{\times}}\left(\int_{^{c}Y}\lambda_{b}^{\ast}(f)\left(T\right)\cdot \left|\omega_{c}\right|\right)\dd_{\psi}c \\
			=&\int_{F^{\times}}\left(\int_{^{bc}Y}f\left(T\right)\cdot (\lambda_{b})_{\ast}\left(\left|\omega_{c}\right|\right)\right)\dd_{\psi}c \\
			=&\int_{F^{\times}}\left(\int_{^{c}Y}f\left(T\right)\cdot (\lambda_{b})_{\ast}\left(\left|\omega_{c\cdot b^{-1}}\right|\right)\right)\left|b\right|^{-1}\dd_{\psi}c.
		\end{align*}
		Comparing the above equations, we obtain 
		$$
		(\lambda_{b})_{\ast}\left(\left|\omega^{(a)}_{c\cdot b^{-1}}\right|\right)=\left|b \right|\cdot \left|\widetilde{a}\cdot b\cdot \widetilde{ab}^{-1}\right|^{-4}\cdot\left|\omega_{c}^{(a)}\right|.
		$$

	\end{proof}

	\begin{lem} \label{lemma4.2}
		For $ f \in C_c^{\infty}\left(	\sqcup_aV_2^{(a)}\otimes W^{(a)}\setminus \sqcup_aY_0^{(a)}\right)^{+}$, one has
		$$
		\int_{\sqcup_aV_2^{(a)}\otimes W^{(a)}}f(T)\dd T=\int_{F^{\times}}\left(\int_{{^{1}Y}}\widetilde{f}\left(c\ast T\right)\cdot \left|\omega_{1}\right|\right)\cdot \frac{\dd_{\psi}c}{\left|c\right|},
		$$
		where $\widetilde{f}(T)=\left|\widetilde{a}\cdot c\cdot \widetilde{ac}^{-1}\right|^{4}\cdot f(T)$ when $T \in V_2^{(a)}\otimes W^{(a)}$.

		Hence the map 
		$$
		\varrho^{\prime\prime}:	L^2\left(\sqcup_aV_2^{(a)}\otimes W^{(a)},\dd T \right)\longrightarrow L^2\left(F^{\times}\times {^{1}Y},\dd_{\psi}^{\times}b\times \left|\omega_{1}\right|\right)
		$$
		given by
		$$
		\varrho^{\prime \prime}(f)|_{Y^{(a)}_{\widetilde{a}}}(b,T)=\left|\widetilde{a}\cdot b \cdot \widetilde{ab}^{-1}\right|^{2}\cdot f\left( b\ast T \right)=m(b\cdot I_2,1)\cdot f|_{Y^{(a)}_{\widetilde{a}}}
		$$
		is an isometric isomorphism.
	\end{lem}
	\begin{proof}
		For $ f \in C_c^{\infty}\left(	\sqcup_aV_2^{(a)}\otimes W^{(a)}\setminus \sqcup_aY_0^{(a)}\right)^{+}$, one has
		\begin{align*}
			&\int_{\sqcup_aV_2^{(a)}\otimes W^{(a)}}f(T)\dd T  \\
			=&\int_{F^{\times}}\left(\int_{^{c}Y}f\left(T\right)\cdot \left|\omega_{c}\right|\right)\dd_{\psi}c \\
			=&\int_{F^{\times}}\left(\int_{^{c}Y}\widetilde{f}\left(T\right)\cdot (\lambda_c)_{\ast}\left(\left|\omega_{1}\right|\right)\right)\cdot \frac{\dd_{\psi}c}{\left|c\right|} \\
			=&\int_{F^{\times}}\left(\int_{{^{1}Y}}\lambda_{c}^{\ast}\left(\widetilde{f}\right)\left(T\right)\cdot \left|\omega_{1}\right|\right)\cdot \frac{\dd_{\psi}c}{\left|c\right|} \\
			=&\int_{F^{\times}}\left(\int_{{^{1}Y}}\widetilde{f}\left(c\ast T\right)\cdot \left|\omega_{1}\right|\right)\cdot \frac{\dd_{\psi}c}{\left|c\right|} 
		\end{align*}
		where $\widetilde{f}(T)=\left|\widetilde{a}\cdot c\cdot \widetilde{ac}^{-1}\right|^{4}\cdot f(T)$ when $T \in V_2^{(a)}\otimes W^{(a)}$.
	\end{proof}
	

	Then one can conclude that
	\begin{lem} \label{lemma4.3}
		The map
		$$
		\varrho^{\prime}: \Omega(1)\cong L^2\left(\sqcup_{a} V_2^{(a)}\otimes W^{(a)},\dd T \right)^{+}\cong L^2\left(N\backslash \overline{B},\psi;\dd_{\psi}^{\times}y\right)\widehat{\otimes}L^2\left({^{1}Y}, \left|\omega_{1}\right| \right)^{+}
		$$
		given by
		$$
		\varrho^{\prime}(f)(b,T)=b\cdot f(T)
		$$
		is a $\overline{B}\times \overline{M}_{\psi}\times \PGSp(W)$-isometry.
	\end{lem}

	We can identify $\SU(2)\backslash \PGSp(W)$ with $^{1}Y\slash \{\pm 1 \}$ by
	$$
	Z_{W}\cdot g \cdot g_c \mapsto \widetilde{c}\cdot g_{c}^{-1}\cdot g^{-1}\cdot T_1 \in Y^{(c)}_{\widetilde{c}}\slash \{\pm 1 \}
	$$
	where the definition of $T_1$ is given in Subsection \ref{subsection: 3.1a}, and $g \in \Sp(W)$.

	
	We equip the $\PGSp(W)$-action on $L^2\left(\SU(2)\backslash \PGSp(W) \right)$ by right translation; we equip the $\overline{M}_{\psi}$-action on $L^2\left(\SU(2)\backslash \PGSp(W) \right)$ by left multiplication via the identification
	$$
	\overline{M}_{\psi} \longrightarrow  \Delta F^{\times}\backslash \left(\GL_2\times \GL_2 \right)_{\mathrm{det}}
	$$
	given by
	$$
	A \mapsto P\left( A,  A\right).
	$$

	Similarly to Lemma \ref{geometric: lemma1}, one obtains
	\begin{lem} \label{lem4.4}
		There is an $\overline{M}_{\psi}\times \PGSp(W)$-isomorphism
		$$
		\varrho: L^2\left(^{1}Y \right)^{+}\cong L^2\left(\SU(2)\backslash \PGSp(W) \right)
		$$
		given by 
		$$
		\Phi=\oplus_a\Phi_a\mapsto f(Z_{W}\cdot g\cdot g_a):=\left|\tilde{a}\right|^{2}\cdot \Phi_{a}(\widetilde{a}\cdot g_a^{-1}\cdot g^{-1}\cdot T_1)=Z_{W}\cdot g\cdot g_a \cdot \Phi(T_1),
		$$
		where $g \in \Sp(W)$.
	\end{lem}
	We will fix a $\PGSp(W)$-invariant measure $\frac{\dd g}{\dd s}$ on $\SU(2)\backslash \PGSp(W)$ such that 
	$$
	L^2\left(^{1}Y \right)^{+}\cong L^2\left(\SU(2)\backslash \PGSp(W) \right)
	$$
	is isometric.
	


	\begin{lem}
		There is a $\overline{B}\times \overline{M}_{\psi}\times \PGSp(W)$-isometric isomorphism
		$$
		\Omega(1)\cong L^2\left( N\backslash \overline{B},\psi\right)\widehat{\bigotimes}L^2\left(\SU(2)\backslash \PGSp(W) \right).
		$$
	\end{lem}
	
	\begin{lem}
		There is a $\overline{B}\times \overline{M}_{\psi}\times \PGSp(W)$-isometric isomorphism
		$$
		\Omega(1)\cong L^2\left( N\backslash \overline{B},\psi\right)\widehat{\bigotimes} \int_{\widehat{\PGSO}(V)^{\temp}}\sigma_2\widehat{\boxtimes}\theta(\sigma_1\widehat{\boxtimes}\sigma_2)\dd\mu_{\PGSO(V)}(\sigma_1\widehat{\boxtimes}\sigma_2).
		$$
	\end{lem}
	\begin{proof}
		By the spectral decomposition of Weil representation, one knows that as $\PGSO(V)\times \PGSp(W)$-representation
		$$
		\Omega(1)\cong\int_{\widehat{\PGSO}(V)^{\temp}}(\sigma_1\widehat{\boxtimes} \sigma_2)\widehat{\boxtimes} \theta(\sigma_1\widehat{\boxtimes} \sigma_2)\dd\mu_{\PGSO(V)}(\sigma_1\widehat{\boxtimes} \sigma_2).
		$$ 
		By \cite[Lemma 4.4]{LM} and \cite[Proposition 2.14.3]{BP2}, one also has a $\overline{B}$-isometry
		$$
		\rest|_{\overline{B}}\circ\beta_{\sigma_1}^{(N,\psi)\backslash \PGL_2}:\sigma_1\cong L^2\left( N\backslash \overline{B},\psi \right),
		$$
		where the map $\beta_{\sigma_1}^{(N,\psi)\backslash \PGL_2}$ comes from the Plancherel decomposition of $L^2(N\backslash \PGL_2,\psi)$ and $\rest|_{\overline{B}}$ is the restriction map.
		Then as a $\overline{B}\times \overline{M}_{\psi}\times \PGSp(W)$-representation, one has
		$$
		\Omega(1)\cong L^2\left( N\backslash \overline{B} ,\psi \right)\widehat{\bigotimes}\int_{\widehat{\PGSO}(V)^{\temp}}\sigma_2\widehat{\boxtimes} \theta(\sigma_1\widehat{\boxtimes} \sigma_2)\dd\mu_{\PGSO(V)}(\sigma_1\widehat{\boxtimes} \sigma_2)
		$$

	\end{proof}
	\begin{prop}
		As an $\overline{M}_{\psi}\times \PGSp(W)$-representations,
		$$
		L^2\left(\SU(2)\backslash \PGSp(W) \right)\cong \int_{\widehat{\PGSO}(V)^{\temp}}\sigma_2\widehat{\boxtimes}\theta(\sigma_1\widehat{\boxtimes}\sigma_2)\dd\mu_{\PGSO(V)}(\sigma_1\widehat{\boxtimes}\sigma_2)
		$$
	\end{prop}

	\subsection{\bf Transfer of test functions} \label{subsection: transfer}
	In this subsection, we define the maps $p_1$, $q_1$, $p$ and $q$ which will be used in the following subsections. Then we define the transfer of test functions which will be used in Chapter \ref{chapter5}.

	
	Let $\mathrm{pr}_1: \Omega(1)\cong \bigoplus_{a\in F^{\times 2}\backslash F^{\times}}\omega_{\psi_a}^{+}\longrightarrow \omega_{\psi}^{+}$ be the projection map.
	Let $p_1$ be the map
	$$
	p_1: \Omega(1)\cong \bigoplus_{a\in F^{\times 2}\backslash F^{\times}}\omega_{\psi_a}^{+}\longrightarrow C^{\infty}\left(N\backslash \PGSO(V),\psi\right)
	$$
	given by
	$$
	p_1(\Phi)(h):=\mathrm{pr}_1\left(h\cdot \Phi \right)(T_1) .
	$$
	Let $\mathcal{S}\left(N\backslash \PGSO(V),\psi \right):=\mathrm{Im}(p_1)$. 
	\begin{lem}
		For any $\Phi\in \Omega(1)$ and $h \in \PGSO(V)$, the function
		$$
		f: \overline{T}_{E} \longrightarrow  \mathbb{C}
		$$
		given by
		$$
		f(t,h)=p_1(\Phi)\left(m(t,\det(t))\cdot h \right)
		$$
		satisfies
		$$
		\int_{\overline{T}_{E}}\left|f(t,h)\right|\dd t<\infty.
		$$
	\end{lem}
	\begin{proof}
		This is clear when $\overline{T}_{E}$ is compact. Suppose $\overline{T}_{E}=F^{\times}$. For
		$$
		s(y)=\begin{pmatrix}
			y& \\
			&1
		\end{pmatrix} \in \overline{T}_{E},
		$$
		one has
		$$
		f(s(y),h)=\left|y\right|^2\cdot \mathrm{pr}_1\left(m(I_2,y)\cdot h\cdot\Phi\right)(T_1\cdot s(y)).
		$$
		For $a \in F^{\times 2}\backslash F^{\times}$, let $y=\tilde{a}\cdot x^2$. Then
		\begin{align*}
			f(s(\tilde{a}\cdot x^2),h)=&\left|\tilde{a}\cdot x^2\right|^2\cdot \mathrm{pr}_1\left( m(I_2,x^{2})\cdot m(I_2,\tilde{a})\cdot h\cdot\Phi\right)(T_1\cdot s(y)) \\
			=&\left|\tilde{a}\right|^2\cdot  \mathrm{pr}_1\left(m(I_2,\tilde{a})\cdot h\cdot\Phi \right)\left(x^{-1}\cdot T_1\cdot s(\tilde{a}\cdot x^2)\right).
		\end{align*}
		Here 
		$$
		T(x):=x^{-1}\cdot T_1\cdot s(\tilde{a}\cdot x^2)=\begin{pmatrix}
			0&0 \\
			\widetilde{a}\cdot x&0 \\
			0&x^{-1} \\
			0&0
		\end{pmatrix},
		$$
		and 
		$$
		\Phi_1:=\left|\tilde{a}\right|^2\cdot  \mathrm{pr}_1\left(m(I_2,\tilde{a})\cdot h\cdot\Phi \right)\in \omega_{\psi}^{+}.
		$$
		We need only to prove
		$$
		\int_{F^{\times}}\left|\Phi_1(T(x))\right|\dd^{\times} x <\infty.
		$$
		When $F$ is non-Archimedean,  $\Phi_1 \in C_c^{\infty}\left(V_2\otimes W\right)$. Consider the projection map $\mathrm{pr}_{i,j}$ from $M_{4\times 2}(F)$ to its ($i,j$)-component. Then $\mathrm{pr}_{i,j}(\supp(\Phi_1))$ is compact in $F$. Thus we have 
		$$
		\left|\mathrm{pr}_{i,j}(\supp(\Phi_1)) \right|<C
		$$
		for some positive number $C> 1$. In the non-Archimedean case,
		\begin{align*}
			\int_{F^{\times}}\left|\Phi_1(T(x))\right|\dd^{\times} x=&\int_{C^{-1}<\left|x\right|<C} \left|\Phi_1(T(x))\right|\dd^{\times} x \\
			=&\int_{C^{-1}<\left|x\right|<C} \frac{\left|\Phi_1(T(x))\right|}{\left|x\right|}\dd x \\
			<&\infty.
		\end{align*}

		When $F$ is Archimedean, $\Phi_1 \in S\left(V_2\otimes W\right)$. Thus $\left(\left|x\right|^2+\left|x\right|^{-2}\right)\cdot \Phi_1(T(x))< C$ for some positive constant $C$. When $\left|x\right|\to \infty$, one has 
		$$
		\left|\Phi_1(T(x))\right|<\frac{C}{\left|x\right|^2};
		$$
		when $\left|x\right|\to 0$, one has
		$$
		\left|\Phi_1(T(x))\right|<C\cdot \left|x\right|^{2}.
		$$

		Then we can conclude that 
		$$
		\int_{F^{\times}}\left|\Phi_1(T(x))\right|\dd^{\times} x <\infty.
		$$

	\end{proof}
	
	Let $p_2$ be the map
	$$
	p_2: \mathcal{S}\left(N\backslash \PGSO(V),\psi\right) \longrightarrow C^{\infty}\left(N\times \overline{T}_{E}\backslash \PGSO(V),\psi\right)
	$$
	given by
	$$
	p_2(p_1(\Phi))(h):=\int_{\overline{T}_{E}}p_1(\Phi)\left(m(t,\det(t))\cdot h \right)  \dd t
	$$
	We denote the composite map $p_2\circ p_1$ by $p$.
	
	Let $q_1$ be the map
	$$
	q_1: \Omega(1)\cong \bigoplus_{a\in F^{\times 2}\backslash F^{\times}}\omega_{\psi_a}^{+}\longrightarrow C^{\infty}\left(\SU(2)\backslash \PGSp(W) \right)
	$$
	given by
	$$
	q_1(\Phi)(g):=\mathrm{pr}_1\left(g\cdot\Phi\right)(T_1)
	$$
	\begin{lem}
		$$
		\mathrm{Im}(q_1)=S\left(\SU(2)\backslash \PGSp(W)\right),
		$$
		where $S\left(\SU(2)\backslash \PGSp(W)\right)$ is the space of Schwartz-Bruhat functions of $\SU(2)\backslash \PGSp(W)$.
	\end{lem}
	Let $q_2$ be the map
	$$
	q_2: S\left(\SU(2)\backslash \PGSp(W)\right) \longrightarrow S\left(\U(2)\backslash \PGSp(W)\right)
	$$
	given by
	$$
	q_2\circ q_1(\Phi)(g):=\int_{\overline{T}_{E}}\mathrm{pr}_1\left(P(t,t)\cdot g\cdot\Phi\right)(T_1) \dd t.
	$$
	Let $L$ be the left translation action of $\overline{T}_{E}$ on $S\left(\SU(2)\backslash \PGSp(W)\right)$ given by
	$$
	L(t)f(g):=f\left(P(t,t)^{-1}\cdot g\right).
	$$
	Let $q=q_2\circ q_1$.

	Let $X=(N,\psi)\times \overline{T}_{E}\backslash \PGSO(V)$ and $Y=\U(2)\backslash \PGSp(W)$. Then we have the following diagram
	\begin{equation} \label{local: transfer}
		\xymatrix{
			&\Omega(1) \ar[ld]_{p_1}   \ar[rd]^{q_1}\\
			\mathcal{S}({N,\psi}\backslash \PGSO(V))    \ar[d]_{p_2} &        &S(\SU(2)\backslash \PGSp(W)) \ar[d]^{q_2} \\
			\mathcal{S}^{\heartsuit}(X) & & S(Y)                  }
	\end{equation}
	Note that $q_2$-maps and $q$-maps land in usual Schwartz space $S(Y)$, and $p_2$-maps and $p$-maps do not.


	Now we can define the correspondence between the two spaces.
	\begin{defn} \label{transfer: defn2}
		We say that $\varphi_1 \in  \mathcal{S}^{\heartsuit}(X)$ and $\varphi_2 \in S(Y)  $ are in correspondence if there exists $\Phi \in \Omega(1)$ such that $p(\Phi)=\varphi_1$ and $q(\Phi)=\varphi_2$.
	\end{defn}
	If $\varphi_1$ and $\varphi_2$ are in correspondence, we will call $\varphi_2$ ($\varphi_1$) is a transfer of $\varphi_1$ ($\varphi_2$). Since the maps $p$ and $q$ are surjective, we know
	\begin{lem} \label{transfer: prop2}
		Every $\varphi_1 \in  \mathcal{S}^{\heartsuit}(X)$ has a transfer $\varphi_2 \in S(Y)  $ and vice versa.
	\end{lem}
	We also have
	\begin{lem}If $\overline{T}_{E}$ is compact, then we have 
		$$
		\mathcal{S}^{\heartsuit}(X) \subset \mathcal{C}(X). 
		$$
	\end{lem}
	
	\begin{proof}
		
		Up to conjugating $\overline{T}_{E}$, we may assume that
		$$
		\overline{M}_{\psi}= \PGL_2=\overline{T}_{E}\cdot S \cdot K,
		$$
		where $K$ is some compact subset of $\PGL_2$ (\cite[Lemma 5.3.1]{SV}, \cite{BO}), and  
		$$
		S=\left\{\begin{pmatrix}
			y& \\
			& 1
		\end{pmatrix}\middle\vert\, y \in F^{\times}\right\}=\bigsqcup_{a \in F^{\times}\slash F^{\times 2}}\left\{s(\widetilde{a}\cdot y^2)=\begin{pmatrix}
			\widetilde{a}\cdot y^{2}& \\
			& 1
		\end{pmatrix}\middle\vert\, y \in F^{\times}\right\}. 
		$$
		Let 
		$$
		f(g):=p\left(g\cdot \Phi\right) \in C^{\infty}\left(X\right),
		$$
		where $g \in \PGSO(V)$.
		
		We consider the growth of 
		$$
		f\left(s(x^2)\right)=p\left(s(x^2)\cdot \Phi\right)=p\left(h(x)\cdot \Phi\right),
		$$
		where $\Phi\in \Omega(1)$ and 
		$$
		h(x)=\begin{pmatrix}
			x& \\
			& x^{-1}
		\end{pmatrix} \in M_{\psi} \subset \GSO(V).
		$$
		For $\Phi \in \Omega(1)$, one has
		\begin{align*}
			f\left(s(x^2)\right)=&\int_{\overline{T}_{E}} \mathrm{pr}_1\left(t\cdot h(x)\cdot \Phi \right)(T_1)\dd t \\
			=& \int_{\overline{T}_{E}} \left|\det(t) \right|^2 \cdot \mathrm{pr}_1\left( m(1,\det(t)) \cdot \Phi \right)(T_1\cdot t\cdot h(x))\dd t \\
			=&\sum_{a \in F^{\times}\slash F^{\times 2}}\int_{\overline{T}_{E},\det(t)\in a\cdot F^{\times 2}} \mathrm{pr}_1\left(m(1,\widetilde{a})\cdot \Phi \right)(T_1\cdot t^{\ast}\cdot h(x)) \dd t 
		\end{align*}
		where $t \in T_{E} \subset  M_{\psi}$ and when $\det(t) \in a \cdot F^{\times 2}$,
		$$
		t^{\ast}:=b \cdot t 
		$$
		with $b^{-2}=\widetilde{a}^{-1}\cdot \det(t)$. Note that $\det(t^{\ast})=\widetilde{a}$. There is a positive constant $C$ such that for all entries $a_{ij}$ of $t^{\ast}$ where $t \in T_{E}$, one has $\left|a_{ij}\right|<C$, and there exists an entry $a_0 \in \{a_{21},a_{22}\}$ satisfying $\left|a_0\right|>C^{-1}$ by the compactness of $\overline{T}_{E}$ and $\det(t^{\ast})=\widetilde{a}$. When $\left|x\right|\rightarrow 0$, by the rapidly decreasing property of $\mathrm{pr}_1\left(m(1,\widetilde{a})\cdot \Phi \right) \in \omega_{\psi}$, one has
		$$
		\mathrm{pr}_1\left(m(1,\widetilde{a})\cdot \Phi \right)(T_1\cdot t^{\ast}\cdot h(x)) \leqslant C_r\cdot x^{r},
		$$
		where $r$ is a big positive integer and $C_r>C$ is some positive constant depending on $r$ but not depending on $t$. When $\left|x\right|\rightarrow \infty$, by the rapidly decreasing property of $\mathrm{pr}_1\left(m(1,\widetilde{a})\cdot \Phi \right) \in \omega_{\psi}$, one has
		$$
		\mathrm{pr}_1\left(m(1,\widetilde{a})\cdot \Phi \right)(T_1\cdot t^{\ast}\cdot h(x)) \leqslant C^{\prime}_{r^{\prime}}\cdot x^{-r^{\prime}},
		$$
		where $r^{\prime}$ is a big positive integer and $C_{r^{\prime}}<C^{-1}$ is some positive constant depending on $r^{\prime}$ but not depending on $t$.

		We have 
		$$
		\overline{\GL_2^{(1)}}=\PGL_2=N\cdot S\cdot K.
		$$
		We consider the growth of 
		\begin{align*}
			f\left(s(x^2) \right)=&\int_{\overline{T}_{E}}  \mathrm{pr}_1\left(t\cdot  m(x^2\cdot I_2, x^2 ) \cdot \Phi \right)(T_1)\dd t \\
			&\int_{\overline{T}_{E}}  \mathrm{pr}_1\left(t\cdot  m(x\cdot I_2, 1 ) \cdot \Phi \right)(T_1)\dd t \\
			=& \int_{\overline{T}_{E}} \left|x^2 \cdot \det(t) \right|^2 \cdot \mathrm{pr}_1\left( m(1,\det(t)) \cdot \Phi \right)(x\cdot T_1\cdot t)\dd t \\
			=&\sum_{a \in F^{\times}\slash F^{\times 2}}\int_{\overline{T}_{E},\det(t)\in a\cdot F^{\times 2}}\left|x\right|^4\cdot \mathrm{pr}_1\left(m(1,\widetilde{a})\cdot \Phi \right)(x\cdot T_1\cdot t^{\ast}) \dd t  
		\end{align*}
		here $s(x^2)\in \GL_2^{(1)} \subset \GSO(V)$ and $\Phi \in \Omega(1)$. When $\left|x\right|\rightarrow \infty$, 
		$$
		\left|x\right|^4\cdot \mathrm{pr}_1\left(m(1,\widetilde{a})\cdot \Phi \right)(x\cdot T_1\cdot t^{\ast}) \leqslant C_{r^{\prime \prime}}^{\prime \prime}\cdot x^{-r^{\prime\prime}},
		$$
		where $r^{\prime\prime}$ is a big enough positive integer and $C_{r^{\prime \prime}}^{\prime \prime}<C^{-1}$ is some positive constant depending on $r^{\prime\prime}$ but not depending on $t$. When $\left|x\right|\rightarrow 0$, there are positive constants $D,D^{\prime}$ such that
		$$
		\left|f(s(x^2))\right|\leqslant D\cdot \left|x\right|^4 \leqslant D^{\prime}\cdot \Xi^{\PGL_2}(s(x^2)).
		$$

		Thus we can conclude that 
		$$
		f \in \mathcal{C}(X).
		$$



	\end{proof}

	\subsection{\bf Plancherel formula of $L^2\left(\SU(2)\backslash \PGSp(W)\right)$}

	From the proof of the Plancherel decomposition of $L^2\left(\SU(2) \backslash \PGSp(W)\right) $,  we know
	\begin{prop} \label{commutativediagram: 1}
		For $\mu_{\PGSO(V)}$-almost all $\sigma\cong \sigma_1\widehat{\boxtimes}\sigma_2 \in \widehat{\PGSO}(V)^{\temp}=\widehat{\PGL}_2^{\temp}\boxtimes \widehat{\PGL}_2^{\temp}$, there is a commutative diagram up to $\mathbb{S}^{1}$:
		\[ \begin{tikzcd}
			\Omega(1) \arrow{rrr}{\theta_{\sigma}} \arrow[swap]{d}{q_1} &&& \sigma\widehat{\boxtimes} \theta(\sigma) \arrow{d}{\ell^{(N,\psi)\backslash \PGL_2}_{\sigma_1}} \\%
			S\left(\SU(2)\backslash \PGSp(W)\right) \arrow{rrr}{\alpha_{\sigma}^{\SU(2)\backslash \PGSp(W)}}&&&  \sigma_2\widehat{\boxtimes} \theta(\sigma)
		\end{tikzcd}
		\]
		where $\theta_{\sigma}$, $\alpha^{\SU(2)\backslash \PGSp(W)}_{\sigma}$ and $\ell^{(N,\psi)\backslash \PGL_2}_{\sigma_1}$ are determined by the Plancherel decompositions of $\Omega(1)$, $L^2\left(\SU(2)\backslash \PGSp(W)\right)$ and $L^2\left(N\backslash \PGL_2,\psi\right)$ separately (see Subsection \ref{subsection: bernsteinmap}).
	\end{prop}

	\begin{proof}
		From the action of the Weil representation in Subsection \ref{subsection2.2b}, one has
		$$
		Z_W\cdot g\cdot g_a \cdot \oplus_{c}\Phi_c(T)=\oplus_c \left|\widetilde{c}\cdot \widetilde{ac}\right|^{2}\cdot \Phi_{ac}(d\cdot g^{-1}\cdot g_a^{-1}\cdot T)
		$$
		where $d$ is such that the following equality holds
		$$
		\widetilde{c}\cdot \widetilde{ac}=d^{2}\cdot \widetilde{a}^{-1}.
		$$
		Then we have
		\begin{align*}
			q_1(\Phi)(Z_W\cdot g\cdot g_a )=&\mathrm{pr}_1\left(Z_W\cdot g\cdot g_a \cdot \Phi \right)(T_1)\\
			=&\varrho(\mathrm{Rest}|_{^{1}Y}\Phi)(Z_W\cdot g\cdot g_a),
		\end{align*}
		where the definition of $\varrho$ is in Lemma \ref{lem4.4}.
		
		From the proof of the Plancherel decomposition of $L^2\left(\SU(2)\backslash \PGSp(W) \right)$ in Subsection \ref{subsection4.1b}, one knows that 
		\begin{align*}
			\Phi \in \Omega(1)\cong &L^2\left( N\backslash \overline{B} ,\psi \right)\widehat{\bigotimes}\int_{\widehat{\PGSO}(V)^{\temp}}\sigma_2\widehat{\boxtimes} \theta(\sigma_1\widehat{\boxtimes} \sigma_2)\dd\mu_{\PGSO(V)}(\sigma_1\widehat{\boxtimes} \sigma_2) 
		\end{align*}
		is represented by 
		$$
		\sigma \mapsto \rest|_{B}\circ\beta_{\sigma_1}^{(N,\psi)\backslash \PGL_2}\circ \theta_{\sigma}(\Phi)
		$$
		and
		$$
		\sigma \mapsto \alpha_{\sigma}^{\SU(2)\backslash \PGSp(W)}\circ\varrho\circ  \varrho^{\prime\prime}(\Phi) 
		$$
		where the definition of $\varrho^{\prime\prime}$ is in Lemma \ref{lemma4.2}.
		Thus one can conclude that up to $\mathbb{S}^{1}$, the following commutative diagram
		\[ \begin{tikzcd}
			\Omega(1) \arrow{rrr}{\theta_{\sigma}} \arrow[swap]{d}{\varrho\circ \mathrm{Rest}|_{^{1}Y}} &&& \sigma\widehat{\boxtimes} \theta(\sigma) \arrow{d}{\ell^{(N,\psi)\backslash \PGL_2}_{\sigma_1}} \\%
			S\left(\SU(2)\backslash \PGSp(W)\right) \arrow{rrr}{\alpha_{\sigma}^{\SU(2)\backslash \PGSp(W)}}&&&  \sigma_2\widehat{\boxtimes} \theta(\sigma)
		\end{tikzcd}
		\]
		holds. Then this proposition follows.
	\end{proof}
	We will fix $\theta_{\sigma}$ such that the commutative diagram
	\[ \begin{tikzcd}
		\Omega(1) \arrow{rrr}{\theta_{\sigma}} \arrow[swap]{d}{\varrho\circ \mathrm{Rest}|_{^{1}Y}} &&& \sigma\widehat{\boxtimes} \theta(\sigma) \arrow{d}{\ell^{(N,\psi)\backslash \PGL_2}_{\sigma_1}} \\%
		S\left(\SU(2)\backslash \PGSp(W)\right) \arrow{rrr}{\alpha_{\sigma}^{\SU(2)\backslash \PGSp(W)}}&&&  \sigma_2\widehat{\boxtimes} \theta(\sigma)
	\end{tikzcd}
	\]
	holds for $\mu_{\PGSO(V)}$-almost all $\sigma\cong \sigma_1\widehat{\boxtimes}\sigma_2 \in \widehat{\PGSO}(V)^{\temp}$.

	We have fixed the measures $\dd T$ on $\sqcup_{a}V_2^{(a)}\otimes W^{(a)}$, $\dd n$ on $N$, and $\frac{\dd g}{\dd s}$ on $\SU(2)\backslash \PGSp(W)$ in Subsection \ref{subsection4.1.a} and Subsection \ref{subsection4.1b}. From Proposition \ref{commutativediagram: 1}, we may conclude the spectral formula of the inner product of $L^2\left(\SU(2) \backslash \PGSp(W)\right) $. 
	\begin{prop} \label{zf5: plancherelformula1}
		For $\Phi_1,\Phi_2 \in \Omega(1)^{\infty}$ and $\sigma\cong \sigma_1\widehat{\boxtimes}\sigma_2 \in \widehat{\PGSO}(V)^{\temp}=\widehat{\PGL}_2^{\temp}\boxtimes \widehat{\PGL}_2^{\temp}$, one has
		$$
		\langle q_1(\Phi_1), q_1(\Phi_2) \rangle_{\SU(2)\backslash \PGSp(W)}=\int_{\widehat{\PGSO}(V)^{\temp}}J^{\SU(2)\backslash \PGSp(W)}_{\sigma}\left(q_1(\Phi_1),q_1(\Phi_2)\right) \dd \mu_{\PGSO(V)}(\sigma)
		$$
		where the positive semi-definite Hermitian form $J^{\SU(2)\backslash \PGSp(W)}_{\sigma}$ is defined as
		$$
		J^{\SU(2)\backslash \PGSp(W)}_{\sigma}\left(q_1(\Phi_1),q_1(\Phi_2)\right)=\int_{N}^{\ast}\overline{\psi(n)}\cdot J^{\theta}_{\sigma}\left(n\cdot\Phi_1,\Phi_2\right)\dd n.
		$$
	\end{prop}
	\begin{proof}
		\begin{align}
			&J^{\SU(2)\backslash \PGSp(W)}_{\sigma}\left(q_1(\Phi_1),q_1(\Phi_2)\right) \notag \\
			=& \langle \alpha^{\SU(2)\backslash \PGSp(W)}_{\sigma}(q_1(\Phi_1)), \alpha^{\SU(2)\backslash \PGSp(W)}_{\sigma}(q_1(\Phi_2))\rangle_{\sigma_2\widehat{\boxtimes}\theta(\sigma)} \notag  \\
			=&\langle \ell^{(N,\psi)\backslash \PGL_2}_{\sigma_1}(\theta_{\sigma}(\Phi_1)),\ell^{(N,\psi)\backslash \PGL_2}_{\sigma_1}(\theta_{\sigma}(\Phi_2))\rangle_{\sigma_2\widehat{\boxtimes}\theta(\sigma)} \notag \\
			=&\ell^{(N,\psi)\backslash \PGL_2}_{\sigma_1}\otimes\overline{\ell^{(N,\psi)\backslash \PGL_2}_{\sigma_1}}\left(\langle \theta_{\sigma}(\Phi_1),\theta_{\sigma}(\Phi_2) \rangle_{\sigma_2\widehat{\boxtimes}\theta(\sigma)}\right) \notag \\
			=& \int_{N}^{\ast}\overline{\psi(n)}\cdot\langle n\cdot\theta_{\sigma}(\Phi_1),\theta_{\sigma}(\Phi_2) \rangle_{\sigma\boxtimes\theta(\sigma)} \dd n \notag \\
			=& \int_{N}^{\ast} \overline{\psi(n)}\cdot  J^{\theta}_{\sigma}(n\cdot \Phi_1,\Phi_2)   \dd n. \notag
		\end{align}
	\end{proof}
	
	Note that the positive semi-definite Hermitian form $J_{\sigma}^{\theta}$ is determined by the Plancherel formula of $\Omega(1)$ (see Section \ref{zf13: 24}).

	\subsection{\bf Plancherel formula of $L^2\left(Y\right)$}
	
	Suppose $F$ is Archimedean. Let $\sigma$ be a unitary tempered representation of $G$. Then we denote the space of continuous endomorphisms of the space of $\sigma$ by $\mathrm{End}(\sigma)$. This is a Banach space with the operator-norm 
	$$
	\normmm{T}=\sup_{\|e\|=1}\|Te \|,\ \text{for}\ T \in \mathrm{End}(\sigma).
	$$
	$\mathrm{End}(\sigma)$ is a continuous representation of $G\times G$ with the action given by left and right translations.

	Fix a $G$-invariant nondegenerate bilinear form $B$ on $\mathfrak{g}$ such that for every maximal compact subgroup $K$ of $G$, the restriction of $B$ to $\mathfrak{k}$ is negative definite and the restriction to $\mathfrak{k}^{\perp}$ (the orthogonal of $\mathfrak{k}$ with respect to $B$) is positive definite. Now we fix a maximal compact subgroup $K$ of $G$. Then we choose a basis $X_1,...,X_k$ of $\mathfrak{k}$ such that 
	$$
	B(X_i,X_j)=-\delta_{ij} \ \text{for $i,j=1,2,...,k$,}
	$$
	and set 
	$$
	\Delta_{K}=1-X_1^{2}-...-X_{k}^{2} \in \mathcal{U}(\mathfrak{k}).
	$$
	Then $\Delta_{K}$ lies in the center of $\mathcal{U}(\mathfrak{k})$ and does not depend on the basis we choose. 
	
	We denote the subspace of smooth vectors by $\mathrm{End}(\sigma)^{\infty}$ and we equip it with the topology defined by the semi-norms
	$$
	\normmm{T}_{\boldsymbol{u},\boldsymbol{v}}=\normmm{\sigma(\boldsymbol{u})T\sigma(\boldsymbol{v})},\ \text{for}\ \boldsymbol{u},\boldsymbol{v} \in \mathcal{U}(\mathfrak{g}), T \in \mathrm{End}(\sigma)^{\infty}.
	$$
	The topology on $\mathrm{End}(\sigma)^{\infty}$ can be generated by the semi-norms $(\normmm{\cdot}_{\Delta_{K}^{r},\Delta_{K}^{r}})_{r\geqslant 1}$ and $\normmm{\cdot}_{\Delta_{K}^{r},\Delta_{K}^{r}} \leqslant \normmm{\cdot}_{\Delta_{K}^{r^{\prime}},\Delta_{K}^{r^{\prime}}}$ for $r \leqslant r^{\prime}$. When $\varphi\in \mathcal{C}(G)$, one has $\sigma(\varphi) \in \mathrm{End}(\sigma)^{\infty}$ (see \cite[\S 2.2, Pg. 48]{BP1}). 
	\begin{lem} \label{estimate: lemma1}
		Assume $F$ is Archimedean. For $\varphi \in \mathcal{C}\left(\PGSO(V)\right)$, there exists a finite set $\{\varphi_i \in \mathcal{C}(\PGSO(V))\}$ such that 
		$$
		\left|\int_{N}^{\ast}\Tr\sigma((n,g_2))\sigma(\varphi)\cdot \overline{\psi(n)}\dd n\right| \leqslant C \cdot \sum_i \normmm{\sigma(\varphi_i)}_{\Delta_{K}^{r},\Delta_{K}^{r}} \cdot \Xi^{\PGL_2}(g_2)
		$$
		for any irreducible tempered representation $\sigma \in \widehat{\PGSO}(V)^{\temp}$, where $C$ is some positive constant not depending on $\sigma$, $g_2 \in \overline{M}_{\psi}\cong \PGL_2$ and $r$ is sufficiently large.
	\end{lem}
	
	\begin{proof}
		Note that $\Tr\sigma( (\cdot,g_2))\sigma(\varphi) \in \mathcal{C}_{d}^{\omega}\left(\PGL_2\right)$ for any $d>0$ \cite[(2.2.5), Pg.49]{BP1}, where $\mathcal{C}_{d}^{\omega}\left(\PGL_2\right)$ is a Fr\'echet space \cite[Pg.33]{BP1}. 

		By the continuity of $\int_{N}^{\ast}\cdot \dd n$, there exist finite sets $S_1:=\{\boldsymbol{u}_i\}\subset \mathfrak{g}$ and $S_2:=\{\boldsymbol{v}_i\}\subset \mathfrak{g}$ such that
		\begin{align*}
			&\left|\int_{N}^{\ast}\Tr\sigma((n, g_2))\sigma(\varphi)\cdot\overline{\psi(n)}\dd n\right| \\
			\leqslant& C\cdot  \sum_i \sup_{g_1\in \PGL_2}\left|\Tr\sigma((g_1, g_2))\sigma(L(\boldsymbol{u}_i)R(\boldsymbol{v}_i)\varphi) \right|\Xi^{\PGL_2}(g_1)^{-1} \\ 
			\leqslant& C\cdot \sum_i\sup_{g_1\in \PGL_2}\Xi^{\PGSO(V)}((g_1,g_2))\normmm{\sigma(\varphi_i)}_{\Delta_{K}^{r},\Delta_{K}^{r}}\Xi^{\PGL_2}(g_1)^{-1} \quad \text{by \cite[(2.2.5)]{BP1}}\\ 
			\leqslant& C \cdot\sum_i \normmm{\sigma(\varphi_i)}_{\Delta_{K}^{r},\Delta_{K}^{r}} \cdot \Xi^{\PGL_2}(g_2),
		\end{align*}
		where $\varphi_i:=L(\boldsymbol{u}_i)R(\boldsymbol{v}_i)\varphi$, and $r$ is sufficiently large.
	\end{proof}
	Suppose $F$ is non-Archimedean. We have an estimate similar to Lemma \ref{estimate: lemma1}.
	\begin{lem}
		\label{estimate: lemma3}
		When $F$ is non-Archimedean, for any irreducible tempered representation $\sigma \in \widehat{\PGSO}(V)^{\temp}$,  and any $\varphi \in \mathcal{C}\left(\PGSO(V)\right)$, one has
		$$
		\left|\int_{N}^{\ast}\Tr\sigma((n, g_2))\sigma(\varphi)\cdot \overline{\psi(n)}\dd n\right| \leqslant C_{\varphi}  \cdot \Xi^{\PGL_2}(g_2)
		$$
		where $C_{\varphi}$ is some positive constant, $g_2 \in \overline{M}_{\psi}\cong \PGL_2$.
	\end{lem}
	\begin{proof}
		Note that $\Tr\sigma((\cdot, g_2))\sigma(\varphi) \in \mathcal{C}_{d,K^{\prime}}^{\omega}\left(\PGL_2\right)$ for any $d>0$, and $\mathcal{C}_{d,K^{\prime}}^{\omega}\left(\PGL_2\right)$ is a Banach space with the norm $p_{-d}$ \cite[Pg.32]{BP1}, where $K^{\prime}$ is some open-compact subgroup of $K^{\prime}$. 
		Following from \cite[Theorem 2]{CHH}, we have
		$$
		\left|\Tr\sigma((g_1, g_2))\sigma(\varphi)\right| \leqslant C_{\varphi}\cdot\Xi^{G}((g_1,g_2)).
		$$
		By the continuity of the integral $\int_{N}^{\ast}$, one has	
		\begin{align*}
			\left|\int_{N}^{\ast}\Tr\sigma((n, g_2))\sigma(\varphi)\cdot \overline{\psi(n)}\dd n\right| \leqslant& C  \cdot p_{-d}\left(\Tr\sigma((\cdot, g_2))\sigma(\varphi) \right) \\
			=&C  \cdot \sup_{g_1\in \PGL_2} \left|\Tr\sigma((g_1, g_2))\sigma(\varphi)\right| \cdot \Xi^{\PGL_2}(g_1)^{-1}\cdot \sigma_{\PGL_2}(g_1)^{-d}  \\
			\leqslant& C_{\varphi}\cdot\Xi^{\PGL_2}(g_2)\cdot \sup_{g_1\in \PGL_2}\sigma_{\PGL_2}(g_1)^{-d} \\
			\leqslant& C_{\varphi}\cdot\Xi^{\PGL_2}(g_2).
		\end{align*}
	\end{proof}
	
	We take note of the following lemma which will be used in the following section.
	\begin{lem} \label{estimate: lemma2}
		For $\Phi_1, \Phi_2 \in \Omega(1)^{\infty}$, one has
		$$
		J^{\theta}_{\sigma}\left(h\cdot \Phi_1,\Phi_2 \right)=\Tr\sigma(h)\sigma(\varphi^{\vee})
		$$
		where $\varphi(h^{\prime}):=\langle h^{\prime}\cdot \Phi_1,\Phi_2  \rangle_{\Omega(1)}\in \mathcal{C}\left(\PGSO(V)\right)$, and $\varphi^{\vee}(h^{\prime})=\varphi(h^{\prime-1})$.
	\end{lem}
	
	\begin{proof}
		\begin{align*}
			J_{\sigma}^{\theta}\left(h\cdot \Phi_1,\Phi_2 \right)=& \sum_{v\in \ONB(\sigma)}\int_{\PGSO(V)}\langle h^{\prime}h\cdot \Phi_1,\Phi_2  \rangle_{\Omega(1)}\cdot \overline{\langle h^{\prime}\cdot v,v \rangle}_{\sigma}\dd h^{\prime} \\
			=& \Tr\sigma(h)\sigma(\varphi^{\vee}).
		\end{align*}
	\end{proof}


	We obtain the Plancherel formula of $L^2(Y)$:
	\begin{thm} \label{thm: plancherel}
		For $\Phi_1,\Phi_2 \in \Omega(1)^{\infty}$ with $\varphi_i=q(\Phi_i)$, we have
		$$
		\langle\varphi_1,\varphi_2 \rangle_{L^2(Y)}=\int_{\widehat{\PGSO}(V)^{\temp}}\int_{\overline{T}_{E}}\int_{N}^{\ast}\overline{\psi(n)}\cdot J^{\theta}_{\sigma}\left((n, t)\cdot \Phi_1,\Phi_2\right)\dd n \dd t \dd \mu_{\PGSO(V)}(\sigma),
		$$
		where the positive semi-definite Hermitian forms $J^{\theta}_{\sigma}(\cdot,\cdot)$ are determined by the Plancherel decomposition of the Weil representation $\Omega(1)$ (see Section \ref{zf13: 24}), $(n,t) \in N\times \overline{T}_{E}\subset  \overline{\GL}_2^{(1)}\times \overline{M}_{\psi}\cong \PGSO(V)$ (see Subsection \ref{zf13:21a}), the measure on $Y$ is determined by the measure $\frac{\dd g}{\dd s}$ on $\SU(2)\backslash \PGSp(W)$ (see Subsection \ref{subsection4.1b}) and a Haar measure on $\overline{T}_{E}$, the measure $\dd n$ on $N(F)=F$ is self-dual with respect to the pairing $\psi(xy)$ on $F$, and the Plancherel measure $\mu_{\PGSO(V)}$ is determined by a Haar measure $\dd h$ on $\PGSO(V)$.
	\end{thm}
	\begin{proof}
		
		For any $\Phi_1, \Phi_2 \in \Omega(1)^{\infty}$, $q_1(\Phi_1), q_1(\Phi_2) \in S(\SU(2)\backslash \PGSp(W))$, by Proposition \ref{zf5: plancherelformula1}, one has
		$$
		\langle q_1(\Phi_1), q_1(\Phi_2)  \rangle_{\SU(2)\backslash \PGSp(W)}= \int_{\widehat{\PGSO}(V)^{\temp}}\int_{N}^{\ast}J_{\sigma}^{\theta}(n\cdot\Phi_1, \Phi_2)\cdot \overline{\psi(n)} \ \dd n \dd \mu_{\PGSO(V)}(\sigma).    
		$$
		Then one can deduce
		\begin{align}
			&\langle q(\Phi_1),  q(\Phi_2)  \rangle_{\U(2)\backslash \PGSp(W)}\notag \\
			= & \int_{\U(2)\backslash \PGSp(W)}\int_{\overline{T}_E}L(t)q_1\left(\Phi_1\right)(g)\ \dd t \int_{\overline{T}_E} L(t^{\prime})\overline{q_1\left(\Phi_2\right)(g)}\ \dd t^{\prime}\dd g  \notag  \\
			= & \int_{\SU(2)\backslash \PGSp(W)} \int_{\overline{T}_E}L(t)q_1\left(\Phi_1\right)(g)\overline{q_1\left(\Phi_2\right)(g)}\ \dd t \dd g    \notag \\
			= &  \int_{\overline{T}_E} \int_{\SU(2)\backslash \PGSp(W)}L(t)q_1\left(\Phi_1\right)(g)\overline{q_1\left(\Phi_2\right)(g)}\ \dd g \dd t     \notag\quad \text{($\left|q_1(\Phi_1)\right|,\left|q_1(\Phi_2)\right| \in S\left(\SU(2)\backslash \PGSp(W)\right)$)}\\       
			= & \int_{\overline{T}_E}  \langle L(t)q_1(\Phi_1), q_1(\Phi_2)  \rangle_{\SU(2)\backslash \PGSp(W)}\  \dd t  \notag \\
			= & \int_{\overline{T}_E}  \int_{\widehat{\PGSO}(V)^{\temp}}\int_{N}^{\ast}J_{\sigma}^{\theta}((n,t)\cdot\Phi_1, \Phi_2) \overline{\psi(n)}\  \dd n\dd \mu_{\PGSO(V)}(\sigma) \dd t  \notag \\
			=& \int_{\overline{T}_E}\int_{\widehat{\PGSO}(V)^{\temp}}\int_{N}^{\ast}\Tr\left(\sigma((n,t))\sigma(f_{\Phi_1,\Phi_2}^{\vee})   \right)\cdot \overline{\psi(n)}\ \dd n \dd \mu_{\PGSO(V)}(\sigma)\dd t \qquad \text{(Lemma \ref{estimate: lemma2})},    \notag 
		\end{align}
		where 
		$$
		f_{\Phi_1,\Phi_2}(g)=\langle g\cdot  \Phi_1, \Phi_2   \rangle_{\Omega(1)}  \in \mathcal{C}\left(\PGSO(V)\right).
		$$

		When $F$ is Archimedean, by Lemma \ref{estimate: lemma1}, we have
		$$
		\iint_{\overline{T}_{E} \times \widehat{\PGSO}(V)^{\temp}}\left|\int_{N}^{\ast}\Tr\left(\sigma((n,t))\sigma(f_{\Phi_1,\Phi_2}^{\vee})   \right)\cdot \overline{\psi(n)}\dd n \right| \dd \mu_{\PGSO(V)}(\sigma)\dd t <\infty.
		$$	
		

		When $F$ is non-Archimedean, by \cite[Th\'eor\`eme VIII.1.2]{W}, one knows
		$$
		\overline{\{\sigma \in \widehat{G}^{\temp}_{\ind}\mid \Tr\sigma(g)\sigma(f_{\Phi_1,\Phi_2}^{\vee})\neq 0 \}}
		$$
		is a compact set in $\widehat{G}^{\temp}_{\ind}$.
		Also following from Lemma \ref{estimate: lemma3}, we have
		$$
		\left|\int_{N}^{\ast}\Tr\sigma((n,t))\sigma(f_{\Phi_1,\Phi_2}^{\vee})\cdot \overline{\psi(n)}\ \dd n\right| \leqslant C_{\Phi_1,\Phi_2} \cdot \Xi^{\PGL_2}(t).
		$$
		Then we have 
		$$
		\iint_{\overline{T}_{E} \times \widehat{\PGSO}(V)^{\temp}}\left|\int_{N}^{\ast}\Tr\left(\sigma((n,t))\sigma(f_{\Phi_1,\Phi_2}^{\vee})   \right)\cdot \overline{\psi(n)}\dd n \right| \dd \mu_{\PGSO(V)}(\sigma)\dd t <\infty.
		$$

		
		
		Thus one has
		\begin{align*}
			&\langle q(\Phi_1),  q(\Phi_2)  \rangle_{\U(2)\backslash \PGSp(W)}\notag \\
			=&\iint_{\overline{T}_{E} \times \widehat{\PGSO}(V)^{\temp}}\int_{N}^{\ast}\Tr\left(\sigma((n,t))\sigma(f_{\Phi_1,\Phi_2}^{\vee})   \right)\cdot \overline{\psi(n)}\ \dd n  \dd \mu_{\PGSO(V)}(\sigma)\dd t \\
			=&\int_{ \widehat{\PGSO}(V)^{\temp}}\int_{\overline{T}_{E}}\int_{N}^{\ast}\Tr\left(\sigma((n,t))\sigma(f_{\Phi_1,\Phi_2}^{\vee})   \right)\cdot \overline{\psi(n)}\ \dd n \dd t  \dd \mu_{\PGSO(V)}(\sigma) \\
			=&\int_{\widehat{\PGSO}(V)^{\temp}}\int_{\overline{T}_{E}}\int_{N}^{\ast}\overline{\psi(n)}\cdot J^{\theta}_{\sigma}\left((n,t)\cdot \Phi_1,\Phi_2\right)\ \dd n \dd t \dd \mu_{\PGSO(V)}(\sigma).
		\end{align*}
	\end{proof}
	Note that the proof holds for $\overline{T}_{E}\cong F^{\times}$ so that $Y=\GL_2\backslash \PGSp(W)$.

	By Theorem \ref{thm: plancherel}, one has
	\begin{equation} \label{zf2c: 1}
		L^2\left(Y\right)\cong\int_{\widehat{\PGSO}(V)^{\temp}} \mathcal{H}_{\sigma} \dd \mu_{\PGSO(V)}(\sigma),
	\end{equation}
	where $\mathcal{H}_{\sigma} $ is some Hilbert space (may be $0$). Let $\alpha_{\sigma}^{Y}$, $\beta_{\sigma}^{Y}$ and $\ell_{\sigma}^{Y}$ be determined by the isomorphism (\ref{zf2c: 1}) (see Subsection \ref{subsection: bernsteinmap}). Then by the Plancherel formula of $L^2(Y)$, we know that for $\mu_{\PGSO(V)}$-almost all $\sigma \in \widehat{\PGSO}(V)^{\temp}$, one has
	\begin{align}
		&\langle \alpha_{\sigma}^{Y}\circ q(\Phi_1),\alpha_{\sigma}^{Y}\circ q(\Phi_2)  \rangle_{\mathcal{H}_{\sigma}} \notag \\
		=&\int_{\overline{T}_{E}}\int_{N}^{\ast}\overline{\psi(n)}\cdot J^{\theta}_{\sigma}\left((n, t)\cdot \Phi_1,\Phi_2\right)\dd n \dd t  \notag \\
		=&\ell^{X}_{\sigma}\otimes\overline{\ell^{X}_{\sigma}}\left(\langle \theta_{\sigma}(\Phi_1), \theta_{\sigma}(\Phi_2) \rangle_{\theta(\sigma)}\right) \notag \\
		=&\langle \ell^{X}_{\sigma}\circ \theta_{\sigma}(\Phi_1),\ell^{X}_{\sigma}\circ \theta_{\sigma}(\Phi_2) \rangle_{\theta(\sigma)} \label{plancherelformula: eq100}
	\end{align}
	We need to determine what exactly $\mathcal{H}_{\sigma}$ is.
	\begin{lem} \label{pldecomposition: lemma1}
		For almost every $\sigma$, one has
		$$
		\mathcal{H}_{\sigma}\neq 0 \ \text{if and only if}\  \Hom_{N\times \overline{T}_{E}}\left(\sigma,\psi\boxtimes 1\right) \neq 0 
		$$
	\end{lem}
	\begin{proof}
		On one hand, by equation (\ref{plancherelformula: eq100}), one has
		\begin{equation*}
			\ell_{\sigma}^{X} \neq 0 \ \text{if and only if}\  \mathcal{H}_{\sigma}\neq 0
		\end{equation*}
		On the other hand, by the Plancherel decompositions of $L^2\left(N\backslash \PGL_2,\psi \right)$ and $L^2\left(\overline{T}_{E}\backslash \PGL_2 \right)$, and by equation (\ref{plancherelformula: eq100}), one has
		\begin{equation*}
			\ell_{\sigma}^{X} \neq 0 \ \text{if and only if}\  \Hom_{N\times \overline{T}_{E}}\left(\sigma,\psi\boxtimes 1\right) \neq 0 
		\end{equation*}
	\end{proof}
	
	\begin{lem} \label{plancherel: lemma200}
		If $\mathcal{H}_{\sigma} \neq 0$, then we have 
		$$
		\mathcal{H}_{\sigma}\cong \theta(\sigma).
		$$
	\end{lem}
	\begin{proof}
		If $\mathcal{H}_{\sigma} \neq 0$, then $\ell_{\sigma}^{X}\neq 0$, and by diagram (\ref{commutativediagram: 2}), one has $ \alpha_{\sigma}^{Y}\neq 0$.
		
		First take an orthonormal basis $\{\alpha_{\sigma}^{Y}\circ q(\Phi_i)\mid i\in I \}$ of $\mathcal{H}_{\sigma}$. Then by equation (\ref{plancherelformula: eq100}), one knows $\{\ell^{X}_{\sigma}\circ \theta_{\sigma}(\Phi_i)\mid i \in I \}$ is an orthonormal basis of $\theta(\sigma)$. 
		
		Second, if $\Phi,\Phi^{\prime} \in \Omega(1)^{\infty}$ with $\alpha_{\sigma}^{Y}\circ q(\Phi)=\alpha_{\sigma}^{Y}\circ q(\Phi^{\prime})$, then one has $\ell^{X}_{\sigma}\circ \theta_{\sigma}(\Phi)=\ell^{X}_{\sigma}\circ \theta_{\sigma}(\Phi^{\prime})$.
		
		Thus the map 
		$$
		\mathcal{H}_{\sigma}\longrightarrow \theta(\sigma)
		$$
		given by
		$$
		\alpha_{\sigma}^{Y}\circ q(\Phi)\longmapsto \ell^{X}_{\sigma}\circ \theta_{\sigma}(\Phi)
		$$
		is well-defined and $\PGSp(W)$-isometric.
	\end{proof}


	

	By Lemma \ref{pldecomposition: lemma1}, Lemma \ref{plancherel: lemma200} and the isomorphism (\ref{zf2c: 1}), one has the Plancherel decomposition of $L^2(Y)$:
	\begin{thm} \label{Theorem4.21}
		There is a $\PGSp(W)$-isometric isomorphism
		$$
		L^2\left(Y\right)\cong\int_{\widehat{\PGSO}(V)^{\temp}} \Hom_{N\times\overline{T}_{E}}\left(\sigma,\psi\boxtimes 1\right)\cdot \theta(\sigma) \dd \mu_{\PGSO(V)}(\sigma).
		$$	
	\end{thm}

	%
	
	We can obtain:
	\begin{cor} \label{commutativediagram: 2}
		For $\mu_{\PGSO(V)}$-almost all $\sigma \in \widehat{\PGSO}(V)^{\temp}$, there is a commutative diagram :
		\[ \begin{tikzcd}
			\Omega(1) \arrow{r}{\theta_{\sigma}} \arrow[swap]{d}{q} & \sigma\boxtimes \theta(\sigma) \arrow{d}{\ell^{X}_{\sigma}} \\%
			S(Y) \arrow{r}{\alpha^{Y}_{\sigma}}&  \theta(\sigma)
		\end{tikzcd}
		\]
		where $\theta_{\sigma}$, $\alpha^{Y}_{\sigma}$ and $\ell^{X}_{\sigma}$ are determined by the Plancherel decompositions of $\Omega(1)$, $L^2(Y)$ and $L^2(X)$ seperately.
	\end{cor}

	Combining Proposition \ref{commutativediagram: 1} with Corollary \ref{commutativediagram: 2}, one obtains
	\begin{cor} \label{plancherel: cor}
		For $\mu_{\PGSO(V)}$-almost all $\sigma\cong \sigma_1\widehat{\boxtimes}\sigma_2 \in \widehat{\PGSO}(V)^{\temp}=\widehat{\PGL}_2^{\temp}\boxtimes \widehat{\PGL}_2^{\temp} $, there is a commutative diagram:
		\[ \begin{tikzcd}
			S(\SU(2)\backslash \PGSp(W)) \arrow{rr}{\alpha^{\SU(2)\backslash \PGSp(W)}_{\sigma}} \arrow[swap]{d}{q_2} && \sigma_2\widehat{\boxtimes} \theta(\sigma) \arrow{d}{\ell^{\overline{T}_{E}\backslash \overline{M}_{\psi}}_{ \sigma_2}} \\%
			S(Y) \arrow{rr}{\alpha^{Y}_{\sigma}}&&  \theta(\sigma)
		\end{tikzcd}
		\]
	\end{cor}
	
	\begin{cor} \label{cor4.24}
		For $\mu_{\PGSO(V)}$-almost all $\sigma \in \widehat{\PGSO}(V)^{\temp}$, any $\Phi \in \Omega(1)^{\infty}$ and $w \in \theta(\sigma)^{\infty}$, one has
		$$
		\ell^{X}_{\sigma}\left(B_{\sigma}(\Phi,w)  \right)=\langle q(\Phi),\beta_{\sigma}^{Y}(w) \rangle_{Y},
		$$
		where $B_{\sigma}(\Phi,w):=\langle \theta_{\sigma}(\Phi),w \rangle_{\theta(\sigma)}$ (see Subsection \ref{zf13d: 1}).
	\end{cor}
	\begin{proof}
		For any $\Phi \in \Omega(1)^{\infty}$ and $w \in \theta(\sigma)^{\infty}$, one has
		\begin{align}
			\langle q(\Phi),\beta_{\sigma}^{Y}(w) \rangle_{Y}=&\langle \alpha_{\sigma}^{Y}(q(\Phi)), w \rangle_{\theta(\sigma)} \notag \\
			=& \langle \ell^{X}_{\sigma}\left(\theta_{\sigma}(\Phi)\right),w\rangle_{\theta(\sigma)} \notag \\
			=& \ell^{X}_{\sigma}\left( \langle \theta_{\sigma}(\Phi), w \rangle_{\theta(\sigma)}\right) \notag \\
			=&  \ell^{X}_{\sigma}\left( B_{\sigma}(\Phi,w)\right). \notag
		\end{align}
	\end{proof}

\section{\bf Relative character identities}
In this section, we will prove the relative character identities. 


\subsection{\bf Continuity}
In this subsection, we will consider the continuity of relative character.

The following lemma is the pointwise formula version of  \cite[Theorem 6.2.1]{SV}.  
\begin{lem} \label{continuity: lemma1}
	Let $X=H\backslash G$ be strongly tempered. For any $f \in C_c^{\infty}(G)$ with $\varphi(x):=\int_{H}f(h\cdot x)\dd h \in C_c^{\infty}(X)$, one has
	\begin{align*}
		\varphi(g)=&\int_{H}f(h\cdot g)\dd h \\
		=&\int_{H}\int_{\widehat{G}^{\temp}}\Tr\sigma(h\cdot g)\sigma(f^{\vee})\dd \mu_{G}(\sigma)\dd h \\
		=&\int_{\widehat{G}^{\temp}}\int_{H}\Tr\sigma(h\cdot g)\sigma(f^{\vee})\dd h\dd \mu_{G}(\sigma), 
	\end{align*}
	where $f^{\vee}(g):=f(g^{-1})$.
\end{lem}
\begin{proof}
	The last equality holds because 
	\begin{enumerate}
		\item when $F$ is $p$-adic, by \cite[Th\'eor\`em VIII.1.2]{W}, one knows
		$$
		\overline{\{\sigma \in \widehat{G}^{\temp}_{\ind}\mid \Tr\sigma(g)\sigma(f^{\vee})\neq 0 \}}
		$$
		is a compact set in $\widehat{G}^{\temp}_{\ind}$.
		Also following from \cite[Theorem 2]{CHH}, we have
		$$
		\left|\Tr\sigma(h\cdot g)\sigma(f^{\vee})\right| \leqslant C_{f}\cdot\Xi^{G}(h\cdot g),
		$$
		and by \cite[Lemma II.1.2]{W}, we have
		$$
		\Xi^{G}(\cdot g)\in \mathcal{C}^{\omega}(G).
		$$
		Thus 
		$$
		\int_{\widehat{G}^{\temp}}\int_{H}\left|\Tr\sigma(h\cdot g)\sigma(f^{\vee})\right|\dd h\dd \mu_{G}(\sigma)<\infty. 
		$$
		\item when $F$ is Archimedean, by Lemma \ref{estimate: lemma1} and a basic estimate in \cite[(2.6.2)]{BP1}, one has
		$$
		\int_{\widehat{G}^{\temp}}\int_{H}\left|\Tr\sigma(h\cdot g)\sigma(f^{\vee})\right|\dd h\dd \mu_{G}(\sigma)<\infty. 
		$$
	\end{enumerate}
\end{proof}

By \cite[Proposition 2.14.2]{BP2}, one obtains
\begin{lem} \label{continuity: lemma2}
	Let $X=\left(N\backslash G,\psi\right)$ be a Whittaker variety. For any $f \in C_c^{\infty}(G)$ with $\varphi(x):=\int_{N}^{\ast}f(n\cdot x)\cdot \overline{\psi(n)}\dd n\in C_c^{\infty}\left(N\backslash G,\psi\right)$, one has
	$$
	\varphi(g)=\int_{\widehat{G}^{\temp}}\int_{N}^{\ast}\Tr\sigma(n\cdot g)\sigma(f^{\vee})\cdot \overline{\psi(n)}\dd n \dd \mu_{G}(\sigma).
	$$
\end{lem}

By Lemma \ref{continuity: lemma1} and Lemma \ref{continuity: lemma2}, we obtain
\begin{lem}
	\begin{enumerate}
		\item Let $X=H\backslash G$ be strongly tempered. For $\mu_{G}$-almost all $\sigma \in \widehat{G}^{\temp}_{\ind}$, one has
		$$
		\vartheta_{\sigma}^{X}(\varphi)=\int_{H}\Tr\sigma(h)\sigma(f^{\vee})\dd h
		$$
		for any $f \in C_c^{\infty}(G)$ with $\varphi(x):=\int_{H}f(h\cdot x)\dd h \in C_c^{\infty}(X)$.
		\item Let $X=\left(N\backslash G,\psi\right)$ be a Whittaker variety. For $\mu_{G}$-almost all $\sigma \in \widehat{G}^{\temp}_{\ind}$, one has
		$$
		\vartheta_{\sigma}^{X}(\varphi)=\int_{N}^{\ast}\Tr\sigma(n)\sigma(f^{\vee})\cdot \overline{\psi(n)}\dd n
		$$
		for any $f \in C_c^{\infty}(G)$ with $\varphi(x):=\int_{N}^{\ast}f(n\cdot x)\cdot \overline{\psi(n)}\dd n\in C_c^{\infty}\left(N\backslash G,\psi\right)$.
	\end{enumerate}
\end{lem}

By \cite[Lemma 2.3.1]{BP1} and the continuity of the (stable) integral
$$
\int_{H}\cdot \dd h: \mathcal{C}^{\omega}(G)\longrightarrow \mathbb{C}
$$
and
$$
\int_{N}^{\ast}\cdot \overline{\psi(n)}\dd n: \mathcal{C}^{\omega}(G)\longrightarrow \mathbb{C},
$$
one can deduce
\begin{lem} \label{continuity: lemma3}
	\begin{enumerate}
		\item Let $X=H\backslash G$ be strongly tempered. For any $f \in C_c^{\infty}(G)$,
		$$
		\int_{H}\Tr\sigma(h)\sigma(f^{\vee})\dd h
		$$
		is continuous with respect to $\sigma \in \widehat{G}^{\temp}_{\ind}$ in the sense of \cite[Lemma 2.3.1(ii)]{BP1}.
		\item Let $X=(N\backslash G,\psi)$ be a Whittaker variety. For any $f \in C_c^{\infty}(G)$,
		$$
		\int_{N}^{\ast}\Tr\sigma(n)\sigma(f^{\vee})\cdot \overline{\psi(n)}\dd n
		$$
		is continuous with respect to $\sigma \in \widehat{G}^{\temp}_{\ind}$.
	\end{enumerate}
\end{lem}

For the continuity of relative character of $L^2(Y)$:
\begin{lem} \label{continuity: lemma4}
	For $\Phi_1,\Phi_2 \in \Omega(1)^{\infty}$, one knows
	$$
	\int_{\overline{T}_{E}}\int_{N}^{\ast}J^{\theta}_{\sigma}((t, n)\cdot \Phi_1,\Phi_2)\cdot \overline{\psi(n)}\dd n\dd t
	$$
	is continuous with respect to $\sigma \in \widehat{\PGSO}(V)^{\temp}_{\ind}=\widehat{\PGSO}(V)^{\temp}$.
\end{lem}
\begin{proof}
	By Lemma \ref{estimate: lemma2} and \cite[Lemma 2.3.1(i)]{BP1}, one has 
	$$
	J^{\theta}_{\sigma}(h\cdot \Phi_1,\Phi_2)\in \mathcal{C}^{\omega}\left(\PGSO(V)\right).
	$$
	By \cite[Lemma 2.3.1(ii)]{BP1}, one knows
	$J^{\theta}_{\sigma}(h\cdot \Phi_1,\Phi_2)$ is continuous with respect to $\sigma \in \widehat{\PGSO}(V)^{\temp}_{\ind}=\widehat{\PGSO}(V)^{\temp}$.
	By the proof of Theorem \ref{thm: plancherel}, one has 
	$$
	\int_{N}^{\ast}J^{\theta}_{\sigma}((g,n)\cdot \Phi_1,\Phi_2)\cdot \overline{\psi(n)}\dd n \in \mathcal{C}^{\omega}(\PGL_2).
	$$
	By \cite[Lemma 2.14.1]{BP2}, one knows
	$$
	\int_{N}^{\ast}J^{\theta}_{\sigma}((g,n)\cdot \Phi_1,\Phi_2)\cdot \overline{\psi(n)}\dd n
	$$
	is continuous with respect to $\sigma$. 
	Since $\overline{T}_{E}\backslash \PGL_2$ is strongly tempered, the following map
	$$
	\int_{\overline{T}_{E}}\cdot \dd t: \mathcal{C}^{\omega}(\PGL_2)\longrightarrow \mathbb{C}
	$$
	is continuous. Thus
	$$
	\int_{\overline{T}_{E}}\int_{N}^{\ast}J^{\theta}_{\sigma}((t, n)\cdot \Phi_1,\Phi_2)\cdot \overline{\psi(n)}\dd n\dd t
	$$
	is continuous with respect to $\sigma \in \widehat{\PGSO}(V)^{\temp}$.
\end{proof}
For $\varphi \in C_c^{\infty}(X)$ with 
$$
\varphi(x)=\int_{H}f(h\cdot x)\dd h;\qquad \varphi(x)=\int_{N}^{\ast}f(n\cdot x)\cdot \overline{\psi(n)}\dd n
$$
where $f\in C_c^{\infty}(G)$, one has
\begin{equation}\label{continuity: eq101}
	\int_{H}\Tr\sigma(h)\sigma(f^{\vee})\dd h=\ell_{\sigma}^{X}\circ \alpha_{\sigma}^{X}(\varphi); \qquad 		\int_{N}^{\ast}\Tr\sigma(n)\sigma(f^{\vee})\cdot \overline{\psi(n)}\dd n=\ell_{\sigma}^{X}\circ \alpha_{\sigma}^{X}(\varphi)
\end{equation}
for $\mu_{G_X}$-almost all $\sigma$. By Lemma \ref{continuity: lemma3}, these functions are continuous in $\widehat{G}^{\temp}_{\ind}$. For general $\varphi \in \mathcal{C}(X)$, we could not use the expression of relative character in equation \ref{continuity: eq101} anymore. 



\subsection{\bf Relative character identities}
Now we can prove the relative character identities between $L^2(X)$ and $L^2(Y)$.

\begin{thm} \label{relative characteridentity: thm1}
	
	Suppose that 
	\begin{itemize}
		\item $\sigma$ is a tempered irreducible representation of $\PGSO(V)$;
		\item $f_1 \in \mathcal{S}^{\heartsuit}(X)$ and $f_2 \in S(Y)$ are in correspondence, i.e.\ there exists a $\Phi \in \Omega(1)^{\infty}$ such that $p(\Phi)=f_1$ and $q(\Phi)=f_2$;
		\item $\ell^{X}_{\sigma} \in \Hom_{N \times \overline{T}_{E}}\left(\sigma, \psi\boxtimes 1  \right)$ and $\ell^{Y}_{\sigma} \in \Hom_{\U(2)}\left(\theta(\sigma), \mathbb{C}  \right)$ are in correspondence in the sense that $\ell^{X}_{\sigma}$ and $\alpha_{\sigma}^{Y}$ (determined by $\ell^{Y}_{\sigma}$) satisfy the diagram in Corollary \ref{commutativediagram: 2}.
	\end{itemize}
	
	\begin{equation} \label{relative characteridentity: eq1}
		\vartheta_{\sigma}^{X}(f_1)=\vartheta_{\sigma}^{Y}(f_2).
	\end{equation}
\end{thm} 
\begin{proof}
	We proceed the proof in two steps.
	
	$\mathrm{Step\  1}$: We take the Plancherel measure $\mu_{\PGSO(V)}$ of $\PGSO(V)$ as underlying spectral measure of the Plancherel decomposition. Note that for $ \mu_{\PGSO(V)}$-almost all $\sigma \in \widehat{\PGSO}(V)^{\temp}$, there exists a linear functional $l_{\sigma} \in \Hom_{\U(2)}\left(\theta(\sigma), \mathbb{C} \right)$ such that the following diagram is commutative:

	\begin{equation} \label{relativecharacteridentity: diagram1}
		\begin{tikzcd}
			\Omega(1) \arrow{r}{\theta_{\sigma}} \arrow[swap]{d}{p} & \sigma\boxtimes \theta(\sigma) \arrow[dashed]{d}{l_{\sigma}} \\%
			\mathcal{S}^{\heartsuit}(X) \arrow{r}{\alpha^{X}_{\sigma}}& \sigma
		\end{tikzcd}
	\end{equation}
	
	By Proposition \ref{commutativediagram: 2}, we also have the following commutative diagram:
	\begin{equation} \label{relativecharacteridentity: diagram2}
		\begin{tikzcd}
			\Omega(1) \arrow{r}{\theta_{\sigma}} \arrow[swap]{d}{q} & \sigma\boxtimes \theta(\sigma) \arrow{d}{\ell^{X}_{\sigma}} \\%
			S(Y) \arrow{r}{\alpha^{Y}_{\sigma}}& \theta(\sigma)
		\end{tikzcd}
	\end{equation}
	
	Following from the diagram (\ref{relativecharacteridentity: diagram1}) and diagram (\ref{relativecharacteridentity: diagram2}), for $ \mu_{\PGSO(V)}$-almost all $\sigma \in \widehat{\PGSO}(V)^{\temp}$, one has
	$$
	\vartheta_{\sigma}^{X}\circ p-\vartheta_{\sigma}^{Y}\circ q= \ell_{\sigma}^{X}\circ l_{\sigma}\circ \theta_{\sigma}- \ell_{\sigma}^{X}\circ \ell^{Y}_{\sigma}\circ \theta_{\sigma}.
	$$

		$\mathrm{Step \ 2}$: The aim of this step is to demonstrate that for $ \mu_{\PGSO(V)}$-almost all $\sigma \in \widehat{\PGSO}(V)^{\temp}$, one has
		$$
		\vartheta_{\sigma}^{Y}\circ q= \vartheta_{\sigma}^{X}\circ p.
		$$

		Take $\Phi \in \Omega(1)^{\infty}$ with $f_1=p(\Phi)$ and $f_2=q(\Phi)$. Then
		$$
		f_{1}(1)=\int_{\overline{T}_{E}}\mathrm{pr}_1\left(m(t,\det(t))\cdot \Phi\right)(T_1)\dd t
		$$
		and
		\begin{align*}
			f_2(1)=& \int_{\overline{T}_{E}}\mathrm{pr}_1\left(P(t,t)\cdot \Phi\right)(T_1)  \dd t  \\
			=&\int_{\overline{T}_{E}}\mathrm{pr}_1\left(m(t^{-1},\det(t^{-1}))P(t,t)m(t,\det(t))\cdot \Phi\right)(T_1)  \dd t  \\
			=&\int_{\overline{T}_{E}}m(t^{-1},\det(t^{-1}))P(t,t)\cdot \mathrm{pr}_1\left(m(t,\det(t))\cdot \Phi\right)(T_1)  \dd t \\
			=&\int_{\overline{T}_{E}}\mathrm{pr}_1\left(m(t,\det(t))\cdot \Phi\right)(T_1)  \dd t. 
		\end{align*}
		Hence $f_1(1)=f_2(1)$. On the other hand, by Lemma \ref{zf12: pointwiselemma2},
		$$
		f_1(1)=\int_{\widehat{\PGSO}(V)^{\temp}}\vartheta^{X}_{\sigma}\circ p(\Phi)\dd \mu_{\PGSO(V)}(\sigma),
		$$
		and
		$$
		f_2(1)=\int_{\widehat{\PGSO}(V)^{\temp}}\vartheta^{Y}_{\sigma}\circ q(\Phi)\dd \mu_{\PGSO(V)}(\sigma).
		$$
		
		Thus one has
		\begin{align} 
			0=&f_1(1)-f_2(1)  \label{relativecharacteridentity: eq11} \\
			=&\int_{\widehat{\PGSO}(V)^{\temp}}\vartheta^{X}_{\sigma}\circ p(\Phi)-\vartheta^{Y}_{\sigma}\circ q(\Phi)\dd \mu_{\PGSO(V)}(\sigma)  \notag \\
			=&\int_{\widehat{\PGSO}(V)^{\temp}}\ell_{\sigma}^{X}\circ l_{\sigma}\circ \theta_{\sigma}(\Phi)- \ell_{\sigma}^{X}\circ \ell^{Y}_{\sigma}\circ \theta_{\sigma}(\Phi)\dd \mu_{\PGSO(V)}(\sigma).   \notag
		\end{align}
		
		Now we define a linear functional on $ C_c^{\infty}(\PGSO(V))$. For any $\varphi \in C_c^{\infty}(\PGSO(V))$, denote
		\begin{align}
			&L^{\Phi}_{\sigma}(\varphi)\notag \\
			:=&\ell_{\sigma}^{X}\circ l_{\sigma}\circ \theta_{\sigma}\left(\Omega(1)|_{\PGSO(V)}(\varphi)\Phi\right)- \ell_{\sigma}^{X}\circ \ell^{Y}_{\sigma}\circ \theta_{\sigma}(\Omega(1)|_{\PGSO(V)}(\varphi)\Phi) \notag \\
			=&  \ell_{\sigma}^{X}\left( \sigma(\varphi)l_{\sigma}\circ \theta_{\sigma}(\Phi)\right)- \ell_{\sigma}^{X}\left( \sigma(\varphi)\ell^{Y}_{\sigma}\circ \theta_{\sigma}(\Phi)\right)\notag
		\end{align}
		For every $\Phi \in \Omega(1)^{\infty}$, the linear functional $L^{\Phi}_{\sigma}(\cdot)$ is continuous on $C_c^{\infty}(\PGSO(V))$ and it factors through the map
		$$
		\varphi \mapsto \sigma(\varphi).
		$$

		From the equation (\ref{relativecharacteridentity: eq11}), we have
		\begin{equation} \label{zf5: eq101}
			\int_{\widehat{\PGSO}(V)^{\temp}} L^{\Phi}_{\sigma}(\varphi) \dd \mu_{\PGSO(V)}(\sigma)=0
		\end{equation}
		for all $\varphi \in C_c^{\infty}\left(\PGSO(V)\right)$.

		Suppose $F$ is $p$-adic. Let $\boldsymbol{z}$ be any element in the Bernstein center of $\PGSO(V)$. Then there exists a function $\boldsymbol{z}$ on $\widehat{\PGSO}(V)^{\temp}$ such that
		$$
		\sigma(\boldsymbol{z}\varphi)=\boldsymbol{z}(\sigma)\sigma(\varphi)
		$$
		for any $\varphi \in C_c^{\infty}(\PGSO(V))$. These functions $\{\boldsymbol{z}(\cdot)\in C^{\infty}\left( \widehat{\PGSO}(V)^{\temp}\right)\mid \boldsymbol{z}  \}$ are dense in the space of all Schwartz functions (i.e. smooth functions of rapid decay) on $\widehat{\PGSO}(V)^{\temp}$ which is the space $C_c^{\infty}\left(\widehat{\PGSO}(V)^{\temp} \right)$ in $p$-adic case. By equation (\ref{zf5: eq101}), one has
		$$
		\int_{\widehat{\PGSO}(V)^{\temp}} L^{\Phi}_{\sigma}(\boldsymbol{z}\varphi) \dd \mu_{\PGSO(V)}(\sigma)=\int_{\widehat{\PGSO}(V)^{\temp}}\boldsymbol{z}(\sigma)\cdot L^{\Phi}_{\sigma}(\varphi) \dd \mu_{\PGSO(V)}(\sigma)=0.
		$$
		Hence by Stone-Weierstrass theorem, one can deduce for any $\varphi \in C_c^{\infty}\left(\PGSO(V)\right)$,
		$$
		L_{\sigma}^{\Phi}(\varphi)=0
		$$
		for $\mu_{\PGSO(V)}$-almost all $\sigma \in \widehat{\PGSO}(V)^{\temp}$. There exists a countable dense subset $\{\varphi_i\in C_c^{\infty}\left(\PGSO(V)\right)\mid i \in \mathbb{Z}_{+}\} $ in $C_c^{\infty}\left(\PGSO(V)\right)$. For every $\varphi_i$, there is a measure $0$ subset $\Sigma_i \subset \widehat{\PGSO}(V)^{\temp}$ such that $L_{\sigma}^{\Phi}(\varphi_i)=0$ holds for $\sigma \in \widehat{\PGSO}(V)^{\temp}\setminus \Sigma_i$. Then for every $\varphi \in C_c^{\infty}\left(\PGSO(V) \right)$, one knows
		$$
		L_{\sigma}^{\Phi}(\varphi)=0
		$$
		holds for $\sigma \in \widehat{\PGSO}(V)^{\temp}\setminus \cup_i\Sigma_i$. Hence
		$$
		L_{\sigma}^{\Phi}=0
		$$
		holds for $ \mu_{\PGSO(V)}$-almost all $\sigma \in \widehat{\PGSO}(V)^{\temp}$. If $F$ is Archimedean, then one may use a certain subalgebra of Arthur's algebra of multipliers \cite[Lemma 5.7.2]{BP2}. One can also obtain
		$$
		L_{\sigma}^{\Phi}=0
		$$
		for $ \mu_{\PGSO(V)}$-almost all $\sigma \in \widehat{\PGSO}(V)^{\temp}$ for Archimedean case. Thus we can deduce that
		$$
		\vartheta_{\sigma}^{Y}\circ q= \vartheta_{\sigma}^{X}\circ p
		$$
		holds for $\mu_{\PGSO(V)}$-almost all $\sigma \in \widehat{\PGSO}(V)^{\temp}$.
		

	\end{proof}

	From Theorem $\ref{relative characteridentity: thm1}$, we obtain
	\begin{cor} \label{relative characteridentity: cor1}
		For $\mu_{\PGSO(V)}$-almost all $\sigma \in \widehat{\PGSO}(V)^{\temp}$, the following diagram  
		\[ \begin{tikzcd}
			\Omega(1) \arrow{r}{\theta_{\sigma}} \arrow[swap]{d}{p} & \sigma\boxtimes \theta(\sigma) \arrow{d}{\ell^{Y}_{\sigma}} \\%
			\mathcal{S}^{\heartsuit}(X) \arrow{r}{\alpha^{X}_{\sigma}}& \sigma
		\end{tikzcd}
		\]
		is commutative.
	\end{cor}

\section{\bf Calculations on the unramified setting}

\subsection{\bf Unramified setting} \label{subsection6.1a}
Assume in this section
\begin{itemize}
	\item $F$ is a non-archimedean local field with residual characteristic not 2, and $E$ is an unramified quadratic extension of $F$.
	\item $G$ is defined over $\mathcal{O}_{F}$. The measure $\dd g$ of a group $G(F)$ in this chapter is normalized such that $\int_{G(\mathcal{O}_{F})}\dd g=1$. 
\end{itemize}

	\begin{itemize}
		\item Let $\widetilde{K}\subset \GO(V)$ be a hyperspecial maximal compact subgroup of $\GO(V)$, $\widetilde{K}_1:=\widetilde{K}\cap \OO(V)$, $K:=\widetilde{K}\cap \GSO(V)$ and $K_1:=K\cap \SO(V)$.
		
		\item Let $ K^{\prime}\subset \GSp(W)$ be a hyperspecial maximal compact subgroup of $\GSp(W)$ and $K^{\prime}_1:=K^{\prime}\cap \Sp(W)$. 
		
		\item Let $\widetilde{\sigma}$ be an irreducible unitary $\widetilde{K}$-unramified representation of $\PGO(V)$ with $\widetilde{K}$-fixed vector $\widetilde{v}_0$ such that $\langle \widetilde{v}_0, \widetilde{v}_0  \rangle_{\widetilde{\sigma}}=1$.  If $\widetilde{\sigma}|_{\GSO(V)}$ is irreducible, we denote $\widetilde{\sigma}|_{\GSO(V)}$ by $\sigma$ with $\langle \cdot ,\cdot  \rangle_{\sigma}=\langle \cdot,\cdot \rangle_{\widetilde{\sigma}}$. Then $\sigma$ is an irreducible $K$-unramified representation with $K$-fixed vector $v_0:=\widetilde{v}_0$ such that $\langle v_0,v_0 \rangle_{\sigma}=1$. If $\widetilde{\sigma}|_{\GSO(V)}$ is not irreducible, we define $\sigma$ such that $\widetilde{\sigma}|_{\GSO(V)}=\sigma \oplus \boldsymbol{t}\cdot \sigma $, where $\boldsymbol{t} \in \OO(V)\backslash \SO(V)$ and $\boldsymbol{t} \in \widetilde{K}$, and $\sigma$ is a $K$-unramified representation. Equipped with the inner product $\langle \cdot,\cdot \rangle_{\sigma}=2\cdot \langle \cdot ,\cdot \rangle_{\widetilde{\sigma}}|_{\sigma}$, $\sigma$ is an irreducible unitary $K$-unramified representation with $K$-fixed vector $v_0$ such that $\widetilde{v}_0=v_0+\boldsymbol{t}\cdot v_0$. In this case, $\langle v_0,v_0 \rangle_{\sigma}=1$.
	\end{itemize}

	
	\begin{lem} \label{zf98: lemma6.3}
		We have
		$$
		\dim \sigma^{K_1}=1; \qquad  \dim \widetilde{\sigma}^{\widetilde{K}_1}=1.
		$$
		In particular, exactly one irreducible summand of the restriction of $\sigma $ to $\SO(V)$ is $K_1$-unramified; exactly one irreducible summand of the restriction of $\widetilde{\sigma} $ to $\OO(V)$ is $\widetilde{K}_1$-unramified.  
	\end{lem}	
	\begin{proof}
		We prove the $\GSO(V)$-case, since the $\GO(V)$-case is similar.

		For any $\sigma \hookrightarrow \Ind^{\GSO(V)}_{B}\chi$ for some unramified character $\chi$ of a Borel subgroup $B=TN$, one has
		$$
		0\neq \sigma^{K}\subset \sigma^{K_1} \subset \left(\Ind^{\GSO(V)}_{B}\chi\right)^{K_1}.
		$$
		As $K \subset B\cdot K_1$, $\GSO(V)=B\cdot K_1$. Then
		$$
		\dim \left(\Ind^{\GSO(V)}_{B}\chi\right)^{K_1}=1.
		$$
	\end{proof}
	\begin{itemize}
		\item Suppose $\sigma|_{\SO(V)}=\oplus_i\tau_i$ with $\tau_1$ being $K_1$-unramified, and $\widetilde{\sigma}|_{\OO(V)}=\oplus_i\widetilde{\tau}_i$ with $\widetilde{\tau}_1$ being $\widetilde{K}_1$-unramified. Define the inner product $\langle \cdot ,\cdot \rangle_{\tau_1}:=\langle \cdot,\cdot \rangle_{\sigma}|_{\tau_1}$ and $\langle \cdot ,\cdot \rangle_{\widetilde{\tau}_1}:=\langle \cdot,\cdot \rangle_{\widetilde{\sigma}}|_{\widetilde{\tau}_1}$. By Lemma \ref{zf98: lemma6.3}, $v_0 \in \tau_1$ which is $K_1$-fixed and $\langle v_0,v_0  \rangle_{\tau_1}=1$; $\widetilde{v}_0 \in \widetilde{\tau}_1$ which is $\widetilde{K}_1$-fixed and $\langle \widetilde{v}_0, \widetilde{v}_0 \rangle_{\widetilde{\tau}_1}=1$.
	\end{itemize} 

		\begin{itemize}
			\item If $\widetilde{\tau}_1|_{\SO(V)}$ is irreducible, then let  $\widetilde{\tau}_1|_{\SO(V)}=\tau_1$. If $\widetilde{\tau}_1|_{\SO(V)}$ is not irreducible, then let $\widetilde{\tau}_1|_{\SO(V)}=\tau_1 \oplus \boldsymbol{t}\cdot \tau_1 $. 
			\item  Suppose $\psi$ has conductor $\mathcal{O}_{F}$ so that the associated measure $\dd x$ of $F$ gives $\mathcal{O}_{F}$ volume 1.
			For any $a\in \mathcal{O}_{F}^{\times 2}\backslash \mathcal{O}_{F}^{\times}$, we have the spherical function $\Phi_{0,a}$ in $\omega_{\psi_a}^{+}$ (see Subsection \ref{subsection: measure: weil}). Let $\Phi_0:=\oplus_{a\in \mathcal{O}_{F}^{\times 2}\backslash \mathcal{O}_{F}^{\times}}\Phi_{0,a} \in \Omega(1)$. It is $K\times K^{\prime}$-fixed vector in $\Omega(1)$ and $\langle \Phi_0,\Phi_0 \rangle_{\Omega(1)}=\frac{1}{2}\cdot \left(\oplus_{a\in \mathcal{O}_{F}^{\times 2}\backslash \mathcal{O}_{F}^{\times}}\langle \Phi_{0,a},\Phi_{0,a}  \rangle_{\omega_{\psi_a}}\right)=1$. 
		\end{itemize}
		\begin{defn}
			Set
			$$
			f_0:=p(\Phi_0)\in \mathcal{S}^{\heartsuit}(X);\qquad \phi_0:=q(\Phi_0)\in S(Y).
			$$
			Then we call $f_0\in \mathcal{S}^{\heartsuit}(X)$ and $\phi_0\in S(Y)$ the basic functions.
		\end{defn} 
		Then $f_0$ is $K$-invariant and $\phi_0=1_{\U(2)(\mathcal{O}_{F})\backslash \PGSp(W)(\mathcal{O}_{F})}$ is $K^{\prime}$-invariant. 
		\begin{lem}
			The basic functions $f_0$ and $\phi_0$ are in correspondence.
		\end{lem}
		\begin{proof}
			Since $f_0=p(\Phi_0)$ and $\phi_0=q(\Phi_0)$, by the definition of correspondence, one knows the basic functions $f_0$ and $\phi_0$ are in correspondence.
		\end{proof}

		

		\subsection{\bf The values of $\left|\ell_{\sigma}^{Y}(w_0)\right|$}\label{subsection6.1b}
		We continue to adopt the notations and conventions in Subsection \ref{subsection6.1a}. Recall from Subsection \ref{subsection2.3b} and Section \ref{zf13: 24}
		$$
		Z_{\sigma}^{\PGSO(V)}(\Phi,\Phi^{\prime},v,v^{\prime})=\int_{\PGSO(V)}\langle \Omega(1)(h)\Phi,\Phi^{\prime}  \rangle_{\Omega(1)}\cdot \overline{\langle \sigma(h)v,v^{\prime} \rangle}_{\sigma}\  \dd^{\PGSO(V)} h,
		$$
		where $\Phi,\Phi^{\prime}\in \Omega(1)$ and $v,v^{\prime}\in \sigma$;
		$$
		Z_{\sigma}^{\SO(V)}(\Phi,\Phi^{\prime},v,v^{\prime}):=\int_{\SO(V)}\langle \omega_{\psi}(h)\Phi,\Phi^{\prime}  \rangle_{\omega_{\psi}}\cdot \overline{\langle \sigma(h)v,v^{\prime} \rangle}_{\sigma}\  \dd^{\SO(V)} h,
		$$
		where $\Phi,\Phi^{\prime}\in \omega_{\psi}$ and $v,v^{\prime}\in \sigma$;
		for $\tau$ an irreducible $K_1$-unramified representation,  
		$$
		Z_{\tau}^{\SO(V)}(\Phi,\Phi^{\prime},v,v^{\prime}):=\int_{\SO(V)}\langle \omega_{\psi}(h)\Phi,\Phi^{\prime}  \rangle_{\omega_{\psi}}\cdot \overline{\langle \tau(h)v,v^{\prime} \rangle}_{\tau}\  \dd^{\SO(V)} h,
		$$
		where $\Phi,\Phi^{\prime}\in \omega_{\psi}$ and $v,v^{\prime}\in \tau$.

		For $\widetilde{\tau}$ an irreducible $\widetilde{K}_1$-unramified representation, define 
		$$
		Z_{\widetilde{\tau}}^{\OO(V)}(\Phi,\Phi^{\prime},v,v^{\prime}):=\int_{\OO(V)}\langle \omega_{\psi}(h)\Phi,\Phi^{\prime}  \rangle_{\omega_{\psi}}\cdot \overline{\langle \widetilde{\tau}(h)v,v^{\prime} \rangle}_{\widetilde{\tau}}\  \dd^{\OO(V)} h,
		$$
		where $\Phi,\Phi^{\prime}\in \omega_{\psi}$ and $v,v^{\prime}\in \widetilde{\tau}$.

		\begin{lem} \label{zf98: lemma6.5}
			We have
			$$
			Z_{\sigma}^{\PGSO(V)}(\Phi_0,\Phi_0,v_0,v_0)=Z_{\sigma}^{\SO(V)}(\Phi_{0,1},\Phi_{0,1},v_0,v_0)=Z_{\tau_1}^{\SO(V)}(\Phi_{0,1},\Phi_{0,1},v_0,v_0).
			$$
		\end{lem}
		\begin{proof}
			\begin{align*}
				&Z_{\sigma}^{\PGSO(V)}(\Phi_0,\Phi_0,v_0,v_0) \\
				=&\int_{\PGSO(V)}\langle \Omega(1)(h)\Phi_0,\Phi_0  \rangle_{\Omega(1)}\cdot \overline{\langle \sigma(h)v_0,v_0 \rangle}_{\sigma}\  \dd^{\PGSO(V)} h  \\
				=&\sum_{a \in F^{\times 2}\backslash F^{\times}}\int_{\{\pm 1\}\backslash \SO(V)}\langle \Omega(1)(h_a\cdot h)\Phi_0,\Phi_0  \rangle_{\Omega(1)}\cdot \overline{\langle \sigma(h_a\cdot h)v_0,v_0 \rangle}_{\sigma}\  \dd^{\PGSO(V)} h; \qquad \text{(1)}  \\
				=&\sum_{a \in \mathcal{O}_F^{\times 2}\backslash \mathcal{O}_F^{\times}}\int_{\{\pm 1\}\backslash \SO(V)}\langle \Omega(1)(h_a\cdot h)\Phi_0,\Phi_0  \rangle_{\Omega(1)}\cdot \overline{\langle \sigma(h_a\cdot h)v_0,v_0 \rangle}_{\sigma}\  \dd^{\PGSO(V)} h; \qquad \text{(2)}  \\
				=&2\cdot \int_{\{\pm 1\}\backslash \SO(V)}\langle \Omega(1)( h)\Phi_0,\Phi_0  \rangle_{\Omega(1)}\cdot \overline{\langle \sigma( h)v_0,v_0 \rangle}_{\sigma}\  \dd^{\PGSO(V)} h; \qquad \text{(3)}  \\
				=&\int_{ \SO(V)}\langle \Omega(1)( h)\Phi_0,\Phi_0  \rangle_{\Omega(1)}\cdot \overline{\langle \sigma( h)v_0,v_0 \rangle}_{\sigma}\  \dd^{\SO(V)} h; \qquad \text{(4)}  \\
				=&\int_{ \SO(V)}\sum_{a \in \mathcal{O}_{F}^{\times 2}\backslash \mathcal{O}_{F}^{\times}}\langle \Omega(1)( h)\Phi_{0,a},\Phi_{0,a}  \rangle_{\Omega(1)}\cdot \overline{\langle \sigma( h)v_0,v_0 \rangle}_{\sigma}\  \dd^{\SO(V)} h  \\
				=&2\cdot \int_{ \SO(V)}\langle \Omega(1)( h)\Phi_{0,1},\Phi_{0,1}  \rangle_{\Omega(1)}\cdot \overline{\langle \sigma( h)v_0,v_0 \rangle}_{\sigma}\  \dd^{\SO(V)} h; \qquad \text{(5)}  \\
				=&\int_{ \SO(V)}\langle \omega_{\psi}( h)\Phi_{0,1},\Phi_{0,1}  \rangle_{\omega_{\psi}}\cdot \overline{\langle \sigma( h)v_0,v_0 \rangle}_{\sigma}\  \dd^{\SO(V)} h  \\
				=&\int_{ \SO(V)}\langle \omega_{\psi}( h)\Phi_{0,1},\Phi_{0,1}  \rangle_{\omega_{\psi}}\cdot \overline{\langle \sigma( h)v_0,v_0 \rangle}_{\tau_1}\  \dd^{\SO(V)} h  \\
				=&Z_{\tau_1}^{\SO(V)}(\Phi_{0,1},\Phi_{0,1},v_0,v_0),
			\end{align*}
			where (1) is because $\PGSO(V)=\sqcup_{a \in F^{\times 2}\backslash F^{\times}}h_a\cdot \{\pm 1\}\backslash \SO(V)$; (2) is because if $a \notin \mathcal{O}_{F}^{\times}$, $h_a\cdot \Phi\in \oplus_{b \notin \mathcal{O}_{F}^{\times 2}\backslash \mathcal{O}_{F}^{\times}}\omega_{\psi_{b}}^{+}$ whereas $\Phi \in \oplus_{b \in \mathcal{O}_{F}^{\times 2}\backslash \mathcal{O}_{F}^{\times}}\omega_{\psi_{b}}^{+}$; (3) is because the contribution of the two terms for $a \in \mathcal{O}_{F}^{\times 2}\backslash \mathcal{O}_{F}^{\times}$ are same; (4) is because one replaced the domain of integration by $\SO(V)$ and the measure by $\dd^{\SO(V)}h$; (5) is because the two terms have the same contribution.
		\end{proof}
		
		\begin{lem}
			We have 
			$$
			Z^{\OO(V)}_{\widetilde{\tau_1}}(\Phi_{0,1},\Phi_{0,1},\widetilde{v}_0,\widetilde{v}_0)=Z^{\SO(V)}_{\tau_1}(\Phi_{0,1},\Phi_{0,1},v_0,v_0)
			$$
		\end{lem}
		\begin{proof}
			If $\widetilde{\tau_1}|_{\SO(V)}=\tau_1 \oplus \boldsymbol{t}\cdot \tau_1$ where $\boldsymbol{t}\in \widetilde{K}$, then
			\begin{align*}
				&Z^{\OO(V)}_{\widetilde{\tau_1}}(\Phi_{0,1},\Phi_{0,1},\widetilde{v}_0,\widetilde{v}_0) \\
				=&\int_{\OO(V)}\langle \omega_{\psi}(h)\Phi_{0,1},\Phi_{0,1}  \rangle_{\omega_{\psi}}\cdot \overline{\langle \widetilde{\tau_1}(h)\widetilde{v}_0, \widetilde{v}_0 \rangle}_{\widetilde{\tau_1}}\  \dd^{\OO(V)} h \\
				=&\int_{\SO(V)}\langle \omega_{\psi}(h)\Phi_{0,1},\Phi_{0,1}  \rangle_{\omega_{\psi}}\cdot \overline{\langle \widetilde{\tau_1}(h)\widetilde{v}_0, \widetilde{v}_0 \rangle}_{\widetilde{\tau_1}}\  \dd^{\OO(V)} h \ + \\
				&\int_{\SO(V)}\langle \omega_{\psi}(h\cdot \boldsymbol{t})\Phi_{0,1},\Phi_{0,1}  \rangle_{\omega_{\psi}}\cdot \overline{\langle \widetilde{\tau_1}(h\cdot \boldsymbol{t})\widetilde{v}_0, \widetilde{v}_0 \rangle}_{\widetilde{\tau_1}}\  \dd^{\OO(V)} h  \\
				=& 2\cdot \int_{\SO(V)}\langle \omega_{\psi}(h)\Phi_{0,1},\Phi_{0,1}  \rangle_{\omega_{\psi}}\cdot \overline{\langle \widetilde{\tau_1}(h)\widetilde{v}_0, \widetilde{v}_0 \rangle}_{\widetilde{\tau_1}}\  \dd^{\OO(V)} h; \qquad \text{($\omega_{\psi}(\boldsymbol{t})\Phi_{0,1}=\Phi_{0,1}$; and $\widetilde{\tau_1}(\boldsymbol{t})\widetilde{v}_0=\widetilde{v}_0$)} \\
				=&\int_{\SO(V)}\langle \omega_{\psi}(h)\Phi_{0,1},\Phi_{0,1}  \rangle_{\omega_{\psi}}\cdot \overline{\langle \widetilde{\tau_1}(h)\widetilde{v}_0, \widetilde{v}_0 \rangle}_{\widetilde{\tau_1}}\  \dd^{\SO(V)} h; \qquad \text{($\dd^{\SO(V)} h=2\cdot \dd^{\OO(V)}h$)}  \\
				=&\int_{\SO(V)}\langle \omega_{\psi}(h)\Phi_{0,1},\Phi_{0,1}  \rangle_{\omega_{\psi}}\cdot \left(\overline{\langle \widetilde{\tau_1}(h) v_0, v_0 \rangle}_{\widetilde{\tau_1}}+\overline{\langle \widetilde{\tau_1}(h\cdot \boldsymbol{t}) v_0, \widetilde{\tau_1}(\boldsymbol{t})v_0 \rangle}_{\widetilde{\tau_1}}\right)\  \dd^{\SO(V)} h \\
				=&2\cdot \int_{\SO(V)}\langle \omega_{\psi}(h)\Phi_{0,1},\Phi_{0,1}  \rangle_{\omega_{\psi}}\cdot \overline{\langle \widetilde{\tau_1}(h)v_0, v_0 \rangle}_{\widetilde{\tau_1}}\  \dd^{\SO(V)} h \\
				=&\int_{\SO(V)}\langle \omega_{\psi}(h)\Phi_{0,1},\Phi_{0,1}  \rangle_{\omega_{\psi}}\cdot \overline{\langle \tau_1(h)v_0, v_0 \rangle}_{\tau_1}\  \dd^{\SO(V)} h \\
				=&Z^{\SO(V)}_{\tau_1}(\Phi_{0,1},\Phi_{0,1},v_0,v_0).
			\end{align*}
			
			If $\widetilde{\tau_1}|_{\SO(V)}$ is irreducible, then
			\begin{align*}
				&Z^{\OO(V)}_{\widetilde{\tau_1}}(\Phi_{0,1},\Phi_{0,1},\widetilde{v}_0,\widetilde{v}_0) \\
				=&\int_{\OO(V)}\langle \omega_{\psi}(h)\Phi_{0,1},\Phi_{0,1}  \rangle_{\omega_{\psi}}\cdot \overline{\langle \widetilde{\tau_1}(h)\widetilde{v}_0, \widetilde{v}_0 \rangle}_{\widetilde{\tau_1}}\  \dd^{\OO(V)} h \\
				=&\int_{\SO(V)}\langle \omega_{\psi}(h)\Phi_{0,1},\Phi_{0,1}  \rangle_{\omega_{\psi}}\cdot \overline{\langle \widetilde{\tau_1}(h)\widetilde{v}_0, \widetilde{v}_0 \rangle}_{\widetilde{\tau_1}}\  \dd^{\OO(V)} h\  +\\
				&\int_{\SO(V)}\langle \omega_{\psi}(h\cdot \boldsymbol{t})\Phi_{0,1},\Phi_{0,1}  \rangle_{\omega_{\psi}}\cdot \overline{\langle \widetilde{\tau_1}(h\cdot \boldsymbol{t})\widetilde{v}_0, \widetilde{v}_0 \rangle}_{\widetilde{\tau_1}}\  \dd^{\OO(V)} h  \\
				=& 2\cdot \int_{\SO(V)}\langle \omega_{\psi}(h)\Phi_{0,1},\Phi_{0,1}  \rangle_{\omega_{\psi}}\cdot \overline{\langle \widetilde{\tau_1}(h)\widetilde{v}_0, \widetilde{v}_0 \rangle}_{\widetilde{\tau_1}}\  \dd^{\OO(V)} h \\
				=&\int_{\SO(V)}\langle \omega_{\psi}(h)\Phi_{0,1},\Phi_{0,1}  \rangle_{\omega_{\psi}}\cdot \overline{\langle \widetilde{\tau_1}(h)\widetilde{v}_0, \widetilde{v}_0 \rangle}_{\widetilde{\tau_1}}\  \dd^{\SO(V)} h \\
				=&\int_{\SO(V)}\langle \omega_{\psi}(h)\Phi_{0,1},\Phi_{0,1}  \rangle_{\omega_{\psi}}\cdot \overline{\langle \tau_1(h) v_0, v_0 \rangle}_{\tau_1}\  \dd^{\SO(V)} h \\
				=&Z^{\SO(V)}_{\tau_1}(\Phi_{0,1},\Phi_{0,1},v_0,v_0).
			\end{align*}
		\end{proof}


		By \cite[Proposition 3]{LR}, one obtains
		\begin{lem} \label{zf98: lemma6.}
			We have
			\begin{align*}
				Z_{\sigma}^{\PGSO(V)}(\Phi_0,\Phi_0,v_0,v_0)=&Z_{\tau_1}^{\SO(V)}(\Phi_{0,1},\Phi_{0,1},v_0,v_0)\\
				=&Z^{\OO(V)}_{\widetilde{\tau}_1}(\Phi_{0,1},\Phi_{0,1},\widetilde{v}_0,\widetilde{v}_0) \\
				=&\frac{ L\left(1,\widetilde{\sigma},\mathrm{std}\right)}{ \zeta_{F}(2)\cdot \zeta_{F}(4)} \\
				=& \frac{ L\left(1,\sigma_1\times \sigma_2\right)}{ \zeta_{F}(2)\cdot \zeta_{F}(4)}
			\end{align*}
			
		\end{lem}
		\begin{proof}
			\begin{align*}
				\frac{ L\left(1,\widetilde{\sigma},\mathrm{std}\right)}{ \zeta_{F}(2)\cdot \zeta_{F}(4)}=&\int_{\OO(V)}\langle \omega_{\psi}(h)\Phi_{0,1},\Phi_{0,1} \rangle_{\omega_{\psi}}\cdot \overline{\langle \widetilde{\sigma}(h)\widetilde{v}_0,\widetilde{v}_0 \rangle}_{\widetilde{\sigma}}\dd h \\
				=&\int_{\OO(V)}\langle \omega_{\psi}(h)\Phi_{0,1},\Phi_{0,1} \rangle_{\omega_{\psi}}\cdot \overline{\langle \widetilde{\tau}_1(h)\widetilde{v}_0,\widetilde{v}_0 \rangle}_{\widetilde{\tau}_1}\dd h \\
				=&Z^{\OO(V)}_{\widetilde{\tau}_1}(\Phi_{0,1},\Phi_{0,1},\widetilde{v}_0,\widetilde{v}_0) .
			\end{align*}

		\end{proof}
		By Lemma \ref{zf98: lemma6.}, one knows $0\neq A_{\sigma}(\Phi_0,v_0)\in \theta(\sigma)$ is $K^{\prime}$-fixed ( for definition of $A_{\sigma}$, see Subsection \ref{subsection22a} ). Then $\theta(\sigma)$ is a $K^{\prime}$-unramified representation of $\PGSp(W)$ with some $K^{\prime}$-fixed vector $w_0$ such that $\langle w_0,w_0 \rangle_{\theta(\sigma)}=1$.

		Recall we have fixed the map in Subsection \ref{subsection22a}
		$$
		\theta_{\sigma}:\Omega(1)\longrightarrow \sigma\boxtimes\theta(\sigma).
		$$
		Then we have
		$$
		\theta_{\sigma}\left(\Phi_0\right)=c_{\sigma}\cdot v_0\otimes w_0 \in \sigma\boxtimes\theta(\sigma)
		$$
		for some nonzero constant $c_{\sigma}$. Actually, we have
		\begin{align*}
			Z^{\PGSO(V)}_{\sigma}\left(\Phi_0,\Phi_0,v_0,v_0 \right)=&\langle A_{\sigma}\left( \Phi_0,v_0\right),A_{\sigma}\left( \Phi_0,v_0\right) \rangle_{\theta(\sigma)} \\
			=&\langle \langle \theta_{\sigma}\left(\Phi_0\right),v_0 \rangle_{\sigma},\langle\theta_{\sigma}\left(\Phi_0\right),v_0 \rangle_{\sigma} \rangle_{\theta(\sigma)} \\
			=&\langle c_{\sigma}\cdot w_0, c_{\sigma}\cdot w_0  \rangle_{\theta(\sigma)} \\
			=&\left|c_{\sigma}\right|^2.
		\end{align*}
		Similarly, we have
		$$
		\langle B_{\sigma}\left( \Phi_0,w_0\right),B_{\sigma}\left( \Phi_0,w_0\right) \rangle_{\sigma} =\left|c_{\sigma}\right|^2.
		$$
		We can use the similar computation of Lemma \ref{zf98: lemma6.5} to obtain the following lemma.
		\begin{lem} \label{lem6.8}
			We have
			$$
			Z_{\sigma}^{\PGSO(V)}(\Phi_{0},\Phi_{0},v_0,v_0)=4\cdot Z_{\sigma}^{\PGSO(V)}(\Phi_{0,1},\Phi_{0,1},v_0,v_0).
			$$
			Thus
			$$
			A_{\sigma}(\Phi_0,v_0)=\gamma_1\cdot A_{\sigma}(\Phi_{0,1},v_0)
			$$
			and
			$$
			B_{\sigma}(\Phi_0,w_0)=\gamma_2\cdot B_{\sigma}(\Phi_{0,1},w_0)
			$$
			with $\left|\gamma_1\right|=2$ and $\left|\gamma_2\right|=2$.
		\end{lem}
		
		
		Then we have
		$$
		\alpha_{\sigma}^{Y}(\phi_0)=a_{\sigma}^{Y}\cdot w_0
		$$
		for some constant $a_{\sigma}^{Y}$.


		Let $\kappa_{F}$ be the residue field of $F$ and $q=\# \kappa_{F}$. 
		Recall we have fixed the measure $\lambda\cdot \frac{\dd^{\PGSp(W)} g}{\dd s\cdot \dd t}$ on $\U(2)(F)\backslash \PGSp(W)(F)$ such that the following equation hold:
		\begin{align*}
			&\int_{\U(2)(\mathcal{O}_{F})\backslash \PGSp(W)(\mathcal{O}_{F})}\int_{\overline{T}_{E}(\mathcal{O}_{F})} \ \dd t \lambda\cdot \frac{\dd^{\PGSp(W)} g}{\dd s\cdot \dd t}\\
			=&\int_{\SU(2)(\mathcal{O}_{F})\backslash \PGSp(W)(\mathcal{O}_{F})} \ \lambda\cdot \frac{\dd^{\PGSp(W)}g }{\dd s}\\
			=&\int_{^{1}Y(\mathcal{O}_{F})}\left|\omega_1\right| \\
			=&2 \int_{Y_1(\mathcal{O}_{F})}\  \frac{1}{2}\cdot \frac{\left|\omega_{\Sp(W)}\right|}{\left|\omega_{\SU(2)}\right|} \qquad \text{(by equation (\ref{innerproductomega1}))}\\
			=&  (1-q^{-4}),
		\end{align*}
		where $\omega_{\Sp(W)}$ is a $\Sp(W)$-invariant differential of top degree which has nonzero reduction on the special fiber; and $\omega_{\SU(2)}$ is a $\SU(2)$-invariant differential of top degree which has nonzero reduction on the special fiber. Thus $\lambda=1-q^{-4}$. For definitions of $^{1}Y$ and $Y_1$, see Subsection \ref{subsection4.1.a}.
		
		One has
		\begin{align*}
			a_{\sigma}^{Y}=&\langle \alpha_{\sigma}^{Y}(\phi_0), w_0 \rangle_{\theta(\sigma)}  \\
			=&\langle \phi_0, \beta_{\sigma}^{Y}(w_0) \rangle_{Y} \\
			=&\int_{\U(2)(F)\backslash \PGSp(W)(F)}1_{\U(2)(\mathcal{O}_{F})\backslash \SO(5)(\mathcal{O}_{F})}(g)\cdot \overline{\beta_{\sigma}^{Y}(w_0)(g)}\ \lambda\cdot \frac{\dd^{\PGSp(W)} g}{\dd s\cdot \dd t}  \\
			=&\int_{\U(2)(\mathcal{O}_{F})\backslash \PGSp(W)(\mathcal{O}_{F})} \overline{\ell_{\sigma}^{Y}(w_0)}\  \lambda\cdot \frac{\dd^{\PGSp(W)} g}{\dd s\cdot \dd t} \qquad \text{( $w_0$ is $\PGSp(W)(\mathcal{O}_{F})$-invariant )} \\
			=&\overline{\ell_{\sigma}^{Y}(w_0)}\cdot\lambda.
		\end{align*}

		One also has
		\begin{align*}
			c_{\sigma}\cdot \ell_{\sigma}^{X}(v_0)\cdot w_0=\ell_{\sigma}^{X}\left(\theta_{\sigma}(\Phi_0)\right)=a_{\sigma}^{Y}\cdot w_0,
		\end{align*}
		where the second equality is by the Corollary \ref{commutativediagram: 2}.

		Thus we can conclude
		$$
		\left|   \ell_{\sigma}^{Y}(w_0) \right|^2=Z^{\PGSO(V)}_{\sigma}\left(\Phi_0,\Phi_0,v_0,v_0 \right)\cdot \left|\ell^{X}_{\sigma}(v_0)\right|^2\cdot \lambda^{-2}.
		$$
		
		\begin{lem}
			For every irreducible unramified representation $\sigma\cong\sigma_1\boxtimes\sigma_2 \in \widehat{\PGSO}(V)^{\temp}\cong \widehat{\PGL}_2^{\temp}\boxtimes \widehat{\PGL}_2^{\temp}$, one has
			\begin{align*}
				\left|\ell^{X}_{\sigma}(v_0)\right|^2=&\left|\ell^{N,\psi\backslash \PGL_2}_{\sigma_1}(^{1}v_0)\right|^2\cdot \left|\ell^{\overline{T}_{E}\backslash \PGL_2}_{\sigma_2}(^{2}v_0)\right|^2  \\
				=&\frac{\zeta_{F}(2)^2\cdot L(\frac{1}{2},\mathrm{BC}(\sigma_2),\mathrm{std})}{L(1,\sigma_1,\mathrm{Ad})\cdot L(1,\sigma_2,\mathrm{Ad})\cdot L(1,\chi_{E\slash F},\mathrm{std})}.
			\end{align*}
			where ${^{1}v_0} \in \sigma_1$ and ${^{2}v_0} \in \sigma_2$ are such that $\langle {^{i}v_0}, {^{i}v_0}\rangle_{\sigma_i}=1$ and $v_0={^{1}v_0}\otimes {^{2}v_0}$,  $\mathrm{BC}(\sigma_2)$ denotes the base change of $\sigma_2$ to $\GL_2(E)$, and $\chi_{E/F}$ is a quadratic character of $F^{\times}$ determined by $E$ via the local class field theory. 
		\end{lem}
		\begin{proof}
			By \cite[Proposition 2.14]{LM}, one has 
			$$
			\left|\ell^{N,\psi\backslash \PGL_2}_{\sigma_1}(^{1}v_0)\right|^2=\frac{\zeta_{F}(2)}{L(1,\sigma_1,\mathrm{Ad})}.
			$$
			By \cite[Proposition 7]{W2} or \cite[Theorem 1.2]{II}, one has
			$$
			\left|\ell^{\overline{T}_{E}\backslash \PGL_2}_{\sigma_2}(^{2}v_0)\right|^2=\frac{\zeta_{F}(2)\cdot L(\frac{1}{2},\mathrm{BC}(\sigma_2),\mathrm{std})}{L(1,\sigma_2,\mathrm{Ad})\cdot L(1,\chi_{E\slash F},\mathrm{std})}.
			$$
			In conclusion, we have
			$$
			\left|\ell^{X}_{\sigma}(v_0)\right|^2=\frac{\zeta_{F}(2)^2\cdot L(\frac{1}{2},\mathrm{BC}(\sigma_2),\mathrm{std})}{L(1,\sigma_1,\mathrm{Ad})\cdot L(1,\sigma_2,\mathrm{Ad})\cdot L(1,\chi_{E\slash F},\mathrm{std})}.
			$$
		\end{proof}
		

		Combining with above discussions, one obtains
		\begin{prop}
			$$
			\left|\ell_{\sigma}^{Y}(w_0)\right|^2=\zeta_{F}(2)\zeta_{F}(4)\cdot \frac{L(1,\sigma_1\times \sigma_2)\cdot L(\frac{1}{2},\mathrm{BC}(\sigma_2),\mathrm{std})}{L(1,\sigma_1,\mathrm{Ad})\cdot L(1,\sigma_2,\mathrm{Ad})\cdot L(1,\chi_{E/F},\mathrm{std}) }.
			$$
		\end{prop}

\section{Global Results} \label{sec3}
\subsection{\bf Tamagawa measures}

Let $k$ be a number field with the ring of ad\`eles $\mathbb{A}$. Fix a non-trivial unitary character 
$$
\psi: k\backslash \mathbb{A}\longrightarrow S^{1}.
$$
This character has a factorization 
$$
\psi=\prod_v \psi_v,
$$ 
where $\psi_v$ is a non-trivial unitary character of $k_v$ for each place $v$ of $k$. 

For every place $v$, we fix a self-dual measure $d_{\psi_v}x$ on $k_v$ with respect to the pairing $(x,y):=\psi_{v}(xy)$. Then the volume of the ring of integers $\mathcal{O}_{k_v}$ relative to $d_{\psi_v} x$ is $1$ for almost all place $v$. The product measure $d x:=\prod_v d_{ \psi_v}x $ of $\mathbb{A}$ is independent of $\psi$ and satisfies 
$$
\int_{k\backslash \mathbb{A}}d x=1,
$$
which is called the Tamagawa measure of $\mathbb{A}$.

If $G$ is a smooth algebraic group over $k$, we may consider the adelic group $G(\mathbb{A})=\prod^{\prime}_{v}G(k_v)$ with respect to open compact subgroups $\{G(\mathcal{O}_{k_v}) \}$.

Suppose $\omega_{G}$ is a nonzero invariant differential of top degree on $G$ over $k$. Then for each place $v$ of $k$, the pair $(\omega_{G},\psi_v)$ determines a Haar measure $|\omega_{G,\psi_v}|_v$ of $G(k_v)$. We would like to consider the product measure $\prod_v|\omega_{G,\psi_v}|_v$ on $G(\mathbb{A})$. When $G$ is semi-simple or unipotent, this product measure is well-defined and independent of $\psi$ and $\omega_{G}$, which is called the Tamagawa measure of $G$. If the infinite product is not convergent (e.g.
if $G=G_{m}$), one can still deal with this by introducing "normalization factors".

We fix some measures in the following sections.
\begin{enumerate}
	\item We choose a Haar measure $d^{PGSO(V)}h_v$ on $PGSO(V)(k_v)$ with
	$$
	\int_{PGSO(V)(\mathcal{O}_{k_v})}d^{PGSO(V)}h_v=1
	$$
	for every $p$-adic place $v$ such that
	$$
	d^{PGSO(V)}h=\prod_{v}d^{PGSO(V)}h_v,
	$$ 
	where $d^{PGSO(V)}h$ is the Tamagawa measure of $PGSO(V)(\mathbb{A})$. Similarly, we choose $d^{SO(V)}h_v$ on $SO(V)(k_v)$. By definition, for every function $f \in C_c^{\infty}(\{\pm 1 \}\backslash SO(V)(k_v))$ where $v$ is $p$-adic, we have 
	\begin{equation} \label{compatmeasure1}
		\int_{\{\pm 1 \}\backslash SO(V)(k_v)}f\  \frac{d^{SO(V)}h_v}{\nu}= \int_{\{\pm 1 \}\backslash SO(V)(k_v)}f \ d^{PGSO(V)}h_v,
	\end{equation}
	where $\nu$ is the counting measure of $\{\pm 1 \}$.
	
	\item  We choose a Haar measure $d^{PGSp(W)}g_v$ on $PGSp(W)(k_v)$ with
	$$
	\int_{PGSp(W)(\mathcal{O}_{k_v})}d^{PGSp(W)}g_v=1
	$$
	for every $p$-adic place $v$ such that
	$$
	d^{PGSp(W)}g=\prod_{v}d^{PGSp(W)}g_v,
	$$ 
	where $d^{PGSp(W)}g$ is the Tamagawa measure of $PGSp(W)(\mathbb{A})$. Similarly, we choose $d^{Sp(W)}g_v$ on $Sp(W)(k_v)$. By definition, for every function $f \in C_c^{\infty}(\{\pm 1 \}\backslash Sp(W)(k_v))$ where $v$ is $p$-adic, we have 
	\begin{equation} \label{compatmeasure2}
		\int_{\{\pm 1 \}\backslash Sp(W)(k_v)}f\  \frac{d^{Sp(W)}g_v}{\nu}= \int_{\{\pm 1 \}\backslash Sp(W)(k_v)}f \ d^{PGSp(W)}g_v,
	\end{equation}
	where $\nu$ is the counting measure of $\{\pm 1 \}$.
	
	\item We use the same way to fix a Haar measure $d n_v=d_{\psi_v}x$ on $N(k_v)\cong k_v$. The Tamagawa measure $d n$ on $N(\mathbb{A})$ is given by $d n=\prod_v d n_v$ ( $N$ is defined in Subsection \ref{zf13:21a} ).
	
	\item  We use the same way to fix a Haar measure $d s_v$ on $SU(2)(k_v)=SL_2\times 1(k_v)$. The Tamagawa measure $d s$ on $SU(2)(\mathbb{A})=SL_2\times 1(\mathbb{A})$ is given by $d s=\prod_v d s_v$.
	
	\item Let $k^{\prime}$ be a quadratic field of $k$ to which we associate an anisotropic torus $\overline{T}_{k^{\prime}}$ of $PGL_2$. Then there exists an exact sequence
	$$
	1\longrightarrow k^{\times}\longrightarrow k^{\prime\times}\longrightarrow \overline{T}_{k^{\prime}}\longrightarrow 1
	$$ 
	and
	$$
	1\longrightarrow k_v^{\times}\longrightarrow k^{\prime\times}_{v}\longrightarrow \overline{T}_{k^{\prime}_v}\longrightarrow 1.
	$$
	
	Let $d t$ be the Tamagawa measure on $\overline{T}_{k^{\prime}}(\mathbb{A})$. We equip $k_v^{\times}$ with the Haar measure $\zeta_{k_v}(1)^{-1}\cdot \frac{d_{\psi_v}x}{\left|x\right|_{k_v}}$; $k_v^{\prime\times}$ with the Haar measure $\zeta_{K_v}(1)^{-1}\cdot \frac{d_{\psi_v\circ Tr}x}{\left|x\right|_{K_v}}$. Then we equip $\overline{T}_{k^{\prime}_v}$ with the quotient measure $d t_v$ for almost all $v$ such that
	$$
	d t=\prod_v d t_v.
	$$
	\item We equip $U(2)(\mathbb{A})$ with the Tamagawa measure $d u$. It is determined by the Tamagawa measures $d t$ on $\overline{T}_{k^{\prime}}(\mathbb{A})$ and $d s$ on $SU(2)(\mathbb{A})$ via the following exact sequence
	$$
	1\longrightarrow SU(2)(\mathbb{A})\longrightarrow U(2)(\mathbb{A})\longrightarrow \overline{T}_{k^{\prime}}(\mathbb{A})\longrightarrow 1.
	$$ 
	We equip $U(2)(k_v)$ with the measure $d u_v$ which is determined by the measure $d t_v$ on $\overline{T}_{k^{\prime}}(k_v)$ and $d s_v$ on $SU(2)(k_v)$ via the following exact sequence
	$$
	1\longrightarrow SU(2)(k_v)\longrightarrow U(2)(k_v)\longrightarrow \overline{T}_{k^{\prime}}(k_v)\longrightarrow 1.
	$$
	We have $d u=\prod_v d u_v$.
\end{enumerate}

\subsection{\bf Automorphic forms} 
For a semisimple group $G$ defined over $k$, we write $\mathcal{A}(G)$ for the space of smooth automorphic forms on $G$, and $\mathcal{A}_0(G)$ for the subspace of cusp forms. On $\mathcal{A}_0(G)$, we have the Petersson inner product $\langle \cdot,\cdot \rangle_G$ (defined using the Tamagawa measure $d g$ on $G$ divided by the counting measure on the discrete subgroup $G(k)$):
$$
\langle \phi_1, \phi_2  \rangle_{G}:=\int_{[ G]}\phi_1(g)\overline{\phi_2(g)} d g,
$$
where $[G]:=G(k)\backslash G(\mathbb{A})$. The inner product actually defines a pairing between $\mathcal{A}_0(G)$ and $C^{\infty}_{\mathrm{mod}}([G])$, where $C^{\infty}_{\mathrm{mod}}([G])$ denotes the space of smooth functions on $[G]$ with uniform moderate growth. 

\subsection{\bf Global Weil representation}
In this subsection, we will introduce the global Weil representation $\omega_{\psi}$ of  $SO\left(V\right)(\mathbb{A})\times Sp\left(W\right)(\mathbb{A})$.

Let the global Weil representation $\omega_{\psi}$ be the restricted tensor product 
$$
S\left(V_2\otimes W(k\otimes_{\mathbb{Q}}\mathbb{R})\right)\otimes\left(\otimes_{v<\infty}^{\prime}S\left(V_2\otimes W(k_v)\right)\right)=S(M_{4\times 2}(k\otimes_{\mathbb{Q}}\mathbb{R}))\otimes \left( \otimes_{v<\infty}^{\prime}S(M_{4\times 2}(k_v))\right)
$$
with respect to the characteristic function of $M_{4\times 2}(\mathcal{O}_{k_v})$.

Let the automorphic realization 
$$
\theta: \omega_{\psi}\longrightarrow C_{\mathrm{mod}}^{\infty}\left([SO(V)\times Sp(W)]\right)
$$
be defined by
$$
\theta(\Phi)(h,g):=\sum_{T \in V_2\otimes W(k)}\left(\omega_{\psi}(h,g)\Phi\right)(T)
$$
Let 
$$
R_0(\mathbb{A}):=\{(h,g)\in GSO(V)(\mathbb{A})\times GSp(W)(\mathbb{A}) \mid \lambda_{V}(h)\lambda_{W}(g)=1 \}.
$$
One can extend the action of $SO(V)\times Sp(W)$ on this global Weil representation to $R_0(\mathbb{A})$ in the same way as equation (\ref{equation321}). The automorphic realization $\theta$ of $\omega_{\psi}$ extends to one for $R_0$:
$$
\theta: \omega_{\psi}\longrightarrow C_{\mathrm{mod}}^{\infty}\left([R_0]\right).
$$

\subsection{\bf Global theta liftings} \label{subsection: globaltheta}

In this thesis, we will consider the global theta lifting from $\mathcal{A}_0(PGSO(V))$ to $\mathcal{A}_{0}( PGSp(W))$.

Let $\Sigma\subset \mathcal{A}_{0}\left(PGSO(V) \right)$ be an irreducible unitary cuspidal automorphic representation of $PGSO(V)(\mathbb{A})$. We define an $R_0(\A)$-equivariant map
$$
A_{\Sigma}^{\mathrm{Aut}}:\omega_{\psi}\otimes \overline{\Sigma}\longrightarrow \mathcal{A}\left(PGSp(W)\right)
$$
by
$$
A_{\Sigma}^{\mathrm{Aut}}(\Phi,f)(g)=\int_{[SO(V)]}\theta(\Phi)(h_1\cdot h,g)\overline{f(h_1\cdot h)} d h_1
$$
for any $g \in GSp(W)(\mathbb{A})$ with some $h\in GSO(V)(\mathbb{A})$ such that $\lambda_{V}(h)\cdot\lambda_W(g)=1$, where $d h_1$ is the Tamagawa measure of $SO(V)(\mathbb{A})$, $f \in \Sigma$ and $\Phi \in \omega_{\psi}$.

The image of $A_{\Sigma}^{\mathrm{Aut}}$ is the global theta lift of $\Sigma$, which we will denote by
$$
\Pi=\Theta^{\mathrm{Aut}}(\Sigma)\subset\mathcal{A}\left(PGSp(W)\right).
$$
If $\Sigma\cong \Sigma_1\boxtimes\Sigma_2$ is an irreducible cuspidal representation of $PGSO(V)$, where $\Sigma_1\ncong \Sigma_2$, then one knows
$\Pi$ is cuspidal and non-zero \cite[\S 7.2]{GI}. Hence $\Pi$ is an irreducible cuspidal automorphic representation of $PGSp(W)$.

Conversely, let $\Pi \subset \mathcal{A}_{0}\left(PGSp(W) \right)$ be an irreducible cuspidal automorphic representation of $PGSp(W)(\mathbb{A})$. We define an $R_0(\A)$-equivariant map
$$
B_{\Pi}^{\mathrm{Aut}}:\omega_{\psi}\otimes \overline{\Pi}\longrightarrow \mathcal{A}\left(PGO(V)\right)
$$
by
$$
B_{\Pi}^{\mathrm{Aut}}(\Phi,\varphi)(h)=\int_{[Sp(W)]}\theta(\Phi)(h,g_1\cdot g)\cdot \overline{\varphi(g_1\cdot g)}\ d g_1
$$
for any $h \in GO(V)(\mathbb{A})$ with some $g\in GSp(W)(\mathbb{A})$ such that $\lambda_{W}(g)\cdot\lambda_V(h)=1$, where $d g_1$ is the Tamagawa measure of $Sp(W)(\mathbb{A})$, $\varphi \in \Pi$ and $\Phi \in \omega_{\psi}$.
The image of $B_{\Pi}^{\mathrm{Aut}}$ is the global theta lift of $\Pi$, which we will denote by $\Theta^{\mathrm{Aut}}(\Pi)$. By computing the constant term, one has 
$$
\widetilde{\Sigma}:=\Theta^{\mathrm{Aut}}(\Pi)\subset\mathcal{A}_0\left(GO(V)\right).
$$
Then 
$$
\mathrm{Rest}|_{GSO(V)}\Theta^{\mathrm{Aut}}(\Pi)\subset\mathcal{A}_0\left(PGSO(V)\right), 
$$
where $\mathrm{Rest}|_{GSO(V)}\Theta^{\mathrm{Aut}}(\Pi)$ is the restriction of $\Theta^{\mathrm{Aut}}(\Pi)$ to $GSO(V)$ as functions. 
For any irreducible cuspidal automorphic representation $\widetilde{\Sigma}$ of $PGO(V)(\A)$ with nonzero cuspidal theta lift to $PGSp(W)(\A)$, one has
$$
\mathrm{Rest}|_{GSO(V)}\widetilde{\Sigma}= \Sigma_1\boxtimes\Sigma_2\bigoplus\Sigma_2\boxtimes\Sigma_1,
$$ 
where $\Sigma_1\boxtimes\Sigma_2$ is an irreducible cuspidal automorphic representation of $PGSO(V)(\mathbb{A})$ with $\Sigma_1 \ncong \Sigma_2$ (\cite[\S 1]{HST}). 

For each irreducible cuspidal automorphic representation $\Sigma \subset \mathrm{Rest}|_{GSO(V)}\widetilde{\Sigma}$, we set
$$
B^{Aut}_{\Sigma}:=\mathrm{pr}_{\Sigma}\circ \mathrm{Rest}|_{GSO(V)}B^{Aut}_{\Pi}.
$$
If $\Pi=\Theta^{Aut}(\Sigma)$ is a nonzero irreducible 
cuspidal automorphic representation of $PGSp(W)(\A)$, then
$\Sigma \subset \mathrm{Rest}_{GSO(V)}\Theta^{Aut}(\Pi)$, since the cuspidal spectrum of $PGSO(V)\cong PGL_2\times PGL_2$ is multiplicity free.

\begin{lem}
	For $\Phi \in \omega_{\psi}, f\in \Sigma$ and $\varphi \in \Pi=\Theta^{ Aut}(\Sigma)$,   one has
	\begin{equation} \label{global: eq100}
		\langle A_{\Sigma}^{\mathrm{Aut}}(\Phi,f),\varphi \rangle_{\Pi}=\langle B_{\Pi}^{\mathrm{Aut}}(\Phi,\varphi),f \rangle_{\Sigma}
	\end{equation}
\end{lem}
\begin{proof}
	By an exchange of the order of integration, we have
	\begin{align*} 
		&\langle A_{\Sigma}^{\mathrm{Aut}}(\Phi,f),\varphi \rangle_{\Pi}\\
		=&\int_{Z_{W}(\A) GSp(W)(k)\backslash  GSp(W)(\A)} \int_{[ SO(V)]}\theta(\Phi)(h_1\cdot h,g)\cdot \overline{f(h_1\cdot h)}\ d h_1\ \overline{\varphi(g)}\ d g \notag \\
		=&  \int_{Z^{\bigtriangledown}(\A)R_0(k)\backslash R_0(\A)}\theta(\Phi)(h,g)\cdot \overline{f(h)}\cdot \overline{\varphi(g)} \ d h d g \\
		=&\int_{Z_{V}(\A) GSO(V)(k)\backslash  GSO(V)(\A)} \int_{[ Sp(W)]}\theta(\Phi)(h,g_1\cdot g)\cdot \overline{f(h)}\cdot \overline{\varphi(g_1\cdot g)}\  d g_1 d h \notag \\
		=&\langle B_{\Pi}^{\mathrm{Aut}}(\Phi,\varphi),f \rangle_{\Sigma} \notag
	\end{align*}
	
\end{proof}

\subsection{\bf Global periods}
Let $H \subset G$ be a subgroup so that $X=H\backslash G$ is quasi-affine where $G$ has no split center, and $\chi$ a unitary Hecke character of $H$. Then we may consider the global $(H,\chi)$-period:
$$
P_{H,\chi}: \mathcal{A}_0(G)\longrightarrow \mathbb{C}_{\chi}
$$
given by
$$
P_{H,\chi}(\phi)=\int_{[H]}\overline{\chi(h)}\phi(h)\ d h.
$$
For a cuspidal representation $\Sigma \subset \mathcal{A}_0(G)$, we write $P_{H,\chi,\Sigma}$ for the restriction of $P_{H,\chi}$ to $\Sigma$.

In this thesis, we will see two global periods: $P_{(N,\psi)\times \overline{T}_{k^{\prime}}}$ on $\mathcal{A}_{0}(PGSO(V))$ and $P_{U(2)}$ on $\mathcal{A}_{0}(PGSp(W))$.

We define a map
$$
i: T_{k^{\prime}} \longrightarrow GSO(V)
$$
given by
$$
i(t)=m\left(t,\det(t)\right) 
$$
for $t \in T_{k^{\prime}}$, where $m\left(t,\det(t)\right)$ is given in Subsection \ref{zf13:21a};
define a map
$$
j: T_{k^{\prime}} \longrightarrow GSp(W)
$$
given by
$$
j(t)=P(t,t)
$$
for $t \in T_{k^{\prime}}$, where $P(t,t)$ is given in Subsection \ref{subsection: 3.1a}.

Now we can define the global periods
$$
P_{(N,\psi)\times \overline{T}_{k^{\prime}}}: \mathcal{A}_{0}(PGSO(V))\longrightarrow \mathbb{C}_{\psi}
$$
by
$$
P_{(N,\psi)\times \overline{T}_{k^{\prime}}}(\phi)=\int_{[N]\times [\overline{T}_{k^{\prime}}]}\overline{\psi(n)}\cdot \phi(n\cdot i(t))\ d n d t;
$$
and
$$
P_{U(2)}: \mathcal{A}_{0}(PGSp(W))\longrightarrow \mathbb{C}
$$
by
$$
P_{U(2)}(\phi)=\int_{[U(2)]} \phi(\tau)\ d \tau=\int_{[\overline{T}_{k^{\prime}}]}\int_{[SU(2)]}\phi(s\cdot j(t)) \ d s d t.
$$

\subsection{\bf The maps $\beta_{\Sigma}^{X,\mathrm{Aut}}$}
In this subsection, we shall introduce the global analog of the map  $\beta_{\Sigma_v}^{X_v}$ introduced in the local setting. With $X_{\mathbb{A}}=H(\mathbb{A})\backslash G(\mathbb{A}),$ we have $G(\mathbb{A})$-equivariant map
$$
\theta: C_c^{\infty}(X_{\mathbb{A}},\chi) \longrightarrow C_{\mathrm{mod}}^{\infty}([G])
$$ 
defined by 
$$
\theta(f)(g)=\sum_{x\in X(k)}f(x\cdot g).
$$ 

We have a $G(\mathbb{A})$-equivariant map 
$$
\beta_{\Sigma}^{X,\mathrm{Aut}}: \Sigma \longrightarrow C^{\infty}(X_{\mathbb{A}},\chi)
$$
defined by
$$
\beta_{\Sigma}^{X,\mathrm{Aut}}(\phi)(g)=P_{H,\chi}(g\cdot \phi)
$$

Note that there is a natural pairing between $C_c^{\infty}(X_{\mathbb{A}},\chi)$ and $C^{\infty}(X_{\mathbb{A}},\chi)$ given by
$$
\langle f_1, f_2 \rangle_{X_{\mathbb{A}}}:=\int_{X_{\mathbb{A}}}f_1(x)\overline{f_2(x)}\frac{d g}{d h}.
$$

\subsection{\bf Global transfer of periods}

Let $\omega_{\psi}^{+}:=\otimes^{\prime}\omega_{\psi_v}^{+}\subset \omega_{\psi}$. We define the map
$$
q^{\prime}: \omega_{\psi}^{+}\longrightarrow S\left(\mu_2(\A)\cdot SU(2)(\A)\backslash Sp(W)(\A)\right)
$$
by 
$$
q^{\prime}(\Phi)(g)=\omega_{\psi}(g)\Phi(T_1)
$$ 
for $\Phi \in \omega_{\psi}^{+}$, where $\mu_2$ denotes the center of $Sp(W)$.
Then for $\Phi=\otimes_v \Phi_v \in \omega_{\psi}^{+}$, one has $q^{\prime}(\Phi)= \prod_v q_1(\Phi_v)|_{Sp(W)(\A)}$, where $q_1$ is defined in Subsection \ref{subsection: transfer}.
\begin{lem} \label{lemma7.2}
	For $\Phi \in \omega_{\psi}^{+}$, a cuspidal automorphic representation $\Pi$ of $PGSp(W)(\mathbb{A})$ and $f \in \Pi$, one has
	\begin{equation} \label{global: prop1}
		P_{(N,\psi)\times \overline{T}_{k^{\prime}}}(B_{\Pi}^{\mathrm{Aut}}(\Phi,f))=\langle q^{\prime}(\Phi),\beta_{\Pi}^{\mathrm{Y,\mathrm{Aut}}}(f) \rangle_{Y_1},
	\end{equation}
	where $Y_1\cong SU(2)(\A)\backslash Sp(W)(\A)$.
\end{lem}
\begin{proof}
	
	\begin{align*}
		&P_{(N,\psi)\times \overline{T}_{k^{\prime}}}(B_{\Pi}^{\mathrm{Aut}}(\Phi,f))\\
		=&\int_{[N]\times[\overline{T}_{k^{\prime}}]}\overline{\psi(n)}\int_{[Sp(W)]}\theta(\Phi)(n\cdot i(t^{-1}),g\cdot j(t))\cdot \overline{f(g \cdot j(t))}\ d^{Sp(W)} g d n d t \notag  \\
		=&\int_{[N]\times[\overline{T}_{k^{\prime}}]}\overline{\psi(n)}\int_{[Sp(W)]}\sum_{x \in V_2\otimes W(k)}\omega_{\psi}(n\cdot i(t^{-1}),g\cdot j(t))\Phi(x)\cdot \overline{f(g\cdot j(t))}\ d^{Sp(W)} g d n d t \notag  \\
	   =&\int_{[N]\times[\overline{T}_{k^{\prime}}]}\int_{[Sp(W)]}\sum_{x \in V_2\otimes W(k)}\psi(n\cdot Q(x)-n)\cdot \omega_{\psi}( i(t^{-1}),g\cdot j(t))\Phi(x)\cdot \overline{f(g\cdot j(t))}\ d^{Sp(W)} g d n d t \notag  \\
		= & \int_{[\overline{T}_{k^{\prime}}]}\int_{[Sp(W)]}    \sum_{x \in V_2\otimes W(k);Q(x)=1}\omega_{\psi}(i(t^{-1}),g\cdot j(t))\Phi(x)\cdot \overline{f(g\cdot j(t))}\ d^{Sp(W)} g d t  \notag  \\
		=&  \int_{[\overline{T}_{k^{\prime}}]}\int_{[Sp(W)]} \sum_{\gamma \in SU(2)(k) \backslash Sp(W) (k) }\omega_{\psi}(i(t^{-1}), \gamma \cdot g\cdot j(t))\Phi(T_1)\cdot \overline{f(g\cdot j(t))}\ d^{Sp(W)} g d t  \notag  \\
		=&\int_{[\overline{T}_{k^{\prime}}]}\int_{ SU(2) (k) \backslash Sp(W)(\mathbb{A})}    \omega_{\psi}(i(t^{-1}),g\cdot j(t))\Phi( T_1)\cdot \overline{f(g \cdot j(t))}\ d^{Sp(W)} g d t  \notag   \\
		=&\int_{[\overline{T}_{k^{\prime}}]}\int_{ SU(2)(\A) \backslash Sp(W)(\mathbb{A})}\int_{[SU(2)]}    \omega_{\psi}(i(t^{-1}),s\cdot g \cdot j(t))\Phi( T_1)\cdot \overline{f(s\cdot g \cdot j(t))}\ d s \frac{d^{Sp(W)}g}{d s}d t  \notag   \\
		=&\int_{[\overline{T}_{k^{\prime}}]}\int_{ SU(2)(\A) \backslash Sp(W)(\mathbb{A})}    \omega_{\psi}(i(t^{-1}), g \cdot j(t))\Phi( T_1)\cdot \int_{[SU(2)]}\overline{f(s\cdot g \cdot j(t))}\ d s \frac{d^{Sp(W)}g}{d s}d t  \notag  \qquad \text{($\ast$)} \\
		=&\int_{[\overline{T}_{k^{\prime}}]}\int_{ SU(2)(\A) \backslash Sp(W)(\mathbb{A})}    \omega_{\psi}(i(t^{-1}),  j(t)\cdot g)\Phi( T_1)\cdot \int_{[SU(2)]}\overline{f(s\cdot  j(t)\cdot g)}\ d s \frac{d^{Sp(W)}g}{d s}d t  \notag \qquad \text{($\ast\ast$)}  \\
		=&\int_{[\overline{T}_{k^{\prime}}]}\int_{ SU(2)(\A) \backslash Sp(W)(\mathbb{A})}    \omega_{\psi}(g)\Phi( T_1)\cdot \int_{[SU(2)]}\overline{f(s\cdot  j(t)\cdot g)}\ d s \frac{d^{Sp(W)}g}{d s} d t  \notag   \qquad \text{( $\ast\ast\ast$ )}\\
		=&\int_{ SU(2)(\A) \backslash Sp(W)(\mathbb{A})}    \omega_{\psi}(g)\Phi( T_1)\cdot \int_{[\overline{T}_{k^{\prime}}]}\int_{[SU(2)]}\overline{f(s\cdot  j(t)\cdot g)}\ d s d t \frac{d^{Sp(W)}g}{d s}  \notag   \\
		=&\int_{ SU(2)(\A) \backslash Sp(W)(\mathbb{A})}    \omega_{\psi}(   g)\Phi( T_1)\cdot \int_{[U(2)]} \overline{f(\tau\cdot g)}\ d \tau \frac{d^{Sp(W)}g}{d s} \notag   \\
		=&\langle q^{\prime}(\Phi),\beta_{\Pi}^{\mathrm{Y,\mathrm{Aut}}}(f) \rangle_{Y_1},
	\end{align*}
	where $(\ast)$ is because $\omega_{\psi}(1,s)\Phi(T_1)=\Phi(T_1)$, $(\ast\ast)$ is because we change $g$ by $j(t)\cdot g\cdot j(t^{-1})$, $(\ast\ast\ast)$ is because $\omega_{\psi}(i(t^{-1}),j(t))\Phi(T_1)=\Phi(T_1)$.

\end{proof}

Following from Lemma \ref{lemma7.2} and the fact that $q^{\prime}$ is surjective, we can deduce:
\begin{cor} \label{cor7.2}
	The global period $P_{U(2)} \neq 0$ on a cuspidal representation $\Pi$ of $PGSp(W)(\A)$ if and only if $\Pi=\Theta^{Aut}\left(\Sigma_1 \boxtimes \Sigma_2\right)$ with $\Sigma_1\ncong \Sigma_2$ and $\Sigma_2$ being globally  $\overline{T}_{k^{\prime}}-$distinguished, where $\Sigma_1 \boxtimes \Sigma_2$ is some cuspidal automorphic representation of $PGSO(V)\cong PGL_2\times PGL_2$.
\end{cor}

\subsection{\bf Factorization of global periods}\label{subsection: factorization}

For every irreducible cuspidal automorphic representation $\Sigma \subset \mathcal{A}_{0}\left(GSO(V)\right)$ with trivial central character, we fix an isomorphism
$$
\Sigma\cong \bigotimes_{v}^{\prime}\Sigma_{v}
$$ 
with respect to spherical vectors $\{v_{0,v}  \}$ such that $\langle v_{0,v},v_{0,v} \rangle_{\Sigma_v}=1$ for almost all $v$, where  $\langle \cdot,\cdot \rangle_{\Sigma_v}$ is the fixed inner product on each $\Sigma_v$ such that the Petersson inner product can decompose as:
$$
\langle \cdot,\cdot   \rangle_{\Sigma}=\prod_{v}\langle \cdot,\cdot \rangle_{\Sigma_v}.
$$

Let the global Weil representation 
$$
\omega_{\psi}:=S(V_2\otimes W):=\left(\bigotimes_{v|\infty} S(V_2\otimes W)_{v}\right)\bigotimes \left(\bigotimes_{v<\infty}C_c^{\infty}(V_2\otimes W)_{v} \right)
$$
with respect to the family of spherical vectors $\{\Phi_{0,v,1}\mid \text{$v$ is non-Archimedean and $v\nmid 2$} \}$ ( see Subsection \ref{subsection: measure: weil}). We equip the inner product defined by
$$
\langle \cdot,\cdot \rangle_{\omega_{\psi}}:=\prod_{v}\langle \cdot,\cdot \rangle_{\omega_{\psi_v}}.
$$

We have equipped an inner product $\langle \cdot,\cdot \rangle_{\theta(\Sigma_v)}$ on $\theta(\Sigma_v)$ via the local doubling theta integral ( see Subsection \ref{subsection2.3c} ). We form a global representation by
$$
\Theta^{\mathbb{A}}(\Sigma):=\bigotimes_{v}^{\prime}\theta(\Sigma_v)
$$
with respect to a family of spherical vectors $\{w_{0,v}\}$ such that $\langle w_{0,v},w_{0,v} \rangle_{\theta(\Sigma_v)}=1$ for almost all $v$. The abstract global theta lift inherits an inner product 
$$
\langle \cdot,\cdot  \rangle_{\Theta^{\mathbb{A}}(\Sigma)}:=\prod_{v}\langle  \cdot ,\cdot\rangle_{\theta(\Sigma_{v})}.
$$
We fix an isometric isomorphism
$$
\Theta^{\mathrm{Aut}}(\Sigma)\cong \Theta^{\mathbb{A}}(\Sigma)
$$

Let $ \Sigma=\Sigma_1\boxtimes\Sigma_2 \in \mathcal{A}_0\left(PGSO(V)\right)$ be an irreducible cuspidal automorphic representation of $PGSO(V)(\A)$ with $\Sigma_1\ncong \Sigma_2$ and $\Sigma_2$ being  $\overline{T}_{k^{\prime}}$-distinguished. Let $\Sigma_{1,2}=\Sigma_1\boxtimes\Sigma_2+\Sigma_2\boxtimes\Sigma_1$.

We define the local maps which we will need later. Define
$$
A^{\prime}_{\Sigma_{v}}: \omega_{\psi_v}^{+}\otimes \Sigma_v^{\vee} \longrightarrow \theta(\Sigma_v)
$$
by
$$
A^{\prime}_{\Sigma_v}(\Phi_v, \phi_v)=2\cdot A_{\Sigma_v}(\Phi_v,\phi_v).
$$
Here, we use the factor $2$, because 
\begin{align*}
	\langle A_{\Sigma_v}^{\prime}(\Phi_{0,1,v},v_{0,v}), A_{\Sigma_v}^{\prime}(\Phi_{0,1,v},v_{0,v})\rangle_{\theta(\Sigma_v)}=&\langle 2\cdot A_{\Sigma_v}(\Phi_{0,1,v},v_{0,v}), 2\cdot A_{\Sigma_v}(\Phi_{0,1,v},v_{0,v})\rangle_{\theta(\Sigma_v)} \\
	=& 4\cdot Z_{\Sigma_v}^{PGSO(V)}(\Phi_{0,1,v},\Phi_{0,1,v},v_{0,v},v_{0,v}) \\
	=& Z_{\Sigma_v}^{PGSO(V)}(\Phi_{0,v},\Phi_{0,v},v_{0,v},v_{0,v}), \qquad \text{( by Lemma \ref{lem6.8})}
\end{align*}
and we know the value of $Z_{\Sigma_v}^{PGSO(V)}(\Phi_{0,v},\Phi_{0,v},v_{0,v},v_{0,v})$ in Lemma \ref{zf98: lemma6.}.
Similarly, define the map
$$
B^{\prime}_{\Sigma_v}: \omega_{\psi_v}^{+}\otimes \theta(\Sigma_v)^{\vee} \longrightarrow \Sigma_v
$$
by
$$
B^{\prime}_{\Sigma_v}(\Phi_v,f_v)=2\cdot B_{\Sigma_v}(\Phi_v,f_v).
$$
By the equation (\ref{local: eq300}), for $\Phi_v \in \omega_{\psi_v}^{+}$, $\phi_v \in \Sigma_v$ and $f_v \in \theta(\Sigma_v)$, one has
\begin{align}
	&\langle A_{\Sigma_v}^{\prime}(\Phi_v,\phi_v),f_v  \rangle_{\theta(\Sigma_v)} \notag \\
	=& 2\cdot \langle A_{\Sigma_v}(\Phi_v,\phi_v),f_v  \rangle_{\theta(\Sigma_v)} \notag \\
	=&2\cdot \langle B_{\Sigma_v}(\Phi_v,f_v),\phi_v  \rangle_{\Sigma_v} \notag \\
	=& \langle B_{\Sigma_v}^{\prime}(\Phi_v,f_v),\phi_v  \rangle_{\Sigma_v}. \label{local: eq301}
\end{align}
We give five examples to illustrate how to define the regularized product.
\begin{enumerate}
	\item Let 
	$$
	\ell_{\Sigma_v}^{X_v,\#}:=\lambda_v\cdot \frac{\left|L(1,\Sigma_{1,v},\mathrm{Ad})\right|^{1/2}\cdot \left|L(1,\Sigma_{2,v},\mathrm{Ad})\right|^{1/2}\cdot \left|L(1,\chi_{k^{\prime}_{v}\slash k_v},\mathrm{std})\right|^{1/2}}{\zeta_{k_v}(2)\cdot \left|L(\frac{1}{2},\mathrm{BC}\left(\Sigma_{2,v}\right),\mathrm{std})\right|^{1/2}}\cdot\ell_{\Sigma_v}^{X_v},
	$$
	where $\lambda_v$ satisfies $\left| \lambda_v \right|=1$ and
	$$
	\ell_{\Sigma_v}^{X_v,\#}(v_{0,v})=1
	$$ 
	for almost all $v$. We have $\ell_{\Sigma_v}^{X_v,\#} \in Hom_{N\times \overline{T}_{k^{\prime}_v}}\left(\Sigma_v,\psi_v\boxtimes\mathbb{C}\right)$. Let the normalized adelic product
	\begin{align*}
		\ell_{\Sigma}^{X,\A}:=&\prod_{v}^{\ast}\ell_{\Sigma_v}^{X_v}:= \frac{\zeta_{k}(2)\cdot \left|L(\frac{1}{2},\mathrm{BC}\left(\Sigma_2\right),\mathrm{std})\right|^{1/2}}{\left|L(1,\Sigma_1,\mathrm{Ad})\right|^{1/2}\cdot \left|L(1,\Sigma_2,\mathrm{Ad})\right|^{1/2}\cdot \left|L(1,\chi_{k^{\prime}\slash k},\mathrm{std})\right|^{1/2}} \cdot \\
		& \prod_{v} \ell_{\Sigma_v}^{X_v,\#}, 
	\end{align*}
	where $\mathrm{BC}(\Sigma_2)$ denotes the base change of $\Sigma_2$ to $GL_2(\mathbb{A}_{k^{\prime}})$, $\chi_{k^{\prime}/k}$ is determined by $k^{\prime}$ via class field theory.
	\item 
	Let 
	\begin{equation*}
		A_{\Sigma_v}^{\#}:=\lambda_v\cdot  \frac{\zeta_{k_v}(2)^{1/2}\cdot \zeta_{k_v}(4)^{1/2}}{\left|L(1,\Sigma_{1,v}\times \Sigma_{2,v})\right|^{1/2}}\cdot A^{\prime}_{\Sigma_v}
	\end{equation*}
	where $\lambda_v$ satisfies $\left| \lambda_v \right|=1$ and
	$$
	A_{\Sigma_v}^{\#}(\Phi_{0,v,1},v_{0,v})=w_{0,v}
	$$ 
	for almost all $v$ ( by Lemma \ref{lem6.8} ). 
	Let the normalized adelic product
	\begin{align*}
		A_{\omega_{\psi}^{+},\Sigma}^{\A}:=\prod_{v}^{\ast}A^{\prime}_{\Sigma_v}=& \frac{\left|L(1,\Sigma_1\times \Sigma_2)\right|^{1/2}}{\zeta_{k}(2)^{1/2}\cdot \zeta_{k}(4)^{1/2}} \cdot \prod_{v} A_{\Sigma_v}^{\#}.
	\end{align*}
	\item 
	Let 
	\begin{align*}
		B_{\Sigma_v}^{\#}:=&\lambda_v\cdot   \frac{\zeta_{k_v}(2)^{1/2}\cdot \zeta_{k_v}(4)^{1/2}}{\left|L(1,\Sigma_{1,v}\times \Sigma_{2,v})\right|^{1/2}}\cdot B^{\prime}_{\Sigma_v}, 
	\end{align*}
	where $\lambda_v$ satisfies $\left| \lambda_v \right|=1$ and
	$$
	B_{\Sigma_v}^{\#}(\Phi_{0,v,1},w_{0,v})=v_{0,v}
	$$ 
	for almost all $v$ ( by Lemma \ref{lem6.8} ). Let the normalized adelic product
	\begin{align*}
		B_{\omega_{\psi}^{+}, \Sigma}^{\A}:=\prod_{v}^{\ast}B^{\prime}_{\Sigma_v}=& \frac{\left|L(1,\Sigma_1\times \Sigma_2)\right|^{1/2}}{\zeta_{k}(2)^{1/2}\cdot \zeta_{k}(4)^{1/2}} \cdot \prod_{v} B_{\Sigma_v}^{\#}.
	\end{align*}

	\item 
	Let 
	\begin{align*}
		\ell_{\Sigma_v}^{Y_v,\#}:=&\lambda_v\cdot   \frac{1}{\zeta_{k_v}(2)^{1/2}\cdot \zeta_{k_v}(4)^{1/2}}\cdot \\
		&\frac{\left|L(1,\Sigma_{1,v},\mathrm{Ad})\right|^{1/2}\cdot \left|L(1,\Sigma_{2,v},\mathrm{Ad})\right|^{1/2}\cdot \left|L(1,\chi_{k^{\prime}_v\slash k_v},\mathrm{std})\right|^{1/2}}{\left|L(1,\Sigma_{1,v}\times \Sigma_{2,v})\right|^{1/2}\cdot \left|L(1/2,\mathrm{BC}(\Sigma_{2,v}),\mathrm{std})\right|^{1/2}}\cdot \ell_{\Sigma_v}^{Y_v},
	\end{align*}
	where $\lambda_v$ satisfies $\left| \lambda_v \right|=1$ and
	$$
	\ell_{\Sigma_v}^{Y_v,\#}(w_{0,v})=1
	$$ 
	for almost all $v$. We have $\ell_{\Sigma_v}^{Y_v,\#} \in Hom_{U(2)(k_v)}\left(\theta(\Sigma_v),\mathbb{C} \right)$.   
	Let the normalized adelic product
	\begin{align*}
		&\ell_{\Sigma}^{Y,\A}:=\prod_{v}^{\ast}\ell_{\Sigma_v}^{Y_v}\\
		&:= \zeta_{k}(2)^{1/2}\cdot \zeta_{k}(4)^{1/2}\cdot\frac{\left|L(1,\Sigma_1\times\Sigma_2)\right|^{1/2}\cdot \left|L(1/2,\mathrm{BC}(\Sigma_{2}),\mathrm{std})\right|^{1/2}}{\left|L(1,\Sigma_1,\mathrm{Ad})\right|^{1/2}\cdot \left|L(1,\Sigma_2,\mathrm{Ad})\right|^{1/2}\cdot \left|L(1,\chi_{k^{\prime}\slash k},\mathrm{std})\right|^{1/2}} \cdot \prod_{v}	\ell_{\Sigma_v}^{Y_v,\#}.
		\end{align*}

\item Let 
\begin{align*}
	Z_{\Sigma_v}^{SO(V),\#}=&\lambda_v \cdot 	\frac{\zeta_{k_v}(2)\cdot \zeta_{k_v}(4)}{L(1,\Sigma_{1,v}\times \Sigma_{2,v})}\cdot Z_{\Sigma_v}^{SO(V)}, 
	\end{align*}
where $\lambda_v$ satisfies $\left|\lambda_v\right|=1$ and
$$
Z_{\Sigma_v}^{SO(V),\#}\left(\Phi_{0,1},\Phi_{0,1},v_{0,v},v_{0,v}\right)=1
$$
for almost all $v$. Let the normalized adelic product 
\begin{equation*}
	\prod_{v}^{\ast}Z_{\Sigma_v}^{SO(V)}= \frac{L(1,\Sigma_1\times\Sigma_2)}{\zeta_{k}(2)\cdot \zeta_{k}(4)}\cdot \prod_v Z_{\Sigma_v}^{SO(V),\#}.
\end{equation*}

\end{enumerate}

There are constants $c(\Sigma) \in \mathbb{C}^{\times}$ such that
$$
P_{(N,\psi)\times \overline{T}_{k^{\prime}},\Sigma}=c(\Sigma)\cdot\ell_{\Sigma}^{X,\mathbb{A}}  
$$
and
$$
\beta_{\Sigma}^{X,\mathrm{Aut}}=c(\Sigma)\cdot \beta_{\Sigma_{12}}^{X,\mathbb{A}}.  
$$

Similarly, we have $a(\Sigma),b(\Sigma) \in \mathbb{C}^{\times}$ such that
$$
A_{\Sigma}^{\mathrm{Aut}}|_{\omega_{\psi}^{+}}=a(\Sigma)\cdot A_{\omega_{\psi}^{+}, \Sigma}^{\mathbb{A}}   
$$
and
$$
B_{\Sigma}^{\mathrm{Aut}}|_{\omega_{\psi}^{+}}=b(\Sigma)\cdot B_{\omega_{\psi}^{+}, \Sigma}^{\mathbb{A}}.  
$$

We will show the following lemma which we will need later.
\begin{lem} \label{lastchance}
For $ \Phi_v\in \omega_{\psi_v}^{+}$ and $ f_v \in \theta(\Sigma_v)$, we have the following equality
$$
\langle 2\cdot q(\Phi_v), \beta_{\Sigma_v}^{Y_v}(f_v)\rangle_{Y_v}=\langle q_1(\Phi_v), \beta_{\Sigma_v}^{Y_v}(f_v)\rangle_{Y_{1,v}},
$$
where $Y=U(2)(k_v)\backslash PGSp(W)(k_v)$, $Y_1=SU(2)(k_v)\backslash Sp(W)(k_v)$, and $q, q_1$ are defined in Subsection \ref{subsection: transfer}.
\end{lem}
\begin{proof}
\begin{align*}
	&\langle 2\cdot q(\Phi_v), \beta_{\Sigma_v}^{Y_v}(f_v)\rangle_{Y_v}\\
	=&\int_{U(2)(k_v)\backslash PGSp(W)(k_v)}2\cdot \int_{\overline{T}_{k^{\prime}_v}} \mathrm{pr}_1(\Omega(1)\left(P(t_v,t_v)\cdot g_v \right)\Phi_v)(T_1)\cdot \overline{\beta_{\Sigma_v}^{Y_v}(f_v)(g_v)} \ d t_v \frac{d^{PGSp(W)} g_v}{d t_v\cdot d s_v} \\
	=&\int_{U(2)(k_v)\backslash PGSp(W)(k_v)}2\cdot \int_{\overline{T}_{k^{\prime}_v}} \mathrm{pr}_1(\Omega(1)\left(P(t_v,t_v)\cdot g_v \right)\Phi_v)(T_1)\cdot \overline{\beta_{\Sigma_v}^{Y_v}(f_v)\left(P(t_v,t_v)\cdot g_v\right)} \ d t_v \frac{d^{PGSp(W)} g_v}{d t_v\cdot d s_v} \\
	=&2\cdot \int_{SU(2)(k_v)\backslash PGSp(W)(k_v)}\mathrm{pr}_1(\Omega(1)\left( g_v \right)\Phi_v)(T_1)\cdot \overline{\beta_{\Sigma_v}^{Y_v}(f_v)\left( g_v\right)} \ \frac{d^{PGSp(W)} g_v}{d s_{v}} \\
	=&2\cdot \int_{\{\pm 1\}\cdot SU(2)(k_v)\backslash Sp(W)(k_v)}\mathrm{pr}_1(\Omega(1)\left( g_v \right)\Phi_v)(T_1)\cdot \overline{\beta_{\Sigma_v}^{Y_v}(f_v)\left( g_v\right)}\ \frac{d^{PGSp(W)} g_v}{d s_{v}} \\
	=&2\cdot \int_{\{\pm 1\}\cdot SU(2)(k_v)\backslash Sp(W)(k_v)}\mathrm{pr}_1(\Omega(1)\left( g_v \right)\Phi_v)(T_1)\cdot \overline{\beta_{\Sigma_v}^{Y_v}(f_v)\left( g_v\right)}\ \frac{d^{Sp(W)} g_v}{\nu\cdot d s_{v}} \\
	=&\int_{ SU(2)(k_v)\backslash Sp(W)(k_v)}\Omega(1)\left( g_v \right)\Phi_v(T_1)\cdot \overline{\beta_{\Sigma_v}^{Y_v}(f_v)\left( g_v\right)} \ \frac{d^{Sp(W)} g_v}{d s_{v}} \\
	=&\langle q_1(\Phi_v), \beta_{\Sigma_v}^{Y_v}(f_v)\rangle_{Y_{1,v}},
\end{align*}
where $\nu$ is the counting measure of $\{\pm 1 \}$.
\end{proof}

\begin{prop}
Consider a cuspidal representation $\Pi$ of $PGSp(W)(\mathbb{A})$ which is the theta lift of $\Sigma=\Sigma_1\boxtimes\Sigma_2+\Sigma_2\boxtimes\Sigma_1$, where $\Sigma_1\ncong \Sigma_2$.
For $f \in \Pi$, we have
$$
P_{\mathrm{U}(2),\Pi}(f)=\overline{c(\Sigma_{1}\boxtimes \Sigma_2)\cdot b(\Sigma_1\boxtimes\Sigma_2)}\cdot \ell_{\Sigma_{1}\boxtimes\Sigma_2}^{Y,\mathbb{A}}(f)+ \overline{c(\Sigma_{2}\boxtimes \Sigma_1)\cdot b(\Sigma_2\boxtimes\Sigma_1)}\cdot \ell_{\Sigma_{2}\boxtimes\Sigma_1}^{Y,\mathbb{A}}(f) .
$$
\end{prop}
\begin{proof}

For every $\Phi \in \omega_{\psi}^{+}$ and $f \in \Pi$, one has
\begin{align*}
	&\langle q^{\prime}(\Phi),\beta_{\Pi}^{Y,\mathrm{Aut}}(f)\rangle_{Y_1} \notag \\
	=& P_{(N,\psi)\times \overline{T}_{k^{\prime}},\Sigma_{1,2}}\left(B_{\Pi}^{\mathrm{Aut}}(\Phi,f)\right)  \notag  \qquad \text{( by Lemma \ref{lemma7.2})}\\
	=&P_{(N,\psi)\times \overline{T}_{k^{\prime}},\Sigma_{1}\boxtimes\Sigma_2}\left(B_{\Sigma_1\boxtimes\Sigma_2}^{\mathrm{Aut}}(\Phi,f)\right)+P_{(N,\psi)\times \overline{T}_{k^{\prime}},\Sigma_{2}\boxtimes\Sigma_1}\left(B_{\Sigma_2\boxtimes\Sigma_1}^{\mathrm{Aut}}(\Phi,f)\right) \notag \\
	&\text{( $\Sigma_{1,2}=\Sigma_1\boxtimes\Sigma_2+ \Sigma_2\boxtimes\Sigma_1$ )} \\
	=&  c(\Sigma_{1}\boxtimes\Sigma_2)b(\Sigma_1\boxtimes\Sigma_2)\ell_{\Sigma_{1}\boxtimes\Sigma_2}^{X,\mathbb{A}}\left(B^{\A}_{\omega_{\psi}^{+},\Sigma_1\boxtimes\Sigma_2}(\Phi_v,f_v)\right)+ \\ &c(\Sigma_{2}\boxtimes\Sigma_1)b(\Sigma_2\boxtimes\Sigma_1)\ell_{\Sigma_{2}\boxtimes\Sigma_1}^{X,\mathbb{A}}\left(B^{\A}_{\omega_{\psi}^{+},\Sigma_2\boxtimes\Sigma_1}(\Phi_v,f_v)\right)\notag \\
	=& c(\Sigma_1\boxtimes\Sigma_2)b(\Sigma_1\boxtimes\Sigma_2) \prod_v^{\ast}\langle 2\cdot q(\Phi_v), \beta_{\Sigma_{1,v}\boxtimes\Sigma_{2,v}}^{Y_v}(f_v)\rangle_{Y_v}+ \\
	&c(\Sigma_2\boxtimes\Sigma_1)b(\Sigma_2\boxtimes\Sigma_1) \prod_v^{\ast}\langle 2\cdot q(\Phi_v), \beta_{\Sigma_{2,v}\boxtimes\Sigma_{1,v}}^{Y_v}(f_v)\rangle_{Y_v}  \notag  \qquad\text{( by Corollary \ref{cor4.24})}\\
	=& c(\Sigma_1\boxtimes\Sigma_2)b(\Sigma_1\boxtimes\Sigma_2) \prod_v^{\ast}\langle q_1(\Phi_v), \beta_{\Sigma_{1,v}\boxtimes\Sigma_{2,v}}^{Y_v}(f_v)\rangle_{Y_{1,v}}+ \\
	&c(\Sigma_2\boxtimes\Sigma_1)b(\Sigma_2\boxtimes\Sigma_1) \prod_v^{\ast}\langle q_1(\Phi_v), \beta_{\Sigma_{2,v}\boxtimes\Sigma_{1,v}}^{Y_v}(f_v)\rangle_{Y_{1,v}}  \notag \qquad \text{( by Lemma \ref{lastchance})} \\
	= &\langle  q^{\prime}(\Phi),\overline{c(\Sigma_1\boxtimes\Sigma_2)b(\Sigma_1\boxtimes\Sigma_2)}\beta_{\Sigma_1\boxtimes\Sigma_2}^{Y,\mathbb{A}}(f)+\overline{c(\Sigma_2\boxtimes\Sigma_1)b(\Sigma_2\boxtimes\Sigma_1)} \beta_{\Sigma_2\boxtimes\Sigma_1}^{Y,\mathbb{A}}(f) \rangle_{Y_1(\mathbb{A})}. \notag
\end{align*}

Since the linear span of $\{q^{\prime}(\Phi)\mid \Phi \in \omega^{+}_{\psi}\}$ is $S\left(\mu_2(\A)\cdot SU(2)(\A)\backslash Sp(W)(\mathbb{A})\right)$, and $\beta_{\Sigma_1\boxtimes\Sigma_2}^{Y,\mathbb{A}}, \beta_{\Sigma_2\boxtimes\Sigma_1}^{Y,\mathbb{A}}$, and $ \beta_{\Pi}^{Y,\mathrm{Aut}}$ are in $S\left(\mu_2(\A)\cdot SU(2)(\A)\backslash Sp(W)(\A)\right)$,  we have
$$
\beta_{\Pi}^{Y,\mathrm{Aut}}(f)=\overline{c(\Sigma_1\boxtimes\Sigma_2)b(\Sigma_1\boxtimes\Sigma_2)}\beta_{\Sigma_1\boxtimes\Sigma_2}^{Y,\mathbb{A}}(f)+ \overline{c(\Sigma_2\boxtimes\Sigma_1)b(\Sigma_2\boxtimes\Sigma_1)}\beta_{\Sigma_2\boxtimes\Sigma_1}^{Y,\mathbb{A}}(f).
$$
By evaluating at $T_1$ on both side, we obtain
$$
P_{\mathrm{U}(2),\Pi}(f)=\overline{c(\Sigma_1\boxtimes\Sigma_2)\cdot b(\Sigma_1\boxtimes\Sigma_2)}\cdot \ell_{\Sigma_1\boxtimes\Sigma_2}^{Y,\mathbb{A}}(f)+\overline{c(\Sigma_2\boxtimes\Sigma_1)\cdot b(\Sigma_2\boxtimes\Sigma_1)}\cdot \ell_{\Sigma_2\boxtimes\Sigma_1}^{Y,\mathbb{A}}(f).
$$
\end{proof}

By \cite[Lemma 7.11]{GI}, one has the Rallis inner product formula
\begin{lem} \label{Rallisinnerproduct}
For any $\Phi=\otimes^{\prime}_v\Phi_v\in \omega_{\psi}$ and $f=\otimes^{\prime}_{v}f_v\in \Sigma$, one has
\begin{equation}\label{zf13: Rallisinnerproduct}
	\langle A^{Aut}_{\Sigma}(\Phi,f),A^{Aut}_{\Sigma}(\Phi,f) \rangle_{\Pi}=\prod^{\ast}_{v}Z_{\Sigma_v}^{SO(V)}\left(\Phi_v,\Phi_v,f_v,f_v\right).
\end{equation}
\end{lem}
We also have a local version of Lemma \ref{Rallisinnerproduct}.
\begin{lem}
For any $\Phi_v\in \omega_{\psi_v}^{+}$ and $f_v\in \Sigma_v$, one has
\begin{equation}\label{zf13: Rallisinnerproductlocal}
	\langle A^{\prime}_{\Sigma_v}(\Phi_v,f_v),A^{\prime}_{\Sigma_v}(\Phi_v,f_v) \rangle_{\theta(\Sigma_v)}=Z_{\Sigma_v}^{SO(V)}\left(\Phi_v,\Phi_v,f_v,f_v\right).
\end{equation}
\end{lem}
\begin{proof}
\begin{align*}
	&\langle A^{\prime}_{\Sigma_v}(\Phi_v,f_v),A^{\prime}_{\Sigma_v}(\Phi_v,f_v) \rangle_{\theta(\Sigma_v)}  \\
	=&\langle 2\cdot A_{\Sigma_v}(\Phi_v,f_v),2\cdot A_{\Sigma_v}(\Phi_v,f_v) \rangle_{\theta(\Sigma_v)}  \\
	=&4\cdot Z_{\Sigma_v}^{PGSO(V)}(\Phi_v,\Phi_v,f_v,f_v) \qquad \text{( by equation (\ref{local: eq2}))}\\
	=&4\cdot \int_{PGSO(V)(k_v)}\langle \Omega(1)(h_v)\Phi_v,\Phi_v \rangle_{\Omega(1)}\cdot \overline{\langle \Sigma_v(h_v)f_v,f_v \rangle}_{\Sigma_v}\ d^{PGSO(V)} h_v \\
	=& 4 \cdot \int_{\{\pm 1 \}\backslash SO(V)(k_v) }\langle \Omega(1)(h_v)\Phi_v,\Phi_v \rangle_{\Omega(1)}\cdot \overline{\langle \Sigma_v(h_v)f_v,f_v \rangle}_{\Sigma_v}\ d^{PGSO(V)} h_v \qquad \text{($\ast$)}\\
	=& 4 \cdot \int_{\{\pm 1 \}\backslash SO(V)(k_v) }\langle \Omega(1)(h_v)\Phi_v,\Phi_v \rangle_{\Omega(1)}\cdot \overline{\langle \Sigma_v(h_v)f_v,f_v \rangle}_{\Sigma_v}\ \frac{d^{SO(V)} h_v}{\nu} \qquad \text{( by equation (\ref{compatmeasure1}))}\\
	=& 2 \cdot \int_{ SO(V)(k_v) }\langle \Omega(1)(h_v)\Phi_v,\Phi_v \rangle_{\Omega(1)}\cdot \overline{\langle \Sigma_v(h_v)f_v,f_v \rangle}_{\Sigma_v}\ d^{SO(V)} h_v \qquad \text{($\ast\ast$)}\\
	=&\int_{ SO(V)(k_v) }\langle \omega_{\psi_v}(h_v)\Phi_v,\Phi_v \rangle_{\omega_{\psi_v}}\cdot \overline{\langle \Sigma_v(h_v)f_v,f_v \rangle}_{\Sigma_v}\ d^{SO(V)} h_v \qquad \text{( by equation (\ref{innerproductomega1}))}\\
	=&Z_{\Sigma_v}^{SO(V)}\left(\Phi_v,\Phi_v,f_v,f_v\right) \qquad \text{( by equation (\ref{zeta2}))},
\end{align*}
where $(\ast)$ is because the support of $\langle \Omega(1)(h_v)\Phi_v,\Phi_v  \rangle_{\Omega(1)}$ lies in $\{\pm 1 \}\backslash SO(V)(k_v) \subset PGSO(V)(k_v)$; $(\ast\ast)$ is by the definition of counting measure $\nu$.
\end{proof}
\begin{prop}
We have $a(\Sigma)=b(\Sigma)$, $\left|a(\Sigma)\right|=1$, and $\left|c(\Sigma) \right|=\frac{1}{2}$.  
\end{prop}
\begin{proof}
Comparing the global equation (\ref{global: eq100}) with the local equation (\ref{local: eq301}), we have
$$
a(\Sigma)=b(\Sigma).
$$
Comparing the Rallis inner product formula (\ref{zf13: Rallisinnerproduct}) with its local analog of equation (\ref{zf13: Rallisinnerproductlocal}), we have
$$
\left|a(\Sigma)\right|=1
$$
The value $\left|c(\Sigma)\right|=\frac{1}{2}$ can be determined by \cite[Corollary 4.2]{LM} and \cite[Proposition 7]{W2}.
\end{proof} 

As a consequence, we have the main global result of this thesis. 
\begin{thm} \label{coefficient: prop}


For any irreducible cuspidal automorphic representation $\Sigma=\Sigma_1\boxtimes\Sigma_2$ of $PGSO(V)(\A)\cong PGL_2(\A)\times PGL_2(\A)$ with $\Sigma_1\ncong \Sigma_2$ and $f \in \Pi=\Theta^{Aut}(\Sigma)$,
$$
P_{\mathrm{U}(2),\Pi}(f)=d_1  \cdot \ell_{\Sigma_{1}\boxtimes\Sigma_2}^{Y,\mathbb{A}}(f)+d_2 \cdot \ell_{\Sigma_{2}\boxtimes\Sigma_1}^{Y,\mathbb{A}}(f),
$$
where $\left|d_1\right|=\left|d_2\right|=\frac{1}{2}$.

\end{thm}





\end{document}